\documentclass[a4paper,UKenglish,cleveref, autoref, thm-restate]{lipics-v2021}

\usepackage{graphicx} 
\usepackage{cleveref}
\usepackage[dvipsnames]{xcolor}
\usepackage{tikz}
\usetikzlibrary{decorations.pathreplacing,angles,quotes}
\usepackage{booktabs}

\DeclareMathOperator{\dist}{dist}
\DeclareMathOperator{\td}{td}

\newcommand{\ext}{{\it ext}}

\newcommand{\NP}{{\sf NP}}
\newcommand{\fvs}{{\sc Feedback Vertex Set}}
\newcommand{\col}{{\sc Colouring}}
\newcommand{\stcut}{{\sc Stable Cut}}

\theoremstyle{claimstyle}
\newtheorem{myclaim}{Claim}[theorem]

\usetikzlibrary{decorations.pathreplacing, angles,quotes, calligraphy}
\usetikzlibrary{decorations.pathmorphing}
\usetikzlibrary{fit,positioning}
\tikzstyle{vertex}=[thin,circle,inner sep=0.cm, minimum size=1.7mm, fill=black, draw=black]
 \tikzstyle{edge}=[thick, draw = gray]
 \tikzstyle{br} = [decorate, ultra thick, decoration = {calligraphic brace}]

\title{Colouring Graphs Without a Subdivided H-Graph: A Full Complexity Classification}
\titlerunning{Colouring Graphs Without a Subdivided H-Graph: A Full Complexity Classification}

\author{Tala Eagling-Vose}{Department of Computer Science, Durham University, Durham, UK}{tala.j.eagling-vose@durham.ac.uk}{https://orcid.org/0009-0008-0346-7032}{}
\author{Jorik Jooken}{Department of Computer Science, KU Leuven Campus Kulak-Kortrijk, 8500 Kortrijk, Belgium}{jorik.jooken@kuleuven.be}{https://orcid.org/0000-0002-5256-1921}{supported by a Postdoctoral Fellowship of the FWO (1222524N).}
\author{Felicia Lucke}{ENS Lyon, France}{felicia.lucke@ens-lyon.fr}{https://orcid.org/0000-0002-9860-2928}{supported by EPSRC (EP/X01357X/1) and SNSF Postdoc Mobility Grant 230578.}
\author{Barnaby Martin}{Department of Computer Science, Durham University, Durham, UK}{barnaby.d.martin@durham.ac.uk}{https://orcid.org/0000-0002-4642-8614}{}
\author{Dani\"el Paulusma}{Department of Computer Science, Durham University, Durham, UK}{daniel.paulusma@durham.ac.uk}{https://orcid.org/0000-0001-5945-9287}{supported by Leverhulme Trust (RPG-2024-182) and EPSRC (EP/X01357X/1).}

\authorrunning{T. Eagling-Vose, J. Jooken, F. Lucke, B. Martin and D. Paulusma}

\Copyright{Tala Eagling-Vosa, Jorik Jooken, Felicia Lucke, Barny Martin and Dani\"el Paulusma}


\keywords{colouring, \and forbidden subgraph, \and complexity dichotomy\\[-20pt]}

\category{}
\relatedversion{}
\acknowledgements{}

\nolinenumbers
\hideLIPIcs  

\EventEditors{John Q. Open and Joan R. Access}
\EventNoEds{2}
\EventLongTitle{42nd Conference on Very Important Topics (CVIT 2016)}
\EventShortTitle{CVIT 2016}
\EventAcronym{CVIT}
\EventYear{2016}
\EventDate{December 24--27, 2016}
\EventLocation{Little Whinging, United Kingdom}
\EventLogo{}
\SeriesVolume{42}
\ArticleNo{23}

\begin{document}
\maketitle

\begin{abstract}
We consider {\sc Colouring} on graphs that are $H$-subgraph-free for some fixed graph $H$, which are graphs that do not contain $H$ as a subgraph. To classify the complexity of {\sc Colouring} on $H$-subgraph-free graphs for connected $H$, it remains to consider when $H$ is a tree of maximum degree~$4$ with exactly one vertex of degree~$4$, or a tree of maximum degree~$3$ with at least two vertices of degree~$3$. We let $H$  be a so-called subdivided ``H''-graph, which is either a subdivided $\mathbb{H}_0$: a tree of maximum degree~$4$ that is a star, or a subdivided  $\mathbb{H}_1$:  a tree of maximum degree~$3$ with exactly two vertices of degree~$3$.  We develop new decomposition theorems resulting in polynomial-time algorithms, and in combination with known results, fully classify all cases $\mathbb{H}_0$ and $\mathbb{H}_1$.  To illustrate the wider applicability of our techniques, we also employ them to obtain similar new polynomial-time results for two other classic graph problems: {\sc Stable Cut} and, in part, {\sc Feedback Vertex Set}.
\end{abstract}

\section{Introduction}\label{sec:intro}

For an integer~$k\geq 1$, a {\it $k$-colouring} $c$ of a graph $G=(V,E)$ maps every vertex $u\in V$ to an integer $c(u)\in \{1,2,\ldots,k\}$ such that $c(v)\neq c(w)$ for every two vertices $v$ and $w$ with $vw\in E$. The corresponding decision problem {\sc Colouring} has as input a graph $G$ and integer $k$ and is to decide if $G$ has a $k$-colouring. If $k$ is not part of the input but a fixed constant, we denote this problem as $k$-{\sc Colouring}. It is well known that $3$-{\sc Colouring} is \NP-complete~\cite{Lo73}. This led to an extensive study of the complexity of {\sc Colouring} for special graph classes.

Most studied graph classes are {\it hereditary}, i.e., closed under vertex deletion. In particular, Kr{\'{a}}l et al.~\cite{KKTW01} determined, for every graph $H$, the complexity of {\sc Colouring} for {\it $H$-free graphs}, i.e., graphs that do not contain $H$ as an induced subgraph. For connected~$H$, this classification implies that {\sc Colouring} is polynomial-time solvable if $H$ is an induced subgraph of the $4$-vertex path $P_4$ and \NP-complete otherwise.
The complexity classification becomes much more involved for $(H_1,H_2)$-free graphs, see~\cite{GJPS17}, and is  not yet settled. The same holds for the classification of $k$-{\sc Colouring} for $H$-free graphs; see e.g.~\cite{CHS24,HLS22,JKMM} for recent progress.

We consider {\it monotone} graph classes, which are not only closed under vertex deletion but also under edge deletion, and in particular 
 {\it $H$-subgraph-free} graphs, i.e., graphs that do not contain some graph $H$ as a subgraph. Every $H$-subgraph-free graph is also $H$-free, but the reverse implication only holds if $H$ is a complete graph. 
A  range of graph problems on partitioning, covering, network design, width parameter computation etc.
can be fully classified on $H$-subgraph-free graphs (even on ${\cal H}$-subgraph-free graphs for finite sets of forbidden subgraphs~${\cal H}$) 
if they are (C1) polynomial-time solvable on graphs of bounded treewidth; (C2) 
\NP-complete for {\it subcubic} graphs (graphs of maximum degree at most~$3$), and (C3) stay \NP-complete under edge subdivision of subcubic graphs~\cite{JMOPPSV25}. 
%Such problems are called C123. For a connected graph $H$, every C123-problem is polynomial-time solvable on $H$-subgraph-free graphs if $H$ is a subcubic tree with at most one vertex of degree~$3$ and \NP-complete otherwise.
However, {\sc Colouring} does not satisfy C2, as it is polynomially solvable for subcubic graphs due to Brooks' Theorem.  Our research question is:

\medskip
\noindent
{\it Is it still possible to classify the complexity of {\sc Colouring} for $H$-subgraph-free graphs?}

\medskip
\noindent
We restrict ourselves to {\it connected} $H$, as even 
with this restriction, the classification is much more problematic than
the classification for $H$-free graphs, 
as we now discuss.

\medskip
\noindent
{\bf Known Results.}
Emden-Weinert, Hougardy and Kreuter~\cite{EHK98} (see also~\cite{KKTW01,KL07}) proved that $3$-{\sc Colouring} is \NP-complete for graphs of girth~$g+1$ for every $g\geq 3$,
which implies the same
 for $H$-subgraph-free graphs if $H$ has a cycle $C_g$. Hence, we may assume $H$ is a tree. Garey, Johnson and Stockmeyer~\cite{GJS76} proved that {\sc $3$-Colouring} is \NP-complete for graphs of maximum degree~$4$, so for $H$-subgraph-free graphs if $H$ is a tree of maximum degree at least~$5$.
A standard dynamic programming algorithm
(and also a result on clique-width~\cite{KR03})
shows that {\sc Colouring} is polynomial-time solvable  
for graphs of bounded treewidth, and thus for $H$-subgraph-free graphs if $H$ is a subcubic tree with at most one vertex of degree~$3$~\cite{RS84}.

\begin{figure}[t]
\centering\scalebox{1}{\input{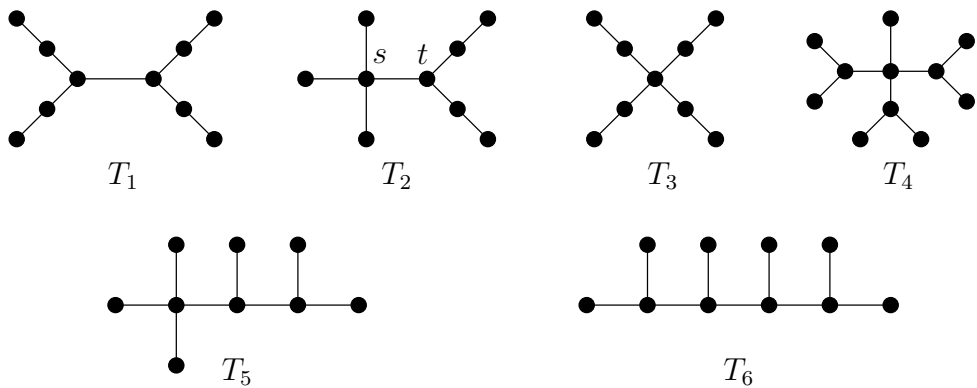}}
\vspace*{-0.2cm}
\caption{The trees $T_1,\ldots, T_6$ from~\cite{GPR15}; note that $T_1=\mathbb{H}_1^{2,2,2,2}$ and $T_3=S_{2,2,2,2}$.
\label{fig:T}}
\vspace*{-0.1cm}
\end{figure}

\begin{figure}[b]
\centering
\begin{tikzpicture}
\begin{scope}[scale=0.8]
	\node[vertex] (v1) at (0,2){};
	\node[vertex, label=left:$u$] (v2) at (0,1){};
	\node[vertex] (v3) at (0,0){};
	\node[vertex] (v4) at (1,2){};
	\node[vertex, label=right:$v$] (v5) at (1,1){};
	\node[vertex] (v6) at (1,0){};
	
	\draw[edge](v1)--(v2);
	\draw[edge](v2)--(v5);
	\draw[edge](v2)--(v3);
	\draw[edge](v4)--(v5);
	\draw[edge](v5)--(v6);
	
\end{scope}

\begin{scope}[shift = {(3,0)},scale=0.8]
	\node[vertex] (v1) at (0,2){};
	\node[vertex, label=left:$u$] (v2) at (0,1){};
	\node[vertex] (v3) at (0,0){};
	\node[vertex] (v4) at (3,2){};
	\node[vertex, label=right:$v$] (v5) at (3,1){};
	\node[vertex] (v6) at (3,0){};
	
	\node[vertex] (u1) at (1,1){};
	\node[vertex] (u2) at (2,1){};
	
	\draw[edge](v1)--(v2);
	\draw[edge](v2)--(v3);
	\draw[edge](v4)--(v5);
	\draw[edge](v5)--(v6);
	\draw[edge](v2)--(u1);
	\draw[edge](u2)--(v5);
	\draw[edge](u1)--(1.3,1);
	\node[color = gray](dots) at (1.65,1){$\dots$};
	
	\draw[br](2.925,0.875)--(0.075,0.875);
	\scriptsize
	\node[](i) at (1.5, 0.6){$i$ edges};	
\end{scope}
\end{tikzpicture}
\caption{The ``H''-graph $\mathbb{H}=\mathbb{H}_1=\mathbb{H}_1^{1,1,1,1}$ and for $i\geq 1$, the graph $\mathbb{H}_i=\mathbb{H}_i^{1,1,1,1}$ obtained from
$\mathbb{H}$ by subdividing $i-1$ times the {\it horizontal} edge $uv$. The other edges of $\mathbb{H}$ are the {\it vertical} edges.}\label{fig-Hstar} 
\end{figure}
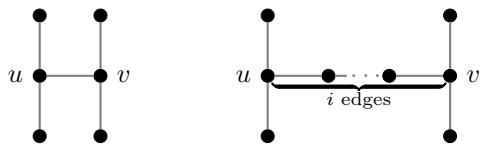

Golovach et al.~\cite{GPR15} started a systematic study on {\sc Colouring} for $H$-subgraph-free graphs. To explain their results, we first consider the trees of maximum degree at most~$4$ displayed in Figure~\ref{fig:T}.  The {\it subdivision} of an edge $e=uv$ in a graph replaces $uv$ with a new vertex $w$ made adjacent  to $u$ and $v$.  For a graph $F$, a {\it subdivided}~$F$ is a graph obtained from $F$ by subdividing each edge zero or more times. The ``H''-graph $\mathbb{H}=\mathbb{H}_1$  is the graph from Figure~\ref{fig-Hstar} that looks like the letter ``H''. For $d\geq c\geq b\geq a\geq 1$ and $i\geq 1$, a subdivided ``H''-graph $\mathbb{H}_i^{a,b,c,d}$ is obtained from $\mathbb{H}$ by subdividing the horizontal edge $i-1$ times and the vertical edges $a-1$, $b-1$, $c-1$ and $d-1$ times, respectively. For example, $\mathbb{H}=\mathbb{H}_1^{1,1,1,1}$ and $T_1=\mathbb{H}_1^{2,2,2,2}$. For $p\geq 0$, let $T_2^p$ be the tree obtained from $T_2$ by subdividing $p$ times the edge $st$, so $T_2^0=T_2$.
Golovach et al.~\cite{GPR15} proved that {\sc $3$-Colouring} is \NP-complete for $H$-subgraph-free graphs if  $H$ is a tree with at least two vertices of degree~$4$, or contains as a subgraph either a subdivided $T_1$, i.e., a tree $\mathbb{H}_i^{a,b,c,d}$ with $i\geq 1$ and $a\geq 2$; or a tree $T_2^p$ with $0\leq p\leq 9$, or a tree from $\{T_3,T_4,T_5,T_6\}$. They also proved that {\sc Colouring} is polynomial-time solvable for $H$-free graphs if $H$ is a forest of maximum degree at most~$4$ on at most seven vertices, which
includes the cases where $H\in \{\mathbb{H}_1^{1,1,1,1},\mathbb{H}_2^{1,1,1,1},\mathbb{H}_1^{1,1,1,2}\}$. 

Johnson et al.~\cite{JMPPSV23} focussed on the case where $H$ is a subdivided star $S_{p,q,r,s}$ of maximum degree~$4$ for some
$s\geq r\geq q\geq p\geq 1$, which is obtained from the $5$-vertex star $K_{1,4}$ by subdividing each edge $p-1$, $q-1$, $r-1$ and $s-1$ times, respectively; for example, $K_{1,4}=S_{1,1,1,1}$ and $T_3=S_{2,2,2,2}$. From the above it follows that {\sc Colouring} on $S_{p,q,r,s}$-subgraph-free graphs is polynomial-time solvable if $(p,q,r,s)\in \{(1,1,1,1),(1,1,1,3),(1,1,2,2)\}$ and \NP-complete, even for $k=3$, if $p\geq 2$. Johnson et al.~\cite{JMPPSV23} extended these results by showing polynomial-time solvability  if $p=q=1$ and $r\geq 1$. In fact, they showed this not only for {\sc Colouring} but for every graph problem that (i) can be solved in polynomial time on subcubic graphs;
(ii) can be solved in polynomial time on graphs of bounded treedepth; and (iii) can be solved on graphs with proper bridges by using a polynomial-time reduction to a family of instances on graphs that are either of bounded treedepth or subcubic. Here, a {\it proper bridge} in a connected graph~$G$ is an edge $e$ that is not incident to a vertex of degree~$1$ such that $G-e$ is disconnected. As other examples of such problems, they gave {\sc (Independent) Feedback Vertex Set}, {\sc Connected Vertex Cover} and {\sc Matching Cut}.

\medskip
\noindent
{\bf Our Focus.}
The graph $S_{1,1,1,1}$ is also known as the graph $\mathbb{H}_0$, which is part of an infinite sequence 
$\mathbb{H}_0,\mathbb{H}_1,\mathbb{H}_2,\ldots$, where $\mathbb{H}_i=\mathbb{H}_i^{1,1,1,1}$ for $i\geq 1$; see also Figure~\ref{fig-Hstar}.
%These graphs play an important role in the complexity classification of C123-problems. The reason is that
We note that the infinite set ${\cal M}=\{C_3,C_4,\ldots,\mathbb{H}_0,\mathbb{H}_1,\mathbb{H}_2,\ldots\}$ is a maximal antichain in the poset of connected graphs under the subgraph relation.
Conditions~C2 and~C3 ensure that for every finite set ${\cal M}'\subseteq {\cal M}$, C123-problems are \NP-complete on ${\cal M}'$-subgraph-free graphs. If  C2 or C3 is not satisfied, as is the case with {\sc Colouring}, we must consider graphs $H$ from ${\cal M}$.
For example, this has been done by Lozin et al.~\cite{LMPPSSV24} for {\sc Hamilton Cycle}, {\sc $k$-Induced Disjoint Paths}, {\sc Star $3$-Colouring} and {\sc $C_3$-Colouring}, which all satisfy C2 but not C3. For each of these, they left challenging open cases of graphs $\mathbb{H}_i$.
As {\sc Colouring} is \NP-complete for $(C_3,\ldots,C_{g-1})$-subgraph-free graphs for every $g\geq 3$~\cite{EHK98},  we only need to focus on graphs from $\{\mathbb{H}_0,\mathbb{H}_1,\mathbb{H}_2,\ldots\}$. However, in order to generalize the known results, we will consider {\it every} subdivided $\mathbb{H}_0$ and {\it every} subdivided $\mathbb{H}_1$.

\medskip
\noindent
{\bf  Our Results.}
We fully classify the complexity of {\sc Colouring} on $H$-subgraph-free graphs if $H$ is a subdivided $\mathbb{H}_0$ (subdivided $S_{1,1,1,1})$, which is a subdivided star of maximum degree~$4$, or subdivided $\mathbb{H}_1$, which is a subcubic tree with exactly two vertices of degree~$3$.
We do this by proving that all the remaining cases are polynomial-time solvable. Combining this with the results from~\cite{GPR15,JMPPSV23} yields: 

\begin{theorem}\label{t-mainh}
For $d\geq c\geq b\geq a\geq 1$ and $i\geq 1$, {\sc Colouring} on $\mathbb{H}_i^{a,b,c,d}$-subgraph-free graphs is polynomial-time solvable if $a=1$ and \NP-complete, even for $k=3$, if $a\geq 2$.
\end{theorem}

\begin{theorem}\label{t-mainstar}
For  $s\geq r\geq q\geq p\geq 1$, {\sc Colouring} on $S_{p,q,r,s}$-subgraph-free graphs is polynomial-time solvable if $p=1$ and \NP-complete, even for $k=3$,  if $p\geq 2$.
\end{theorem}

\noindent
Combining Theorems~\ref{t-mainh} and~\ref{t-mainstar} with the other aforementioned known results~\cite{EHK98,GPR15,RS84} shows that the following two cases are still unresolved:\\[-10pt]
\begin{itemize}
\item $H$ is a tree of maximum degree~$4$ with exactly one vertex of degree~$4$ and at least one vertex of degree~$3$; or
\item  $H$ is a subcubic tree with at least three vertices of degree~$3$. 
\end{itemize}

\noindent
{\bf Methodology.}
As {\sc Colouring} is polynomial-time solvable for subcubic graphs by Brook's Theorem, we may assume the input graph $G$ has a vertex of degree at least~$4$.  We may also assume that $G$ is $2$-connected; else we can consider its $2$-connected components separately. Hence, $G$ has no 
proper bridge.  The approach of Johnson et al.~\cite{JMPPSV23} was to show that a connected $S_{1,1,s,s}$-subgraph-free graph $G$ with a vertex of degree at least~$4$ and no proper bridges has bounded treedepth (note that we may assume $r=s$, as every $S_{1,1,r,s}$-subgraph-free graph is $S_{1,1,s,s}$-free if $s\geq r$). Their idea is to take a vertex $u$ of degree at least~$4$. If $G-u$ has unbounded treedepth, $G-u$ must contain a long path $P$. Due to the $2$-connectivity $u$ can be connected to a middle vertex of $P$ by two vertex-disjoint paths, and this yields an $S_{1,1,s,s}$, a contradiction.

Johnson et al.~\cite{JMPPSV23} explained that $3$-edge connectivity is needed for their approach to work for $S_{1,r,r,r}$-subgraph-free graphs and that it is not clear if a suitably modified graph-structural result can be obtained in this way. They left the case where $H=S_{1,r,r,r}$ as an open problem. We solve this open problem by showing a new structural decomposition theorem. This theorem states that if an $S_{1,r,r,r}$-subgraph-free graph~$G$ of 
maximum degree at least~$4$ contains no proper bridge alongside some other {\it special substructure}, then it must have bounded treedepth. To prove this statement we assume $G$ has unbounded treedepth. By a result of Galvin, Rival, and Sands~\cite{GALVIN19827}, $G$, being $S_{1,r,r,r}$-subgraph-free, must have a long induced path. %(as it contains no large balanced bicliques because it is $S_{1,r,r,r}$-subgraph-free). 
We take a vertex on this path and argue the existence of three long vertex-disjoint paths, giving an $S_{1,r,r,r}$, which is the desired contradiction. It then remains to show that we can preprocess instances of {\sc Colouring} in polynomial time to obtain a polynomial number of instances that satisfy the conditions of our new structural decomposition theorem.

To prove that {\sc Colouring} is polynomial-time solvable for $\mathbb{H}_i^{1,d,d,d}$-subgraph-free graphs we use a similar technique, but the details are different,
We need the $\mathbb{H}_i^{1,d,d,d}$-subgraph-free input graph to have minimum degree at least~$3$, while it may not contain some fan  
(path with dominating vertex) and another special type of substructure.
We emphasize that in both proofs the connections of the fan and other special substructures to the rest of the graph must be in some controlled way, in order to 
%as else the preprocessing to obtain instances satisfying the structural decomposition theorems will not 
preserve the $H$-subgraph-freeness. 

\medskip
\noindent
{\bf Paper Organisation.}
The additional notation and concepts that we need to show our structural results are introduced in Sections~\ref{sec-prelim} and~\ref{sec-TL-jump-chain}. Our structural results for $\mathbb{H}_i^{1,d,d,d}$-subgraph-free graphs and for $S_{1,r,r,r}$-subgraph-free graphs  are proven in Sections~\ref{sec-H} and~\ref{sec-S}, respectively. The algorithms to solve \col{} on both graph classes are presented in Section~\ref{sec-problem-specific}. We show in Section~\ref{s-wider} how our general techniques apply to other graph problems, namely {\sc Stable Cut} for $\mathbb{H}_i^{1,d,d,d}$-subgraph-free graphs
and $S_{1,r,r,r}$-subgraph-free graphs, and {\sc Feedback Vertex Set} for $S_{1,r,r,r}$-subgraph-free graphs. Finally, we finish our paper with some directions for future work in Section~\ref{sec-conc}.

\section{Preliminaries}\label{sec-prelim}
All graphs considered here are finite, simple, and undirected. That is, a \emph{graph} $G=(V,E)$ consists of a finite set $V$ of \emph{vertices} and a set $E \subseteq V^{(2)}$ of \emph{edges}, where $V^{(2)}$ is the set of $2$-element subsets of $V$. We may also denote the set of vertices and edges of $G$ by $V(G)$ and $E(G)$ respectively.
Let $G$ be a graph. We denote the {\it neighbourhood} of a vertex $v\in V(G)$ by $N_G(v)=\{ u\in V(G)\; |\; uv\in E(G)\}$. Note that we may also  write $N(v)$ if the graph is clear from the context. For a subset $S\subseteq V(G)$, we write $N_G(S)=\bigcup_{v \in S}N_G(v)$.
A graph~$H$ is a {\it subgraph} of $G$ if $H$ can be obtained from~$G$ by a sequence of vertex deletions and edge deletions. 
If $H$ is a subgraph of $G$, we write $H \subseteq G$. A graph~$H$ is an {\it induced subgraph} of $G$ if $H$ can be obtained from $G$ by a sequence of vertex deletions.
For a vertex set $S\subseteq V(G)$, we write $G[S]$ to denote the subgraph of $G$ {\it induced by} $S$, that is, the graph obtained from $G$ after deleting the vertices not in $S$. 
We further let $G-S$ denote the graph $G[V \setminus S]$.
For $S = \{u\}$, for some $u \in V(G)$, we also write $G-u$ instead of $G-\{u\}$. 
The {\it contraction} of an edge $e=uv$ in $G$ replaces $u$ and $v$ by a new vertex~$w$ that is adjacent to every vertex in $(N_G(u)\cup N_G(v)) \setminus \{u,v\}$ (without creating parallel edges).

We may refer to a path $P$ with vertices $p_1, \ldots, p_n$ and edges $p_{i-1}p_i$ for $2\leq i\leq n$ by the
sequence $(p_1,p_2,\dots,p_n)$. The \emph{length} of $P$ is its number of edges~$n-1$. 
The \emph{distance} $\dist_G(u,v)$ between two vertices $u$ and $v$ of a graph $G$ is the length of a
shortest path from $u$ to $v$.
An edge $uv \in E(G)$ is a \textit{bridge} if its removal increases the number of connected components.
A bridge is \textit{proper} if neither incident vertex has degree 1. A graph is \textit{quasi-bridgeless} if it has no proper bridge.

An \emph{elimination forest} of a graph $G$ is a rooted
forest $T$ such that $V(G)=V(T)$ and for every $uv\in E(G)$
both $u$ and $v$ are on the same root-to-leaf path of $T$. The \emph{treedepth} $\td(G)$ is the minimum
height of an elimination forest of $G$.

We further highlight the following two results from the literature regarding treedepth. Let $K_{r,s}$ be the complete bipartite graph with partition sets of size $r$ and $s$ respectively.

\begin{theorem}[\cite{NO12}]\label{thm-td-path}
    Let $G$ be a graph of treedepth at least~$d$. Then $G$ has a subgraph isomorphic to a path of length at least~$d$.
\end{theorem}

\begin{theorem}[\cite{GALVIN19827}]\label{thm:longInducedPath}
    For all $r,s,\ell\geq 1$, there is a number $c(r,s,\ell)$
    such that every $K_{r,s}$-subgraph-free graph of treedepth $c(r,s,\ell)$ has an induced $P_\ell$.
\end{theorem}

We also consider ordered lists.
For lists $A = [a_1, \ldots, a_n]$ and $B = [b_1, \ldots, b_m]$, let $A+B = [a_1, \ldots, a_n, b_1, \ldots, b_m]$. 
For paths $P = (p_1, \ldots, p_n) \subseteq G$ and $Q = (q_1, \ldots, q_m) \subseteq G - (V(P) \setminus \{p_n\})$ 
such that $p_n =q_1$, we let $P + Q$ denote the path $(p_1, \ldots, p_n, q_2, \ldots, q_m)$. 

\begin{figure}[t]
    \centering
    % \documentclass[crop,tikz]{standalone}
% \usetikzlibrary{fit,positioning}
% \begin{document}

\begin{tikzpicture}[x=1.5cm, y=1cm]

% Draw 20 vertices as circles with swapped labels
\foreach \i in {1,...,10} {
    \node[draw, circle, inner sep=2pt] (p\i) at (\i,0) {};
    \node[above=1pt] at (p\i) {\(p_{\i}\)};
    \pgfmathtruncatemacro{\negindex}{\i - 11}
    \node[below=1pt] at (p\i) {\(p_{\negindex}\)};
}

% Draw edges
\foreach \i in {1,...,9} {
    \pgfmathtruncatemacro{\j}{\i + 1}
    \draw (p\i) -- (p\j);
}

%\draw[bend right=25, thick] (p5) to (p9);
%\draw[bend left=30, thick] (p4) to (p2);

% Define and draw 4 new vertices above p_5 to p_10 as a path
% \node[draw, circle, inner sep=2pt] (q1) at (6, 0.8) {};
% \node[draw, circle, inner sep=2pt, color=teal,fill] (q2) at (7, 1) {};
% \node[draw, circle, inner sep=2pt, color=teal,fill] (q3) at (8, 1) {};
% \node[draw, circle, inner sep=2pt, color=teal,fill] (q4) at (9, 0.8) {};

% \node[draw, circle, inner sep=2pt, color=teal,fill] (p6) at (6,0) {};
% \node[draw, circle, inner sep=2pt, color=teal,fill] (p10) at (10,0) {};

% \node[draw, circle, inner sep=2pt, color=orange,fill] (p9) at (9,0) {};
% \node[draw, circle, inner sep=2pt, color=orange,fill] (p5) at (5,0) {};

% Draw edges for the path from p_5 to q1 to q2 to q3 to q4 to p_10
% \draw[thick] (p5) -- (q1) -- (q2) -- (q3) -- (q4) -- (p10);
% \draw[thick] (p6) -- (q2);

\node[draw=black, thick,dashed, rectangle, fit=(p1)(p2), inner sep=10pt, color=purple] {};
\node[draw=black, thick,dashed, rectangle, fit=(p4)(p5)(p6)(p7), inner sep=10pt, color=teal] {};
\node[draw=black, thick,dashed, rectangle, fit=(p8)(p9)(p10), inner sep=10pt, color=orange] {};
\end{tikzpicture}
% \end{document}

% Draw rectangular boxes around groups of vertices
% \node[] at (1.5,-1) {$P[1:2] = P[:2] =$};
% \node[] at (1.5,-1.5) {$P[p_1:p_2] = P[:p_2] =$};
% \node[] at (1.5,-2) {$P[:p_{-9}] = P[:-9]$};

% \node[] at (5.5,-1) {$P[5:7] = P[7:5] =$};
% \node[] at (5.5,-1.5) {$P[p_5:p_7] = P[p_7:p_5] =$};
% \node[] at (5.5,-2) {$P[5:-4] = P[-4:6]$};

% \node[] at (9,-1) {$P[8:10] = P[8:] =$};
% \node[] at (9,-1.5) {$P[p_8:p_{10}] = P[p_{10}:] =$};
% \node[] at (9,-2) {$P[p_{-3}:] = P[-3:]$};
    \caption{An illustration of the notation regarding paths. The subpath $P[p_1:p_2] = P[1:2]$ is highlighted by the purple box, $P[p_4:p_7] = P[4:7]$ is highlighted by the teal box and $P[p_8:p_{10}] = P[8:10]$ is highlighted by the orange box. Note that negative indexing may also be used i.e., $P[8:-1] = P[-3:-1] = P[8:10]$.}
    \label{fig:path_notation}
\end{figure}
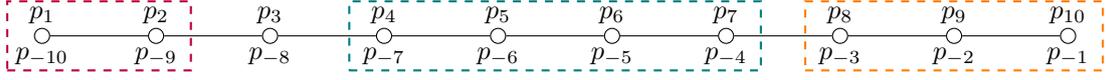

For a path $P = (p_1, \ldots, p_n)$, for every $i \in \{-n, \ldots, -1\}$, we refer to $p_{n+i+1}$ as $p_i$, i.e., $p_{-1} = p_n$. 
See Figure~\ref{fig:path_notation}.
For $i, j \in \{-n, \ldots, n\} \setminus \{0\}$, let $i' = i$, if $i \geq 1$, and $i' = n+i+1$ otherwise. 
Likewise, if $j \geq 1$, let $j' = j$, and $j' = n+j+1$ otherwise. 
If $i' \leq j'$, let $P[i:j] = (p_{i'}, p_{i'+1}, \ldots, p_{j'})$, else let $P[i:j] = (p_{j'}, p_{j'+1}, \ldots, p_{i'})$. 
For vertices $u,v \in V(P)$, if $\dist_P(p_1,u) = i-1$ and $\dist_P(p_1,v) = j-1$, then let $P[u:v] = P[i:j]$. Further, for every integer $i \in \{-n, \ldots, n\} \setminus \{0\}$ and vertex $u \in V(P)$, let $P[:i] = P[1:i]$, $P[i:] = P[i:-1]$, $P[:u] = P[p_1:u]$ and $P[u:] = P[u:p_{-1}]$.

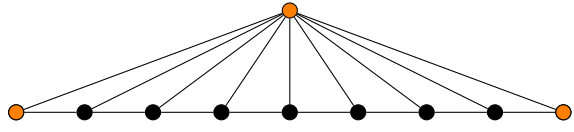
\begin{figure}
    \centering
    \begin{tikzpicture}[scale=0.9, every node/.style={circle, draw, fill=black, inner sep=2pt}]
        % Path vertices v1--v9
        \foreach \i in {1,...,9} {
            \node (v\i) at (\i,0) {};
        }
    
        % Apex vertex v0
        \node[fill=orange] (v0) at (5,1.5) {};
        \node[fill=orange] (v1) at (1,0) {};
        \node[fill=orange] (v9) at (9,0) {};
    
        % Draw path edges
        \foreach \i in {1,...,8} {
            \draw (v\i) -- (v\the\numexpr\i+1\relax);
        }
    
        % Connect apex to all path vertices
        \foreach \i in {1,...,9} {
            \draw (v0) -- (v\i);
        }
    \end{tikzpicture}
    \caption{An example of a fan $F_{10}$ (of order~$10$) with poles drawn in orange.}
    \label{fig:fan}
\end{figure}

The {\it fan} $F_{n}$ is obtained by adding a new vertex adjacent to every vertex of an $(n-1)$-vertex path $P$. We call this new vertex the \textit{centre} of the fan and those vertices at the ends of $P$ the
\textit{ends}
 of the fan. 
We also say that the centre and ends are the three \emph{poles} of the fan, see Figure~\ref{fig:fan}.
We say $G$ contains some \emph{protected fan} of order $n$, if there is some $Q \subseteq V(G)$ such that $G[Q]$ is isomorphic to $F_n$ and for all $v \in Q$, either $v$ is a pole or $N(v) \subseteq Q$.

\section{T-type and L-type Subgraphs, Jumps and Chain Extensions}\label{sec-TL-jump-chain}

The following definitions allow us to formalise the structural results which will be developed in Sections~\ref{sec-H} and \ref{sec-S}.
In particular, the key results of those sections show 
that if a graph does not contain one of those `$H$'-type graphs we are concerned with, i.e., $H$-subgraph-free with $H \in \{\mathbb{H}_m^{1,k,k,k}, S_{1,k,k,k}\}$, for some $m, k \geq 1$, then either the graph contains a large protected fan (see \textit{Preliminaries}), or there exists a pair of vertices whose removal leaves a component of bounded treedepth. As we explain in Section~\ref{sec-problem-specific}, when such a component also satisfies additional properties, we can preprocess the
input graph to an equivalent instance which is also $H$-subgraph-free but contains no such cut vertices. Towards this we define \textit{$T$-type} and \textit{$L$-type} subgraphs.

\begin{definition}[$T$-type subgraph]
A subgraph $G' \subseteq G$ is a \emph{$T$-type subgraph}, for some \emph{treedepth bound} $c \geq 1$, if $G'$ has treedepth at most~$c$ and there exists some \emph{witness set} $S \subseteq V(G) \setminus V(G')$ of size at most~$2$ such that $N(V(G')) \setminus V(G') = S$ and for each $v \in S$, $|N(v) \cap (V(G') \cup S)| \geq 2$.
\end{definition}

\begin{definition}[$L$-type subgraph]\label{def-L-type}
A subgraph $G' \subseteq G$ is an \emph{$L$-type subgraph}, for some \emph{treedepth bound} $c \geq 1$ and some \emph{length bound} $\ell \geq 1$, if $G'$ is a $T$-type subgraph with treedepth at most~$c$, the witness set $S$ has size $2$ and there exists an induced path of length at least~$\ell$ between the two vertices of $S$ in the graph $G[V(G') \cup S]$.
\end{definition}

\noindent
Note that $T$-type and $L$-type graphs may be disconnected.
We now present
Definition~\ref{d-tl} and Lemma~\ref{lem-T-type-smalltd} for \textit{minimal} $T$-type and $L$-type subgraphs, which will be especially useful in Section~\ref{sec-problem-specific}. The reason is that 
they ensure that $T$-type and $L$-type subgraphs can be preprocessed in polynomial time. This will then allow us to assume that they do not exist in our input graph $G$ anymore.

\begin{definition}[Minimal $T$-type ($L$-type) subgraph]\label{d-tl}
$G' \subseteq G$ is a \emph{minimal} $T$-type ($L$-type) subgraph if there is no $G'' \subsetneq G'$ such that $G''$ is a $T$-type ($L$-type) subgraph, with respect to the same treedepth bound (and length).
\end{definition}

\begin{lemma}\label{lem-T-type-smalltd}
    For any $c \geq 1$ and $G' \subseteq G$, if $G'$ is either a minimal $T$-type or $L$-type subgraph with treedepth bound $c$ and witness set $S$, then $\td(G[V(G') \cup S]) \leq 3c+2$. 
\end{lemma}

\begin{proof}
    We note that by definition $\td(G') \leq c$. We first suppose, for contradiction, that $G'$ contains at least~$4$ connected components. As $|S| \leq 2$, there must exist components $C_1, C_2, C_3 \subseteq G'$, such that either both $N(V(C_1)) \setminus V(C_1) = S$ and $N(V(C_2)) \setminus V(C_2) = S$ or there is some $v \in S$ such that $N(V(C_1)) \setminus V(C_1) = \{v\}$ and $N(V(C_2)) \setminus V(C_2) = \{v\}$. Further, if $G'$ is an $L$-type subgraph then we assume that for some component $C \in \{C_1, C_2, C_3\}$, there is an induced path meeting that length bound in $G[V(C) \cup S]$.

    We let $G''$ be the disjoint union of $C_1$, $C_2$ and $C_3$. Since $G'$ contains at least~$4$ connected components, $G'' \subsetneq G$. In the first case it holds by definition that
    $N(V(G'')) \setminus V(G'') = S$. Given each vertex in $S$ has some neighbour in $C_1$ and some neighbour in $C_2$, it follows that for each $v \in S$, $|N(v) \cap (V(G'') \cup S)| \geq 2$. It follows that $G''$ is a $T$-type ($L$-type) subgraph, with respect to the same treedepth bound (and length) and witness set $S$. As $G'' \subsetneq G$ this contradicts the minimality of $G'$. Likewise, in the second case there is some $v \in S$ such that $N(V(G'')) \setminus V(G'') = \{v\}$ and $|N(v) \cap (V(G'') \cup S)| \geq 2$, that is $G''$ is a $T$-type subgraph, with respect to the same treedepth bound and with the witness set $\{v\}$. Again this contradicts the minimality of $G'$.

    That is we now assume that $G'$ contains at most~$3$ connected components. Recall that each of these components has treedepth at most~$c$ and $|S| \leq 2$. As adding a new vertex increases the treedepth by at most one, it follows that $\td(G[V(G') \cup S]) \leq 3c+2$.
\end{proof}

\noindent
The core of both key structural theorems (Theorems~\ref{thrm-H} and~\ref{thrm-S}) will be in finding several long paths with a single common vertex.

Towards this, we identify a path $P$ in our graph such that, for specific vertices in $P$, if removing these vertices disconnects $P$, then we obtain either a $T$-type or an $L$-type subgraph. That is, if no $T$-type or $L$-type subgraph exists, we conclude that there must exist additional paths between specific vertices which are disjoint from~$P$.

These additional disjoint paths will then be used to construct long, internally disjoint paths between designated vertices. In turn, the existence of such paths will be used to show that our graph contains either an $S_{1,k,k,k}$ or an $\mathbb{H}_m^{1,k,k,k}$ subgraph, for $k,m \geq 1$.

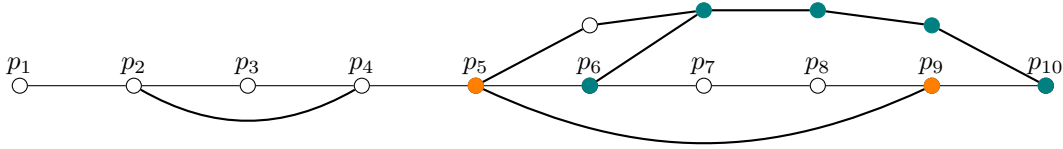
\begin{figure}[ht]
    \centering
    %\documentclass[crop,tikz]{standalone}
%\usetikzlibrary{fit,positioning}
%\begin{document}

\begin{tikzpicture}[x=1.5cm, y=1cm]

% Draw 20 vertices as circles with swapped labels
\foreach \i in {1,...,10} {
    \node[draw, circle, inner sep=2pt] (p\i) at (\i,0) {};
    \node[above=1pt] at (p\i) {\(p_{\i}\)};
    %\pgfmathtruncatemacro{\negindex}{\i - 11}
    %\node[below=1pt] at (p\i) {\(p_{\negindex}\)};
}

% Draw edges
\foreach \i in {1,...,9} {
    \pgfmathtruncatemacro{\j}{\i + 1}
    \draw (p\i) -- (p\j);
}

\draw[bend right=25, thick] (p5) to (p9);
\draw[bend left=30, thick] (p4) to (p2);

% Define and draw 4 new vertices above p_5 to p_10 as a path
\node[draw, circle, inner sep=2pt] (q1) at (6, 0.8) {};
\node[draw, circle, inner sep=2pt, color=teal,fill] (q2) at (7, 1) {};
\node[draw, circle, inner sep=2pt, color=teal,fill] (q3) at (8, 1) {};
\node[draw, circle, inner sep=2pt, color=teal,fill] (q4) at (9, 0.8) {};

\node[draw, circle, inner sep=2pt, color=teal,fill] (p6) at (6,0) {};
\node[draw, circle, inner sep=2pt, color=teal,fill] (p10) at (10,0) {};

\node[draw, circle, inner sep=2pt, color=orange,fill] (p9) at (9,0) {};
\node[draw, circle, inner sep=2pt, color=orange,fill] (p5) at (5,0) {};

% Draw edges for the path from p_5 to q1 to q2 to q3 to q4 to p_10
\draw[thick] (p5) -- (q1) -- (q2) -- (q3) -- (q4) -- (p10);
\draw[thick] (p6) -- (q2);

% \node[draw=black, thick,dashed, rectangle, fit=(p1)(p2), inner sep=10pt, color=purple] {};
% \node[draw=black, thick,dashed, rectangle, fit=(p4)(p5)(p6)(p7), inner sep=10pt, color=teal] {};
% \node[draw=black, thick,dashed, rectangle, fit=(p8)(p9)(p10), inner sep=10pt, color=orange] {};
\end{tikzpicture}
    \caption{An illustration of Definitions~\ref{def-Z^P(u,v)},~\ref{def-jump-out}, and~\ref{def-jump-out-max}. The jump $Z^P(p_6,p_{10})$ is highlighted in teal and the jump $Z^P(p_5,p_9)$ in orange. There is a negative jump out of $(p_4,p_7)$. There are three positive jumps out of $(p_4,p_7)$ and these have endpoints $(p_5,p_{9})$, $(p_5,p_{10})$ and $(p_6,p_{10})$, respectively. Here the maximum positive jump out of $(p_4,p_7)$ has endpoints $(p_6,p_{10})$, that is $(x^+(p_4,p_7), y^+(p_4,p_7)) = (p_6,p_{10})$.}
    \label{fig:jump_notation}
\end{figure}

Our next definitions and notation are illustrated in Figure~\ref{fig:jump_notation}, and will help us to formalise this reasoning and ensure that the constructed paths are indeed internally disjoint.

\begin{definition}[Jump]\label{def-Z^P(u,v)}
Let $G$ be a graph, and let $P = (p_1, \ldots, p_n) \subseteq G$ be a path in $G$. For a pair of vertices $a, b \in V(P)$, a \emph{jump} between $a$ and $b$ (with respect to $P$) is a path from $a$ to $b$ that is both edge-disjoint from $P$ and internally vertex-disjoint from $P$.
\end{definition}

\noindent
While there does not necessarily exist a unique jump between any pair of vertices, for a graph $G$ and path $P \subseteq G$, for every $u, v \in V(P)$, we fix an arbitrary jump between $u$ and $v$ and call this $Z^P(u,v)$, if such a jump exists.

\begin{definition}[Positive/negative jump out of an interval]\label{def-jump-out}
    Let $G$ be a graph, and let $P = (p_1, \ldots, p_n) \subseteq G$ be a path in $G$. 
    
    For $u, v \in V(P)$, suppose $u \in V(P[:v])$. We say that the pair $(a,b)$ are the endpoints of some \emph{jump out of the interval} $(u,v)$ if there is some jump $Z^P(a,b)$ between $a$ and $b$, with $a \in V(P[u:v]) \setminus \{u,v\}$, $b \in V(P) \setminus V(P[u:v])$. If $b \in V(P[v:])$ then we say this is a \emph{positive} jump out of $(u,v)$, else this is a \emph{negative} jump out of $(u,v)$.
\end{definition}

\begin{definition}[Maximum positive/negative jump out of an interval]\label{def-jump-out-max}
    Let $G$ be a graph, and let $P = (p_1, \ldots, p_n) \subseteq G$ be a path in $G$. For $u, v \in V(P)$, without loss of generality $u \in V(P[:v])$. We consider jumps with respect to $P$.

    Suppose there is some positive (negative) jump out of the interval $(u,v)$ with endpoints $(a,b)$ such that $a \in V(P[u:v]) \setminus \{u,v\}$ and $b \in V(P[v:])$ ($b \in V(P[:u])$ if this is a negative jump). This jump is the \emph{maximum} positive (negative) jump out of the interval $(u,v)$, if $b$ maximises $\dist_P(v,b)$ ($\dist_P(u,b)$ if this is a negative jump) and $a$ minimises $\dist_P(a,b)$. 
\end{definition}

\noindent
Note that every maximum positive jump out of an interval has the same endpoints. For a path~$P$ and vertices $u, v \in V(P)$, we let $x^+(u,v)$ and $y^+(u,v)$ denote the endpoints of a maximum positive jump out of $(u,v)$, with $x^+(u,v) \in V(P[u:v]) \setminus \{u,v\}$ and $y^+(u,v) \in V(P) \setminus V(P[u:v])$. Likewise, let $x^-(u,v)$ and $y^-(u,v)$ be the endpoints of a maximum negative jump out of $(u,v)$ with $x^-(u,v) \in V(P[u:v]) \setminus \{u,v\}$ and $y^-(u,v) \in V(P) \setminus {V(P[u:v])}$.

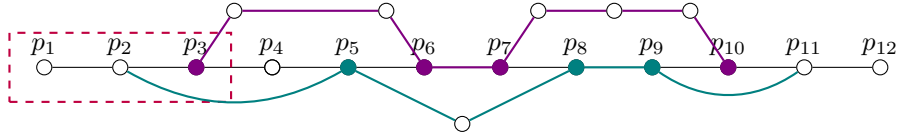
\begin{figure}[ht]
        \centering
        %\documentclass[crop,tikz]{standalone}
%\usetikzlibrary{fit,positioning}
%\begin{document}

\begin{tikzpicture}[x=1cm, y=1cm]

\foreach \i in {1,...,12} {
    \node[draw, circle, inner sep=2pt] (p\i) at (\i,0) {};
    \node[above=1pt] at (p\i) {\(p_{\i}\)};
}

% Draw edges
\foreach \i in {1,...,11} {
    \pgfmathtruncatemacro{\j}{\i + 1}
    \draw (p\i) -- (p\j);
}

\node[draw, circle, inner sep=2pt, color= violet, fill] (p3) at (3, 0) {};
\node[draw, circle, inner sep=2pt] (p4) at (4, 0) {};

\node[draw, circle, inner sep=2pt, color= violet, fill] (p6) at (6, 0) {};
\node[draw, circle, inner sep=2pt, color= violet, fill] (p7) at (7, 0) {};

\node[draw, circle, inner sep=2pt, color= violet, fill] (p10) at (10, 0) {};

\node[draw, circle, inner sep=2pt, color= teal, fill] (p5) at (5, 0) {};
\node[draw, circle, inner sep=2pt, color= teal, fill] (p8) at (8, 0) {};
\node[draw, circle, inner sep=2pt, color= teal, fill] (p9) at (9, 0) {};

\node[draw=black, thick,dashed, rectangle, fit=(p1)(p2)(p3), inner sep=10pt, color=purple] {};
\draw[bend right=30, thick, color=teal] (p2) to (p5);

\node[draw, circle, inner sep=2pt] (q1) at (3.5, 0.75) {};
\node[draw, circle, inner sep=2pt] (q2) at (5.5, 0.75) {};

\draw[thick, color=violet] (p3) -- (q1) -- (q2) -- (p6);
\draw[thick, color=violet] (p6) -- (p7);

\node[draw, circle, inner sep=2pt] (q3) at (6.5, -0.75) {};

\draw[thick, color=teal] (p5) -- (q3) -- (p8);

\node[draw, circle, inner sep=2pt] (q4) at (7.5, 0.75) {};
\node[draw, circle, inner sep=2pt] (q5) at (8.5, 0.75) {};
\node[draw, circle, inner sep=2pt] (q6) at (9.5, 0.75) {};

\draw[thick, color=violet] (p7) -- (q4) -- (q5) -- (q6) -- (p10);
\draw[thick, color=teal] (p8) -- (p9);
\draw[bend right=35, thick, color=teal] (p9) to (p11);

\end{tikzpicture}

%\end{document}
        \caption{An illustration of Definitions~\ref{def:jump_sequence}--\ref{def:path_of_jump}. $I = [(p_2,p_5),(p_3,p_6), (p_5,p_8), (p_7,p_{10}), (p_9,p_{11})]$ is a jump sequence, where $I_{2} = (p_3,p_6)$, $I_{-1} = (p_9,p_{11})$, $I^x_{2} = p_3$ and $I^y_{-1} = p_{11}$. Further, $I=T^+(p_1, p_3)$ is a maximum length positive chain extension of $(p_1, p_3)$, and      although $I$ does not describe a single path, there is an odd and an even path as described in Observation~\ref{obs-chain-paths} shown in teal and magenta respectively.}
        \label{fig:fig-chain_notation}
\end{figure}

\begin{definition}[Jump sequence]\label{def:jump_sequence}
    For a graph $G$ and a path subgraph $P = (p_1, \ldots, p_n)$ we say the list $I = [(x_1,y_1), \ldots, (x_{|I|},y_{|I|})]$ is a \emph{jump sequence} if there is a jump between $x_i$ and~$y_i$ (with respect to $P$) for every $i \in \{1, \ldots, |I| \}$, see Figure~\ref{fig:fig-chain_notation}. 
\end{definition}

\noindent
Suppose $I = [(x_1,y_1), \ldots, (x_{|I|},y_{|I|})]$ is a jump sequence. Throughout this paper, for $i \in \{1, \ldots, |I|\}$ we will let $I^x_i = x_i$, $I^y_i = y_i$ and $I_i = (x_i,y_i)$. For $i \in \{-|I|, \ldots, -1\}$ we will let $I^x_i = x_{|I|+i+1}$, $I^y_i = y_{|I|+i+1}$ and $I_i = (x_{|I|+i+1},y_{|I|+i+1})$, that is, $I_{-1} = I_{|I|} = (x_{|I|},y_{|I|})$.

\begin{definition}[Chain extension]\label{def-chain}
    Let $G$ be a graph with some path subgraph $P = (p_1, \ldots, p_n)$. Let $u,v \in V(P)$. Without loss of generality, we say that $u \in V(P[:v])$.
    We say a jump sequence $I = [(x_1,y_1), \ldots, (x_{|I|},y_{|I|})]$ is a \emph{positive chain extension} of $(u,v)$ (with respect to $P$), if $(x_1, y_1) = (x^+(u, v), y^+(u, v))$ and $(x_i, y_i) = (x^+(u, y_{i-1}), y^+(u, y_{i-1}))$ for every $i \in \{2, \ldots, |I|\}$. See Figure~\ref{fig:fig-chain_notation}. Likewise, $I$ is a \emph{negative chain extension} of $(u,v)$ if $(x_1, y_1) = (x^-(u, v), y^-(u, v))$ and $(x_i, y_i) = (x^-(v, y_{i-1}), y^-(v, y_{i-1}))$.    
\end{definition}

\noindent
Let $P$ be a path in a graph $G$. For a pair of vertices $u,v \in V(P)$, we let $T^+(u,v)$ be the maximum length positive chain extension from $(u,v)$ and let $T^-(u,v)$ be the maximum length negative chain extension from $(u,v)$. We note that $T^+(u,v)$ and $T^-(u,v)$ do not necessarily exist but if they do they are unique. 
Note that as $T^+(u,v)$ is maximal, there is no positive jump out of $(u, T^+(u,v)^y_{-1})$. Likewise, there is no negative jump out of $(T^-(u,v)^y_{-1},v)$.

\begin{definition}[Paths associated with a jump sequence]\label{def:path_of_jump}
    For a jump sequence $I$ we say the path $Z^P(x_1,y_1) \cup \bigcup_{2 \leq i \leq |I|} (P[y_{i-1}:x_{i}] \cup Z^P(x_i,y_i))$ is \emph{associated} with $I$, where, for every $i \in \{1, \ldots, |I| \}$, $Z^P(x_i,y_i)$ is that arbitrary jump fixed between $x_i$ and $y_i$.
\end{definition}

\noindent
We highlight that there is not necessarily a path associated with a given jump sequence. Moreover, depending on which arbitrary jumps have been chosen, the resulting vertices may or may not form a valid path. Nonetheless, whenever it is relevant, we will justify the existence of such a path and explain why it suffices to consider an arbitrary path.

\begin{observation}\label{obs-chain-paths}
    Let $G$ be a graph containing some path $P \subseteq G$. Suppose that for $u,v \in V(P)$, $T$ is a chain extension of $(u,v)$, with respect to $P$. There exist paths associated with the jump sequences $[(T^x_i,T^y_i): i\mod 2=1, 1 \leq i \leq |T|]$ and $[(T^x_i,T^y_i): i\mod 2=0, 1 \leq i \leq |T|]$, which we call the odd and even path of $T$, respectively. Further, if $|T| \geq 2$ the odd and the even path of $T$ are disjoint.
\end{observation}

\begin{proof}
    We will prove this observation by induction on the length of $|T|$. We will assume that $T$ is a positive chain extension. The case where $T$ is a negative chain extension will follow symmetrically. For the notation, we direct the reader to the above definitions and notation regarding jump sequences. Suppose $|T|=2$. The paths $Z^P(T^x_{1}, T^y_{1})$ and $Z^P(T^x_{2}, T^y_{2})$ are disjoint else there exists some jump from $T^x_1$ to $T^y_2$, thus contradicting that the jump between $T^x_1$ and $T^y_1$ was a maximum jump.

    We now assume that for some $\delta \geq 2$, every chain extension of $(u,v)$ with length $\delta$ results in an odd path and an even path which are disjoint. Suppose $|T| = \delta+1$. If for some $i \in \{1,\ldots, |T|-1\}$, the path $Z^P(T^x_{-1}, T^y_{-1})$ intersects the path $Z^P(T^x_{i}, T^y_{i})$ at some vertex which is not $T^x_{-1}$ then there is some jump from $T^x_{i}$ to $T^y_{-1}$, hence the jump between $T^x_{i}$ and $T^y_{i}$ were not maximum.

    As $T - (T^x_{-1}, T^y_{-1})$ is a chain extension of $(u,v)$ with length $\delta$, by assumption the paths associated with the jump sequences $[(T^x_i,T^y_i): i\mod 2=1, 1 \leq i \leq |T|-1]$ and $[(T^x_i,T^y_i): i\mod 2=0, 1 \leq i \leq |T|-1]$ are disjoint. We note that by definition $T^x_{-1}$ and $T^y_{-1}$ are the endpoints of the maximum jump from $(u,T^y_{-2})$ but not of a jump from $(u,T^y_{-3})$ (or $(u,v)$ if $|T| =3$), else the jump between $(T^x_{-3},T^y_{-3})$ (or $(u,v)$ if $|T| =3$) was not maximum. It follows that, $T^x_{-1} \in V(P[T^y_{-3}:T^y_{-2}]) \setminus \{T^y_{-2}\}$ (or $V(P[v:T^y_{-2}]) \setminus \{T^y_{-2}\}$ if $|T| =3$). Suppose $\delta$ is even.
    
    It follows that the path $P[T^y_{-3}:T^x_{-1}] + Z^P(T^x_{-1}, T^y_{-1}) - (T^y_{-3})$ is disjoint from both the odd and the even path of $T - (T^x_{-1}, T^y_{-1})$. That is, if $\delta$ is even, then the odd path of $T - (T^x_{-1}, T^y_{-1})$ extended via $P[T^y_{-3}:T^x_{-1}] + Z^P(T^x_{-1}, T^y_{-1})$ is disjoint from the even path of $T - (T^x_{-1}, T^y_{-1})$, that is the odd and the even paths of $T$ are disjoint. That is, if $\delta$ is even, then our observation holds in the inductive case.
    
    Suppose now that $\delta$ is odd, then the even path of $T - (T^x_{-1}, T^y_{-1})$ extended via $P[T^y_{-3}:T^x_{-1}] + Z^P(T^x_{-1}, T^y_{-1})$ is disjoint from the even path of $T - (T^x_{-1}, T^y_{-1})$, that is the odd and the even paths of $T$ are disjoint. Thus concluding the inductive case and so also the proof of this observation.  
\end{proof}

\section{Subdivided $\mathbb{H}_1$ Graphs: Structural Results}\label{sec-H}

We now develop those structural results when we forbid some subdivided $\mathbb{H}_1$ graphs. As outlined in Section~\ref{sec-TL-jump-chain}, jumps will play a central role in the proofs. Lemmas~\ref{lem-path-ab} and~\ref{lem-cut-pair-H} form the foundation of our reasoning: the first shows that we can either identify a desired subgraph or reason about the endpoints of jumps, while the second guarantees the existence of such jumps in the graph.

\begin{lemma}\label{lem-path-ab}
	For any $m, k \geq 1$, let $G$ be a graph with minimum degree at least~$3$. Suppose $G$ contains vertices $a,b$ and an induced path $P$ with length at least $3m^2+7k+2$, such that:
    \begin{itemize}
        \item $a \in V(P[k+1:-(k+1)])$
        \item $\dist_P(a,b) \geq 2k$
        \item there is a jump from $a$ to $b$ with respect to $P$, see Definition~\ref{def-Z^P(u,v)} 
    \end{itemize}
    Then $G$ contains a protected fan on $m+k+2$ vertices or some $\mathbb{H}_m^{1,k,k,k}$ subgraph.
\end{lemma}
    
\begin{proof}
    Let $G$ be a graph with those properties described in the Lemma statement and let $P = (p_1, \ldots, p_{\ell})$ be that induced path of length at least $3m^2+7k+2$. 
    Given $P$ has length at least $3m^2+7k+2$, after possibly relabelling $a$ and $b$, we claim that at least one of the following must hold. 
    \begin{align*}
        \dist_P(a,b) & \geq m^2+3k+1 \\
        \dist_P(p_{1}, a) & \geq  m^2 + 2k +1 \text{ and } a \in V(P[:b]) \text{ or}\\
        \dist_P(a, p_{-1}) & \geq  m^2 + 2k +1 \text{ and } a \in V(P[b:])
    \end{align*}
    Suppose towards a contradiction that none of the above inequalities hold. We will consider the case where $a \in V(P[:b])$, the case where $a \in V(P[b:])$ will follow symmetrically. Given that $\dist_P(p_{1}, a) \leq  m^2 + 2k$ and $\dist_P(a,b) \leq m^2+3k$ it follows that $\dist_P(p_{1}, b) \leq 2m^2 + 5k$ and so $\dist_P(b, p_{-1}) \geq m^2 + 2k+1$. That is exchanging $a$ and $b$ we find that $a \in V(P[b:])$ and $\dist_P(a, p_{-1}) \geq  m^2 + 2k+1$, a contradiction. We may therefore assume that at least one of the above inequalities hold.

    We will now define disjoint paths $A$, $A'$ and $B$. The path $A$ will have length $m^2 + 2k$, the path $A'$ will have length $k-1$ and the path $B$ will have length $k-2$. Further, the paths $(a) + A$, $(a) + A'$ and $(b) + B$ will each be a subpath of $P$. 
    Let $(a, \ldots, b_{k-1}, \ldots, b_1, b) = P[a:b]$ be the subpath of $P$ between $a$ and $b$ and $B = (b_1, \ldots, b_{k-1})$. 
    If $\dist_P(a,b) \geq m^2+3k+1$, then let $(a,a_1, \ldots, a_{m^2+3k+1}, \ldots, b) = P[a:b]$ be the subpath of $P$ between $a$ and $b$, note that for $i \in \{1, \ldots, \dist_P(a,b)-1\}$, $a_i = b_{\dist_P(a,b)-i}$. In addition, we let $(a,a'_1, \ldots, a'_k, \ldots, p_{-1}) = P[a:]$, if $b \in V(P[:a])$, and $(p_{1}, \ldots, a'_k, \ldots, a'_1, a) = P[:a]$ otherwise. See Figure~\ref{fig-lem-path-ab-1}.
    
    \begin{figure}
        \centering
        \includegraphics[width=0.95\linewidth,page=1]{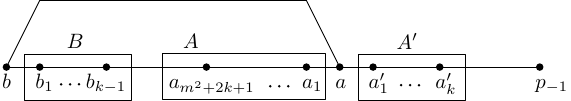}
        \caption{An illustration Lemma~\ref{lem-path-ab}, in particular of the case where $a \in V(P[b:])$ and $\dist_P(a,b) \geq m^2+3k+1$. The paths $A$, $A'$ and $B$ are highlighted alongside the vertices $a_1, \ldots, a_{m^2+3k+1}$, $a'_1, \ldots, a'_k$ and $b_1,\ldots, b_{k-1}$. We note that the vertices $a_1, \ldots$ and $b_1, \ldots$ each have a pair of labels, i.e. $a_1 = b_{\dist_P(a,b)-1}$ and $b_1 = a_{\dist_P(a,b)-1}$.}
        \label{fig-lem-path-ab-1}
    \end{figure}
    
        If $a \in V(P[:b])$ and $\dist_P(p_1, a) \geq m^2 + 2k+1$, then let $(p_{1}, \ldots, a_{m^2 + 2k+1}, \ldots, a_1,a) = P[:a]$ and $(a,a'_1, \ldots, a'_k, \ldots, b) = P[a:b]$. Else, $a \in V(P[b:])$ and $\dist_P(a, p_{-1}) \geq m^2 + 2k+1$. 
    Let $(a, a_1, \ldots, a_{m^2 + 2k+1}, \ldots, p_{-1}) = P[a:]$ and $(a,a'_1, \ldots, a'_k, \ldots, b) = P[a:b]$. Let $A = (a_1, \ldots, a_{m^2 + 2k+1})$ and $A' = (a'_1, \ldots, a'_k)$. We note that the paths $A$, $A'$ and $B$ have the properties outlined.

    Let $Z$ denote that path from $a$ to $b$ in $G - (V(P) \setminus \{a,b\})$. We will first consider the case where $Z$ has length $2$, that is without loss of generality there is a vertex $x \in V(G)$ such that $Z = (a,x,b)$.

\medskip
\noindent {\ensuremath{\vartriangleright}} {\sf \sffamily Claim~\ref{lem-path-ab}.1.}
        Either $G$ contains a $\mathbb{H}_m^{1,k,k,k}$ subgraph or the vertices $\{x, a_{m^2-m+1}, \ldots , a_{m^2+k+1}\}$ form a protected fan on $m+k+2$ vertices. 
    \begin{claimproof}
        To simplify arguments, we will also denote $a$ by $a_0$. Recall that the paths $A$, $A'$ and $B$ are disjoint subpaths of $P$ and have length $m^2 + 2k$, $k-1$ and $k-2$ respectively. We also make the following pair of observations. 
        
        Observation 1: if there is some $i \in \{m, \ldots, m^2 + k+1\}$ 
        such that $a_{i}$ is adjacent to some $y_i \in V(G \setminus P)$ with $y_i \neq x$
        and $a_{i-m}$ is adjacent to $x$, then $G$ contains $\mathbb{H}_m^{1,k,k,k}$ with degree three vertices $a_{i-m}$ and $a_{i}$; the path $(a_{i}, y_i)$ of length $1$; and paths of length at least~$k$ via $(x,b)+B$, $A[i-m:1] + (a) + A'$ and $A[i:i+k]$. 
        
        Observation 2: if there is some $i \in \{m+1, \ldots, m^2 + k+1\}$ such that $a_{i}$ is adjacent to some $y_i \in V(G \setminus P)$ with $y_i \neq x$ and $a_{i-m+1}$ is adjacent to $x$, then $G$ contains $\mathbb{H}_m^{1,k,k,k}$ with degree three vertices $x$, $a_{i}$; the path $(a_{i}, y_i)$ of length $1$; and paths of length at least~$k$ via $(x)+B$, $(x) + A'$ and $A[i:i+k]$.

        We will now apply these observations to show that, if $G$ does not contain $\mathbb{H}_m^{1,k,k,k}$ as a subgraph, then for each $\delta \in \{1, \ldots, m\}$, the vertices $a_{\delta m - (\delta-1)},\ldots,a_{\delta m}$ each have exactly one neighbour outside $P$, namely $x$. Our proof will follow by induction on $\delta$.
        We highlight that, by definition, $G$ has minimum degree at least~$3$ meaning every vertex in $A$ is either adjacent to $x$ or has some neighbour not in $V(P) \cup \{x\}$. 
        As $a_0 = a$ and $a$ is adjacent to $x$, it follows by Observation 1, that if $G$ does not contain $\mathbb{H}_m^{1,k,k,k}$ as a subgraph then $N(a_m) \setminus V(P) = \{x\}$. That is our claim holds in the base case.
        Suppose now $N(\{a_{\delta m - (\delta-1)},\ldots,a_{\delta m}\}) \setminus V(P) = \{x\}$ 
        for some $\delta \in \{1, \ldots, m-1\}$. For every $i \in \{(\delta+1) m - (\delta-1),\ldots, (\delta+1) m\}$,
        by assumption $N(a_{i-m}) \setminus V(P) = \{x\}$. That is by Observation 1, if $G$ does not contain $\mathbb{H}_m^{1,k,k,k}$ as a subgraph then $N(a_i) \setminus V(P) = \{x\}$. Likewise, for $i = (\delta+1) m - ((\delta+1)-1)$, by assumption $N(a_{i-m+1}) \setminus V(P) = \{x\}$. That is by Observation 2, if $G$ does not contain $\mathbb{H}_m^{1,k,k,k}$ as a subgraph then $N(a_i) \setminus V(P) = \{x\}$. That is, the vertices $a_{(\delta+1) m - ((\delta+1)-1)},\ldots,a_{(\delta+1) m}$ each have exactly one neighbour outside~$P$, namely $x$, thus concluding our inductive step.

        Hence, if $G$ does not contain $\mathbb{H}_m^{1,k,k,k}$ as a subgraph, then the vertices $a_{m^2 - m + 1},\ldots,a_{m^2}$ each have exactly one neighbour outside $P$, namely $x$. 
        We will now also show that, if $G$ does not contain $\mathbb{H}_m^{1,k,k,k}$ as a subgraph, then for every $i \in \{m^2+1, \ldots, m^2+k+1\}$, $N(a_i) \setminus V(P) = \{x\}$, that is the vertices $\{x, a_{m^2-m+1}, \ldots, a_{m^2+k+1}\}$ form a protected fan on $m+k+2$ vertices. 
        Assume towards a contradiction, that there is some minimum $i \in \{m^2+1, \ldots, m^2+k+1\}$ such that $N(a_i) \setminus V(P) \neq \{x\}$. 
        As $i$ is minimal and $N(a_{m^2 - m + 1},\ldots,a_{m^2}) \setminus V(P) = \{x\}$, it follows that $N(a_{i-m}) \setminus V(P) = \{x\}$. By Observation 1 either $N(a_i) \setminus V(P) = \{x\}$ or $G$ contains $\mathbb{H}_m^{1,k,k,k}$ as a subgraph.
        By assumption $N(a_i) \setminus V(P) \neq \{x\}$ and $G$ does not contain $\mathbb{H}_m^{1,k,k,k}$ as a subgraph, that is we have a contradiction and find that the vertices $\{x, a_{m^2-m+1}, \ldots, a_{m^2+k+1}\}$ form a protected fan on $m+k+2$ vertices.
    \end{claimproof}

    \noindent
    Claim~\ref{lem-path-ab}.1 implies that if $Z$ has length $2$, then our Lemma holds. We will now show that this can be also used to solve the more general case, i.e., where $Z$ has length greater than $2$. In particular, we will show that either $G$ contains $\mathbb{H}_m^{1,k,k,k}$ as a subgraph or there is some pair $a'$, $b'$ such that $a' \in V(P[k+1:-(k+1)])$, $\dist_P(a',b') \geq 2k$ and there is a path of length $2$ from $a'$ to $b'$ which is internally disjoint from $P$. That is, by Claim~\ref{lem-path-ab}.1 either $G$ contains a $\mathbb{H}_m^{1,k,k,k}$ subgraph or a protected fan on $m+k+2$ vertices.

\medskip
\noindent {\ensuremath{\vartriangleright}} {\sf \sffamily Claim~\ref{lem-path-ab}.2.}
        If $G$ does not contain $\mathbb{H}_m^{1,k,k,k}$ as a subgraph, then there is some $z \in V(Z) \setminus \{a,b\}$ such that $z$ is adjacent to $a_m$ and $a_{m \lceil \frac{2k+m}{m} \rceil}$.
    \begin{claimproof}
        We will first show that $a_m$ is adjacent to some vertex $z \in V(Z) \setminus \{a,b\}$. As $G$ has minimum degree at least~$3$, if $a_m$ is not adjacent to some vertex $z \in V(Z) \setminus \{a,b\}$ then it must have some neighbour $y \in V(G) \setminus (V(P) \cup V(Z))$.
        However, now $G$ contains $\mathbb{H}_m^{1,k,k,k}$ as a subgraph with degree~$3$ vertices $a$ and $a_m$; a path of length $1$ via $(a_m,y)$; and paths of length at least~$k$ via $(a)+Z+B$, $(a)+A'$, $A[m:m+k]$, a contradiction. That is $a_m$ is adjacent to some vertex in $z \in V(Z) \setminus \{a,b\}$.

        We will now show that for every $\delta \in \{1, \ldots, \lceil \frac{2k+m}{m} \rceil\}$, $a_{\delta m}$ is adjacent to $z$. Suppose there is some minimum $\delta \in \{1, \ldots, \lceil \frac{2k+m}{m} \rceil\}$, such that $a_{\delta m}$ is not adjacent to $z$. 
        As $G$ has minimum degree at least~$3$, $a_{\delta m}$ has some neighbour $y \in V(G \setminus P)$, such that $y \neq z$. 
        As $\delta$ is minimal and by definition $a_m$ is adjacent to $z$, it follows that $a_{\delta m-m}$ is adjacent to $z$. We highlight that for either $Z' = Z[a:z]$ or $Z' = Z[z:b]$, $y \notin V(Z')$. 
        Now we observe that $G$ contains $\mathbb{H}_m^{1,k,k,k}$ as a subgraph with degree~$3$ vertices $a_{\delta m-m}$ and $a_{\delta m}$; a path of length~$1$ via $(a_{\delta m},y)$; and paths of length at least~$k$ via $A[\delta m-m:1]+(a)+A'$, $(a_{\delta m-m})+Z'+(b)+B$ and $A[\delta m:\delta m+k]$, a contradiction.
        
        It follows that $a_{\delta m}$ is adjacent to $z$ for every $\delta \in \{1, \ldots, \lceil \frac{2k+m}{m} \rceil\}$, that is $a_{m \lceil \frac{2k+m}{m} \rceil}$ is adjacent to $z$ concluding the proof of this claim.
    \end{claimproof}

    \noindent
    Let $a' = a_m$ and $b' = a_{m\lceil \frac{2k+m}{m} \rceil}$. By definition, $a' \in V(P[k+1:-(k+1)])$ and $\dist_P(a',b') \geq 2k$. Further, from Claim~\ref{lem-path-ab}.2,  there is a path of length $2$ from $a'$ to $b'$ which is internally disjoint from $P$. That is, by Claim~\ref{lem-path-ab}.1 either $G$ contains a $\mathbb{H}_m^{1,k,k,k}$ subgraph or a protected fan on $m+k+2$ vertices which concludes the proof of this lemma.
\end{proof}  

\noindent

We also need the following lemma.

\begin{lemma}\label{lem-path-isolated} 
	For any $m, k \geq 1$, let $G$ be a graph with minimum degree at least~$3$. Suppose $G$ contains a vertex $a$ and an induced path $P$ with length at least $3m^2+7k+2$, such that:
    \begin{itemize}
        \item $a \in V(P[k+1:-(k+1)])$
        \item there is some path $Z = (a, z_1, \ldots, z_{2k-1})$ of length at least $2k-1$ in $G \setminus (P-a)$
    \end{itemize}
    Then $G$ contains a protected fan on $m+k+2$ vertices or some $\mathbb{H}_m^{1,k,k,k}$ subgraph.
\end{lemma}

\begin{proof}
    Let $G$ be a graph with those properties described in the Lemma statement and let $P = (p_1, \ldots, p_{\ell})$ be that induced path of length at least $3m^2+7k+2$. Without loss of generality, $\dist_P(a, p_{-1}) \geq \dist_P( p_1, a)$ and so $\dist_P(a, p_{-1}) \geq \frac{3m^2+7k+2}{2}-1$. Let $a_0 = a$, $A' = (a'_1, \ldots, a'_k, \ldots, p_1) = P[a:p_1] - \{a\}$ and $A = (a_0, \ldots, p_{-1}) = P[a:p_{-1}]$.
    
    We now claim that either $G$ contains some $\mathbb{H}_m^{1,k,k,k}$ subgraph or for every $\delta \in \{1, \ldots, \lceil \frac{2k}{m} \rceil\}$, the vertex $a_{\delta m}$ is adjacent to some vertex in $Z -a$. Suppose there is some minimum $\delta \in \{1, \ldots, \lceil \frac{2k}{m} \rceil\}$ such that $a_{\delta m}$ is not adjacent to some vertex in $Z -a$. As $G$ has minimum degree at least~$3$, $a_{\delta m}$ has some neighbour $y \notin V(P) \cup Z$. 
    Further, as $\delta$ is minimal and $a_0 = a$, there is some $i \in \{1, \ldots, 2k-1\}$ such that $a_{\delta m-m}$ is adjacent to $z_i$. If $i \geq k$, then let $Z' = Z[z_i:z_1]$ else let $Z' = Z[z_i:z_{2k-1}]$, note in each case $Z'$ has length at least $k-1$.
    Now $G$ contains $\mathbb{H}_m^{1,k,k,k}$ as a subgraph with degree~$3$ vertices $a_{\delta m-m}$ and $a_{\delta m}$; a path of length $1$ via $(a_{\delta m},y)$; and paths of length at least~$k$ via $A[\delta m-m:] +(a)+A'$, $(a)+Z'$ and $A[\delta m:\delta m+k]$. 
    That is, either $G$ contains some $\mathbb{H}_m^{1,k,k,k}$ subgraph or for every $\delta \in \{1, \ldots, \lceil \frac{2k}{m} \rceil\}$, the vertex $a_{\delta m}$ is adjacent to some vertex in $Z -a$.

    Let $b = a_{m \lceil \frac{2k}{m} \rceil }$. From above, if $G$ does not contain some $\mathbb{H}_m^{1,k,k,k}$ subgraph, then $b$ is adjacent to some vertex in $Z -a$. Note that this implies that there is some path from $a$ to $b$ which is internally disjoint from $P$. By definition $a \in V(P[k+1:-(k+1)])$ and $\dist_P(a,b) \geq 2k$. Note that this implies that there is some path from $a$ to $b$ which is internally disjoint from $P$. It then follows from Lemma~\ref{lem-path-ab} that either $G$ contains a protected fan on $m+k+2$ vertices or some $\mathbb{H}_m^{1,k,k,k}$ subgraph.
\end{proof}

\noindent
Using the above lemma we can now prove the following lemma.

\begin{lemma}\label{lem-cut-pair-H}
	For every $m, k \geq 1$, let $G$ be a graph with minimum degree at least~$3$. Suppose $G$ contains some induced path $P$ with length at least $3m^2+7k+2$ and vertices $u,v \in V(P[k+1:-(k+1)])$ such that $\dist_P(u,v) \geq 3$. Then at least one of the following properties must hold:
    \begin{description}
        \item[i)] $G$ contains a protected fan on $m+k+2$ vertices,
        \item[ii)] $G$ contains some $\mathbb{H}_m^{1,k,k,k}$ subgraph,
        \item[iii)] $G$ contains some $L$-type subgraph with treedepth bound $(4k-3)(\dist_P(u,v)-1)-1$ and length bound $\dist_P(u,v)-2$, see Definition~\ref{def-L-type}, or
        \item[iv)] there is some jump out of $(u,v)$ with respect to $P$, see Definition~\ref{def-jump-out}.        
    \end{description}
\end{lemma}

\begin{proof}
    We direct the reader to Section~\ref{sec-TL-jump-chain} for definitions and notation regarding jumps and $L$-type subgraphs.

    Suppose there exists some graph $G$ and vertices $u,v \in V(G)$ as described in the Lemma statement. We highlight that if iv) does not hold, that is there is no jump out of $(u,v)$, then there is some connected component $C$ in $G- \{u,v\}$ such that $V(P[u:v]) \subseteq V(C) \cup \{u,v\}$ and $V(C) \cap V(P \setminus V(P[u:v]))) = \varnothing$. Suppose $\td(C) \geq (4k-3)(\dist_P(u,v)-1)$, then by Theorem~\ref{thm-td-path}, $C$ contains a path of length at least $(4k-3)(\dist_P(u,v)-1)$. Given that $|V(P) \cap V(C)| = \dist_P(u,v)-1$, it follows that $C - V(P)$ contains a path of length at least $(4k-3)$. Let $Q = (q_1, \ldots, q_{4k-2})$, denote this path. As $C$ is a connected component, it contains some path from $V(P[u:v])$ to $V(Q)$. Let $p \in V(P[u:v]) \setminus \{u,v\}$ and $q_i \in V(Q)$ be vertices such that there is a path from $p$ to $q_i$ in $C$ which is internally disjoint from $P[u:v]$ and $Q$. 
    Without loss of generality we may assume that this path from $p$ to $q_i$ consists of a single edge. Let $Q' = (q_i, \ldots, q_{i+2k-2})$, if $i \leq 2k$, and $Q' = (q_i, \ldots, q_{i-(2k-2)})$ otherwise. As $p \in V(P[k+1:-(k+1)])$ and there is some path $(p) + Q'$ of length at least $2k-1$ in $G \setminus (P-p)$, from Lemma~\ref{lem-path-isolated}, we find that either $G$ contains a protected fan on $m+k+2$ vertices or some $\mathbb{H}_m^{1,k,k,k}$ subgraph. That is either i) or ii) is satisfied.

    It follows that if i), ii) and iv) do not hold, then $\td(C) \leq (4k-3)(\dist_P(u,v)-1)-1$ and $N(V(C)) \setminus V(C) = \{u,v\}$. We now claim that there is some $L$-type subgraph $C' \subseteq C$ with a witness pair $\{\hat{u},\hat{v}\}$. Let $u'$ be that distinct vertex in $N(u) \cap V(P[u:v])$ and $v'$ be that distinct vertex in $N(v) \cap V(P[u:v])$.
    As $G$ has minimum degree at least~$3$ and $C$ is a connected component in $G-\{u,v\}$, it follows that $|N(u') \cap V(C)|, |N(v') \cap V(C)| \geq 2$. If $|N(u) \cap V(C)| \geq 2$, then let $\hat{u} = u$, else, let $\hat{u} = u'$. 
    Likewise, if $|N(v) \cap V(C)| \geq 2$, then let $\hat{v} = v$, else, let $\hat{u} = v'$. If $\hat{u} \neq u$, then $N(u) \cap V(C) = \{u'\}$, symmetrically, if $\hat{v} \neq v$, then $N(v) \cap V(C) = \{v'\}$. It follows that $C - \{\hat{u},\hat{v}\}$ consists of possibly multiple disjoint connected components in the graph $G- \{\hat{u},\hat{v}\}$
    each with treedepth at most that of $C$.
    By definition $|N(u') \cap (V(C) \setminus \{\hat{u},\hat{v}\})|, |N(v') \cap (V(C) \setminus \{\hat{u},\hat{v}\})| \geq 2$, that is, $C - \{\hat{u},\hat{v}\}$ is an $L$-type subgraph with a witness set $\{\hat{u},\hat{v}\}$, thus concluding the proof of this lemma. 
\end{proof}

\noindent
We now have those necessary components to prove our main structural theorem regarding subdivided $\mathbb{H}_1$ graphs.

\begin{theorem}\label{thrm-H}
    For any $k,m \geq 1$, there is some function $c(k, m)$, such that every graph~$G$ with treedepth at least $c(k, m)$ and minimum degree at least~$3$ contains either some protected fan on $m+k+2$ vertices, an $L$-type subgraph with treedepth at most $(4k-3)(8k^2-6k+2m+8)-1$ and length bound $m+k$ or some $\mathbb{H}_m^{1,k,k,k}$ subgraph.
\end{theorem}

\begin{proof}
    Let $G$ be a graph with minimum degree at least~$3$. We note that $K_{4k+m+1,4k+m+1}$ contains $\mathbb{H}_m^{1,k,k,k}$ as a subgraph, that is if $G$ contains $K_{4k+m+1,4k+m+1}$ then our theorem holds. That is we assume that $G$ is $K_{4k+m+1,4k+m+1}$-subgraph-free. 
    By Theorem~\ref{thm:longInducedPath}, there is some function $c(k, m)$ such that if $\td(G) \geq c(k, m)$ then $G$ contains some induced path $P$ with length at least $8k^2-6k+3m^2+8$. As $8k^2-6k+3m^2+8 > 3m^2+7k+2$, the length of $P$ is sufficient for the application of 
    Lemmas~\ref{lem-path-ab}, \ref{lem-path-isolated} and~\ref{lem-cut-pair-H}. 
    We direct the reader to Section~\ref{sec-TL-jump-chain} for definitions and notation regarding $L$-type subgraphs, jumps and chain extensions. 
    In particular, recall that for every pair $a,b \in V(P)$ we have fixed some arbitrary jump between $a$ and $b$, if such a jump exists. We let $Z^P(a,b)$ denote this arbitrary jump.

Let $q$ be the middle vertex of $P$ and $r \neq r'$ be the pair of vertices such that $\dist_P(q,r)=\dist_P(q,r')= \lceil \frac{k+m+2}{2} \rceil$. Note that $\dist_P(r,r') \geq k+m+2$. Whenever we refer to a $L$-type subgraph, it is understood to be with respect to treedepth bound $(4k-3)(8k^2-11k+m+5)-1$ and length bound $m+k$ unless stated otherwise.
    
    \begin{myclaim}\label{clm-chain-ext-h}
        Either $G$ contains some $L$-type subgraph, a protected fan on $m+k+2$ vertices, some $\mathbb{H}_m^{1,k,k,k}$ subgraph or there is some chain extension of $r,r'$ with size at least~$2k$.
    \end{myclaim}
    \begin{claimproof}
        This proof is illustrated by Figure~\ref{fig-H-T-paths}.
        As $\dist_P(q,r)=\dist_P(q,r') \leq \frac{m+k+4}{2}$, $r,r' \in V(P[k+1:-(k+1)])$ and $\dist_P(r,r') \leq m+k+4$. 
        By Lemma~\ref{lem-cut-pair-H} either $G$ contains some protected fan on $m+k+2$ vertices, some $\mathbb{H}_m^{1,k,k,k}$ subgraph, an $L$-type subgraph with treedepth bound at most $(4k-3)(m+k+4-1)-1 \leq (4k-3)(8k^2-6k+2m+8)-1$ and length bound at least $\dist_P(r,r')-2 \geq m+k$, or there is some jump out of $(r,r')$. 
        In the first three cases, this implies that our claim holds.

        That is, we assume that there is some jump out of $(r,r')$. Without loss of generality we may assume this is a positive jump. Let $T^+(r,r')$ be the maximum positive chain extension of $(r,r')$. 
        We note that, by definition, for $i \in \{1,\ldots, |T^+(r,r')|\}$, there is some path from $T^+(r,r')^x_i$ to $T^+(r,r')^y_i$ which is both internally vertex disjoint and edge disjoint from~$P$.
        That is by Lemma~\ref{lem-path-ab}, if $\dist_P(T^+(r,r')^x_i, T^+(r,r')^y_i) \geq 2k$ then $G$ either contains some protected fan on $m+k+2$ vertices or some $\mathbb{H}_m^{1,k,k,k}$ subgraph and our claim holds.

        Assume now that $\dist_P(T^+(r,r')^x_i, T^+(r,r')^y_i) \leq 2k-1$, for all $i \in \{1,\ldots, |T^+(r,r')|\}$. By definition $T^+(r,r')^x_1 \in V(P[r,r']) \setminus \{r,r'\}$ and so $\dist_P(q, T^+(r,r')^x_1) \leq \frac{m+k+4}{2}-1$, it therefore follows that $\dist_P(q, T^+(r,r')^y_i) \leq i(2k-1)+ \frac{m+k+4}{2}-1$ for every $i \in \{1,\ldots, |T^+(r,r')|\}$. That is either $|T^+(r,r')| \geq 2k$, in which case our claim holds, or $\dist_P(q, T^+(r,r')^y_{-1}) \leq (2k-1)^2 + \frac{m+k+4}{2}-1$.
        We highlight that $\dist_P(r, T^+(r,r')^y_{-1}) \geq \dist_P(r,r') \geq m+k+2$ and $\dist_P(r, T^+(r,r')^y_{-1}) \leq (2k-1)^2 + m+k+4 -1$.
        That is by Lemma~\ref{lem-cut-pair-H} either $G$ contains some protected fan on $m+k+2$ vertices, some $\mathbb{H}_m^{1,k,k,k}$ subgraph, an $L$-type subgraph with treedepth bound $(4k - 3)( (2k-1)^2 + m+k+4 -1 - 1) - 1 \leq (4k - 3)(8k^2-6k+2m+8) - 1$, for $k\geq 2$, or there is some jump out of $(r,T^+(r,r')^y_{-1})$. 
        In each of these cases, except for the final, our claim holds. That is, we now assume that there is some jump out of $(r,T^+(r,r')^y_{-1})$. By maximality of $T^+(r,T^+(r,r')^y_{-1})$, this is a negative jump.

        Let $T^-(r,T^+(r,r')^y_{-1})$ be the maximum negative chain extension of $(r,T^+(r,r')^y_{-1})$. Note that by definition, $T^-(r,T^+(r,r')^y_{-1})^y_{1} \in V(P[:r]) \setminus \{r\}$, that is if $T^-(r,T^+(r,r')^y_{-1})^x_{1} \in V(P[r':])$ then $\dist_P( T^-(r,T^+(r,r')^y_{-1})^x_{1}, T^-(r,T^+(r,r')^y_{-1})^y_{1}) \geq 2k$. From Lemma~\ref{lem-path-ab}, this implies that either $G$ contains a protected fan on $m+k+2$ vertices or some $\mathbb{H}_m^{1,k,k,k}$ subgraph, that is we assume that $T^-(r,T^+(r,r')^y_{-1})^x_{1} \in V(P[r:r']) \setminus \{r,r'\}$. That is $T^-(r,T^+(r,r')^y_{-1}) = T^-(r,r')$.

        By  applying the same arguments as for $T^+(r,r')$, we find that   
        either ${|T^-(r,r')| \geq 2k}$, or $G$ contains a protected fan on $m+k+2$ vertices, or $G$ contains some $\mathbb{H}_m^{1,k,k,k}$ subgraph, or $\dist_P(q, T^-(r,r')^y_{-1}) \leq (2k-1)^2 + \frac{m+k+4}{2}-1$. 
        Further, \[m+k+2 \leq \dist_P(T^-(r,r')^y_{-1}, T^+(r,r')^y_{-1}) \leq 2((2k-1)^2 + m+k+4 -1) = 8k^2-6k+2m+8.\] 
        By Lemma~\ref{lem-cut-pair-H}, we either find an $L$-type subgraph, or we obtain that there is some jump out of $(T^-(r,r')^y_{-1},T^+(r,r')^y_{-1})$. 
        Suppose this jump has endpoints $(x,y)$. By maximality of $T^-(r,r')$ and $T^+(r,r')$, it is not the case that both $x,y \in V(P[:r'])$ or that both $x,y \in V(P[r:])$, that is we may assume that one of $x,y \in V(P[r:])$ and the other is in $V(P[:r'])$. 
        Note this implies that $\dist_P(x,y) \geq 2k$.
        Applying Lemma~\ref{lem-path-ab}, we find that $G$ contains either a protected fan on $m+k+2$ vertices or some $\mathbb{H}_m^{1,k,k,k}$ subgraph. This proves the claim.
    \end{claimproof}

    \begin{figure}
        \centering
        \includegraphics[width=0.85\linewidth,page=9]{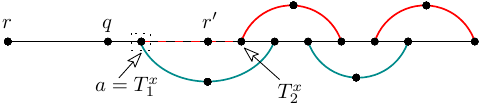}
        \caption{Illustration of the chain extension $T$ from Claim~\ref{clm-chain-ext-h}, with $T = T^+(r,r')$. The odd path of $T$ is shown in solid green, and the even path in solid red. The path $Z_1$ corresponds to the odd path of $T$ and $Z_2$ corresponds to the even path of $T$ extended to $a$. In the figure, $Z_2$ consists of both the solid and dashed red paths.}
        \label{fig-H-T-paths}
    \end{figure}

    \noindent
    From Claim~\ref{clm-chain-ext-h}, either $G$ contains some $L$-type subgraph, a protected fan on $m+k+2$ vertices, some $\mathbb{H}_m^{1,k,k,k}$ subgraph or there is some chain extension $T$ of $r,r'$ with size at least~$2k$. We will now use this chain extension to show that either $G$ contains a protected fan on $m+k+2$ vertices or some $\mathbb{H}_m^{1,k,k,k}$ subgraph. 

    By Observation~\ref{obs-chain-paths}, there exists a pair of disjoint paths described by $T$. As $P$ is an induced path, to recall notation regarding jump sequences, for every $i \in \{1,\ldots, |T|\}$, $Z^P(T_i)$ must have length at least~$2$. That is the odd and even path of $T$ each have length at least $2\cdot \frac{2k}{2} = 2k$.
    Let $a= T^x_1$. Let $Z_1 = (a, z^1_1, \ldots, z^1_{2k-1}, \ldots)$ denote the odd path of $T$ and $Z_2 = (a, z^2_1, \ldots, z^2_{2k-1}, \ldots)$ denote the even path of $T$ extended to $a$ via the path $P[a:T^x_2]$. If $T = T^+(r,r')$, we let $A = (a, a_1, \ldots, p_1)$ be that subpath of $P$ between $a$ and $p_1$. If $T = T^-(r,r')$, let $(a, a_1, \ldots, p_{-1})$ be that subpath of $P$ between $a$ and $p_{-1}$. Hence, $V(A) \cap Z_1 \cap Z_2 = \{a\}$.

    \begin{figure}
        \centering
        \includegraphics[width=1\linewidth,page=11]{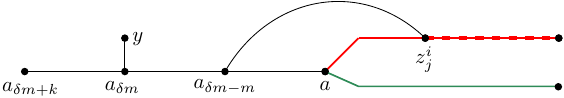}
        \caption{The path $Z_1$ is shown in solid red and the path $Z_2$ is shown in solid green. $a_{\delta m-m}$ is adjacent to some $z^i_j$, in the case depicted $i=1$ and $j < k$. The path $Z'$ is shown in dashed red and the path $Z'' = Z_2$, that is, it is it can be seen in solid green. }
        \label{fig-H-proof-Z-paths}
    \end{figure}

    Let $a_0 = a$. We claim that either $G$ contains some $\mathbb{H}_m^{1,k,k,k}$ subgraph or for every $\delta \in \{1, \ldots, \lceil \frac{2k}{m} \rceil\}$, the vertex $a_{\delta m}$ has some neighbour in $(Z_1 \cup Z_2) - \{a\}$. Suppose there is some minimum $\delta \in \{1, \ldots, \lceil \frac{2k}{m} \rceil\}$ such that $a_{\delta m}$ is not adjacent to some vertex in $(Z_1 \cup Z_2) - \{a\}$.
    As $G$ has minimum degree at least~$3$, $a_{\delta m}$ has some neighbour $y \notin V(P) \cup V(Z_1) \cup V(Z_2)$.
    Further, as $\delta$ is minimal and $a_0 = a$, there is some $i \in \{1, 2\}$ and $j \in \{1, \ldots, 2k-1\}$ such that $a_{\delta m-m}$ is adjacent to $z^i_j$. If $j \geq k$, then let $Z' = Z_i[1:j]$ else let $Z' = Z_i[j:2k-1]$, note in each case $Z'$ has length at least $k-1$.
    Let $Z'' = Z'_1$, if $i =2$, and $Z'' = Z'_2$ otherwise. See Figure~\ref{fig-H-proof-Z-paths}. 
    Now $G$ contains $\mathbb{H}_m^{1,k,k,k}$ as a subgraph with degree~$3$ vertices $a_{\delta m-m}$ and $a_{\delta m}$; a path of length $1$ via $(a_{\delta m},y)$; and paths of length at least~$k$ via $A[\delta m-m:a] +Z''$, $(a_{\delta m-m})+Z'$ and $A[\delta m:\delta m+k]$. 
    That is, either $G$ contains some $\mathbb{H}_m^{1,k,k,k}$ subgraph or for every $\delta \in \{1, \ldots, \lceil \frac{2k}{m} \rceil\}$, the vertex $a_{\delta m}$ is adjacent to some vertex in $Z'_1 \cup Z'_2$.

    Let $b = a_{m \lceil \frac{2k}{m} \rceil }$. From above, if $G$ does not contain some $\mathbb{H}_m^{1,k,k,k}$ subgraph, then $b$ is adjacent to some vertex $z$ in $Z'_1 \cup Z'_2$. By definition, there is some $i \in \{1, \ldots, |T|\}$ such that $z \in V(Z^P(T_i))$. It follows that there is some path $Z^P(b, T^y_i)$. As $\dist_P(b, T^y_i) \geq 2k$, it follows from Lemma~\ref{lem-path-ab} that either $G$ contains a protected fan on $m+k+2$ vertices or some $\mathbb{H}_m^{1,k,k,k}$ subgraph. This concludes the proof of this theorem.
\end{proof}

\section{Subdivided Stars ($\mathbb{H}_0$): Structural Results}\label{sec-S}

This section is dedicated to subdivided star graphs. That is, we will prove the following theorem:

\begin{theorem}\label{thrm-S}
    For any $k \geq 1$, let $c = 16(2k-1)(k-1)$. If $G$ is proper bridgeless, contains some degree~$4$ vertex, has treedepth at least $8(7k^3+15k^2-\frac{4k}{9}+3)^2+6$ and no $T$-type subgraph, with treedepth bound~$c$, then $G$ contains $S_{1,k,k,k}$ as a subgraph.
\end{theorem}

Although the arguments to prove this result have several similarities with those used for subdivided $\mathbb{H}_0$ graphs, there are also crucial differences. In particular, unlike Theorem~\ref{thm-col-H}, here we require not only that the long paths intersect in a single vertex, but also that this vertex has degree at least~$4$. 

To show Theorem~\ref{thrm-S} we begin with the following structural result from the literature.

\begin{theorem}[\cite{JMPPSV23}]\label{thm-S_11kk}
    Let $q, r \geq 1$. The subclass of connected $S_{1,1,q,r}$-subgraph-free graphs that are not subcubic and are quasi-bridgeless has treedepth at most $2(q+r+3)^2+6$.
\end{theorem}

\noindent
As $G$ is proper bridgeless, has maximum degree at least~$4$ and treedepth at least $8(\ell+3)^2+6$, for $\ell = 7k^3+15k^2-\frac{4k}{9}$, it follows from Theorem~\ref{thm-S_11kk} that $G$ contains some $S_{1,1, \ell,\ell}$ as a subgraph. Let $x$ be the centre of this $S_{1,1, \ell,\ell}$. Note that $x$ is the middle vertex of some path $P$ of length $2\ell$.

The central idea of the proof
of Theorem~\ref{thrm-S} 
is to show that there exists some vertex $x$ with degree at least $4$ and three paths each with length at least~$2k$ sharing the single common vertex $x$. Observation~\ref{obs-S_3->S_1kkk} will then be used to find some $S_{1,k, k,k}$.

\begin{observation}\label{obs-S_3->S_1kkk}
    For every $k \geq 1$, if $G$ contains some $S_{2k,2k,2k}$ subgraph with centre~$x$ and $\deg(x) \geq 4$, then $G$ contains $S_{1,k,k,k}$ as a subgraph.
\end{observation}

\begin{proof}
    Let $G$ be a graph that contains a subgraph isomorphic to $S_{2k,2k,2k}$. Let $x$ denote the centre of this subdivided star, and for each $i \in \{1,2,3\}$ and $j \in \{1, \ldots, 2k\}$, let $f^i_j$ denote the $j$th vertex along the $i$th branch.
    
    If $x$ has degree at least~$4$, then $x$ has some neighbour $y \notin \{f^1_1,f^2_1,f^3_1\}$. Suppose $y \neq f^i_j$ for any $i \in \{1,2,3\}$ and $j \in \{2, \ldots, k\}$. Then $G$ contains $S_{1,k,k,k}$ with centre $x$, one branch of length $1$ by the edge $xy$ and three branches of length $k$ by the paths $(f^1_1, \ldots, f^1_k)$, $(f^2_1, \ldots, f^2_k)$ and $(f^3_1, \ldots, f^3_k)$, respectively.

    Therefore, we may assume that $y = f^i_j$ for some $i \in \{1,2,3\}$ and $j \in \{2, \ldots, k\}$. Without loss of generality, let $i=1$. In this case, $G$ again contains $S_{1,k,k,k}$ as a subgraph, with centre~$x$, one branch of length $1$ via the edge $xf^1_1$, and three branches of length $k$ given by the paths $(f^1_{j}, \ldots, f^1_{j+k})$, $(f^2_1, \ldots, f^2_k)$ and $(f^3_1, \ldots, f^3_k)$, respectively.
\end{proof}

\noindent
 Lemma~\ref{lem-cut-pair-S} is crucial to construct the three paths, which must intersect at a vertex of degree at least~$4$, implying several new conditions
in contrast to its counterpart, Lemma~\ref{lem-cut-pair-H}.
 
\begin{lemma}\label{lem-cut-pair-S}
   For every $k \geq 1$, let $G$ be some $S_{1,k,k,k}$-subgraph-free graph. Suppose that $G$ contains some path $P$ with length at least $4k-2$ as a subgraph.

   Suppose there exist $u, v \in V(P[2k+1 : - (2k+1)])$ and some $x \in V(P[u : v]) \setminus \{v\}$ with $|N(x) \setminus V(P)| \geq 2$, such that:
   \begin{description}
       \item[i)] there is no jump from $x$ to some $y \in V(P) \setminus V(P[u:v])$, with respect to $P$ (see Definition~\ref{def-Z^P(u,v)});
       \item[ii)] there exist internally disjoint paths $D^1_{v}, D^2_{v}$ from $x$ to $v$;
       \item[iii)] either $x = u$, or there exist internally disjoint paths $D^1_{u}, D^2_{u}$ from $x$ to $u$;
       \item[iv)] if $x \neq u$, then the paths $D^1_{v}$, $D^2_{v}$ and $D^1_{u}$ are internally disjoint.
   \end{description}
   Then either $G$ contains some $S_{1,k,k,k}$ as a subgraph; $G$ contains some $T$-type subgraph with treedepth bound $16(2k-1)(k-1)$; or there exists some jump out of $(u,v)$.
\end{lemma}

\begin{proof}
    Let $G$ be a graph with a path $P = (p_1, \ldots, p_\ell) \subseteq G$, for some $\ell \geq 4k+3$ and some pair of vertices $u, v$ as described above.

    \begin{myclaim}\label{clm-d-paths}
        If $u \neq x$, then there exist paths $\hat{D}^1_{v}$, $\hat{D}^2_{v}$, $\hat{D}^1_{u}$, $\hat{D}^2_{u}$, such that $\hat{D}^1_{v}$, $\hat{D}^2_{v}$ are internally disjoint $x$-$v$-paths; $\hat{D}^1_{u}$, $\hat{D}^2_{u}$ are internally disjoint $x$-$u$-paths; $\hat{D}^1_{v}$, $\hat{D}^2_{v}$, $\hat{D}^1_{u}$ are internally disjoint; and $\hat{D}^1_{u}$, $\hat{D}^2_{u}$, $\hat{D}^1_{v}$ are internally disjoint.
    \end{myclaim}
    
        \begin{claimproof}
        Suppose $D^2_{u}$ is internally disjoint from either $D^1_{v}$ or $D^2_{v}$. If $D^2_{u}$ is internally disjoint from $D^1_{v}$, then let $\hat{D}^1_{v}= D^1_{v}$ and $\hat{D}^2_{v}= D^2_{v}$. Otherwise, let $\hat{D}^1_{v}= D^2_{v}$ and $\hat{D}^2_{v}= D^1_{v}$. We note that now the paths $\hat{D}^1_{v}$, $\hat{D}^2_{v}$, $\hat{D}^1_{u}= D^1_{u}$, $\hat{D}^2_{u}= D^2_{u}$ have the desired disjointness properties.
        
        We now assume that $D^2_{u}$ intersects both the path $D^1_{v}$ and $D^2_{v}$. Let $d$ be that final vertex in $D^2_{u}$, such that $d \in V(D^1_{v}) \cup V(D^2_{v})$, that is the subpath of $D^2_{u}$ from $d$ to $u$ contains no vertex from $V(D^1_{v}) \cup V(D^2_{v})$ except for $d$. We can now define a path from $x$ to $u$ which intersects only one of $D^1_{v}$ and $D^2_{v}$, if $d \in V(D^2_{v})$ this is via $D^2_{v}[:d] + D^2_{u}[d:]$. Otherwise, this is via $D^1_{v}[:d] + D^2_{u}[d:]$. Let $\hat{D}^1_{v} = D^1_{v}$, $\hat{D}^2_{v} = D^2_{v}$, $\hat{D}^1_{u} = D^1_{u}$ and $\hat{D}^2_{u} = D^2_{v}[:d] + D^2_{u}[d:]$ in the first case and $\hat{D}^1_{v} = D^2_{v}$, $\hat{D}^2_{v} = D^1_{v}$, $\hat{D}^1_{u} = D^1_{u}$ and $\hat{D}^2_{u} = D^1_{v}[:d] + D^2_{u}[d:]$ in the second. We find that the resulting $\hat{D}^1_{v}$, $\hat{D}^2_{v}$, $\hat{D}^1_{u}$, $\hat{D}^2_{u}$ have the desired disjointness properties.
    \end{claimproof}
   
    \noindent
    By Claim~\ref{clm-d-paths}, if $u \neq x$, then we can define paths $\hat{D}^1_{v}$, $\hat{D}^2_{v}$, $\hat{D}^1_{u}$, $\hat{D}^2_{u}$, with the properties outlined above. Let $D^1_{v} = \hat{D}^1_{v}$, $D^2_{v} = \hat{D}^2_{v}$, $D^1_{u} = \hat{D}^1_{u}$, $D^2_{u} = \hat{D}^2_{u}$.
        
    Let $C_1$ be that component of $G - \{u,v\}$ containing the path $D^1_{v} - \{u,v\}$. If $C_1$ contains some vertex $y \in V(P) \setminus V(P[u:v])$, then there is some path from $x$ to $y$ in $G[V(C_1) \cup \{x\}]$, without loss of generality $y$ is the single vertex from $V(P) \setminus V(P[u:v])$ on this path. From i), there is no jump from $x$ to $y$. It follows that this path from $x$ to $y$ must contain some vertex from $P[u:v] - \{u,v\}$, let $x' \in V(P[u:v]) \setminus \{u,v\}$ be that vertex which is closest to $y$. By definition there is a jump between $x'$ and $y$ and this is a jump out of the interval $(u,v)$. It follows, if $C_1$ contains some vertex $y \in V(P) \setminus V(P[u:v])$, the lemma holds.   
    
    Hence, in the following we assume that $C_1$ contains no vertices from $V(P) \setminus V(P[u:v])$. Further, if $D^1_{v}$ has length at least $2k+1$, then the paths $P[:x]$, $D^1_v - \{v\}$ and $D^2_{v} + P[v:]$ form a $S_{2k,2k,2k}$ subgraph with centre $x$. By Observation~\ref{obs-S_3->S_1kkk}, it follows that $G$ contains some $S_{1, k,k,k}$, hence our Lemma holds. Symmetrically, the same follows for $D^2_{v}$, $D^1_{u}$ and $D^2_{u}$, that is, each of $D^1_{v}$, $D^2_{v}$, $D^1_{u}$ and $D^2_{u}$ have length at most $2k$.
    
    We first consider the case where $u = x$, without loss of generality we may assume that $x \in P[:v]$. Suppose $\td(C_1) \geq 8(2k-1)(k-1)$. Then by Theorem~\ref{thm-td-path}, $C_1$ contains a path of length at least $8(2k-1)(k-1)$. Given that $|(V(D^1_{v}) \cup V(D^2_{v})) \cap V(C_1)| \leq 2(2k-1)$, it follows that $C_1 - (V(D^1_{v}) \cup V(D^2_{v}))$ contains a path of length at least $4(k-1)$.
    Let $Q = (q_1, \ldots, q_{4k-3})$, denote this path. As $C_1$ is a connected component, it contains some path from $V(D^1_{v})$ to $V(Q)$. Let $d \in V(D^1_{v}) \cup V(D^2_{v}) \setminus \{u,v\}$ and $q_i \in V(Q)$ be vertices such that there is a path from $d$ to $q_i$ in $C_1$ which is internally disjoint from $D^1_{v}$, $D^2_{v}$ and $Q$. Without loss of generality we may assume that $d \in V(D^1_{v})$ and this path from $d$ to $q_i$ consists of a single edge.
    
    Let $Q' = (q_i, \ldots, q_{i+2k-1})$, if $i \leq 2(k-1)$, and $Q' = (q_i, \ldots, q_{i-(2k-1)})$ otherwise. The paths $P[:x]$, $D^2_{v} \cup P[v:]$ and $D^1_{v}[:d] \cup Q'$ form a $S_{2k,2k,2k}$ subgraph with centre $x$, from Observation~\ref{obs-S_3->S_1kkk} $G$ contains $S_{1,k,k,k}$ as a subgraph and so our lemma holds.
    
    That is, we assume $\td(G[C_1]) < 8(2k-1)(k-1)$. If $(V(D^1_{v}) \cup V(D^2_{v})) \setminus \{x,v\} \subseteq C_1$ then as the paths $D^1_{v}$, $D^2_{v}$ are internally disjoint, $|N(x) \cap (V(C_1) \cup \{v\})|, |N(v) \cap (V(C_1) \cup \{x\})| \geq 2$. It follows that $G[C_1]$ is a $T$-type subgraph with a witness set $\{x,v\} = \{u,v\}$ and so our lemma holds. 
    Else, let $C_2$ be that component of $G - \{x,v\}$ containing the path $\hat{D}^2_{v} - \{x,v\}$. By symmetry $\td(G[C_2]) < 8(2k-1)(k-1)$. Given the paths $D^1_{v}$ and $D^2_{v}$ are internally disjoint, $|N(v) \cap (V(C_1) \cup V(C_2) \cup \{x,v\})| \geq 2$. It follows that $G[V(C_1) \cup V(C_2)]$ is a $T$-type subgraph with a witness set $\{x,v\}=\{u,v\}$. The case where $u=x$ is now complete.

    We now assume $u \neq x$. Suppose $\td(G[C_1]) \geq 16(2k-1)(k-1)$. Then by Theorem~\ref{thm-td-path}, $G[C_1]$ contains a path of length at least $16(2k-1)(k-1)$. Given $(V(D^1_{v}) \cup V(D^2_{v}) \cup V(D^1_{u}) \cup V(D^2_{u})| \leq 4(2k-1)$, it follows that $C_1 - (V(D^1_{v}) \cup V(D^2_{v}) \cup V(D^1_{u}) \cup V(D^2_{u}))$ contains some path of length at least $4(k-1)$. 
    
    Let $Q = (q_1, \ldots, q_{4k-3})$, denote this path. Let $d \in (V(D^1_{v}) \cup V(D^2_{v}) \cup V(D^1_{u}) \cup V(D^2_{u})) \setminus \{u,v\}$ and $q_i \in V(Q)$ be a pair such that there is a path from $d$ to $q_i$ in $C_1$ which is internally disjoint from $D^1_{v}$, $D^2_{v}$, $D^1_{u}$, $D^2_{u}$ and $Q$. Without loss of generality $d \in V(D^1_{v})$. We also assume that this path consists of a single edge, else we replace this edge by the respective path in the below reasoning. Let $Q' = (q_i, \ldots, q_{i+2k-1})$, if $i \leq 2(k-1)$, and $Q' = (q_i, \ldots, q_{i-(2k-1)})$ otherwise. The paths $D^1_{u} \cup P[:u]$, $D^2_{v} \cup P[v:]$ and $D^1_{v}[:d] \cup Q'$ form a $S_{2k,2k,2k}$ subgraph with centre $x$, from Observation~\ref{obs-S_3->S_1kkk}, $G$ contains $S_{1,k,k,k}$ as a subgraph and so our lemma holds. That is, we assume that $\td(G[C_1]) < 16(2k-1)(k-1)$. We note as $x \in C_1$, it follows that $(V(D^1_{v}) \cup V(D^2_{v}) \cup V(D^1_{u}) \cup V(D^2_{u})) \setminus \{u,v\} \subseteq C_1$. As the paths $D^1_{v}$, $D^2_{v}$ are internally disjoint $|N(v) \cap (C_1 \cup \{u\})| \geq 2$ and as the paths $D^1_{u}$, $D^2_{u}$ are internally disjoint $|N(u) \cap (C_1 \cup \{v\})| \geq 2$. That is $G[C_1]$ is a $T$-type subgraph with a witness set $\{x,v\}=\{u,v\}$ and concluding the proof of this lemma.
\end{proof}

\begin{figure}[t]
    \centering    \includegraphics[width=0.85\linewidth,page=15]{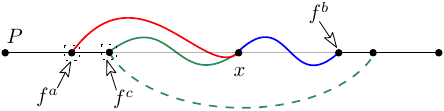}
    \caption{An illustration of Theorem~\ref{thrm-S}, as described in the proof sketch. The active frontier path is drawn in solid red, the inert frontier path is drawn in solid blue and the candidate path is drawn in solid green. In the case depicted, in the inductive step the candidate path is extended via that jump drawn in dashed green. The previous inert frontier becomes the new candidate, the previous active frontier is becomes the new inert frontier and the previous candidate becomes the new active frontier.}
    \label{fig-s-sketch}
\end{figure}

\noindent
We now sketch the proof of Theorem~\ref{thrm-S}.
Recall that our goal is to find three paths of length at least~$2k$ with a single common vertex~$x$. These paths will be constructed inductively. In every step of the induction, we either find some $S_{1,k, k,k}$ subgraph, or a pair of vertices $b^-,b^+ \in V(P)$ such that  $b^-$, $b^+$ and $x$ and $P$ meet the conditions to apply Lemma~\ref{lem-cut-pair-S}.
By Lemma~\ref{lem-cut-pair-S}, we either find some $S_{1,k, k,k}$, some $T$-type subgraph (with the desired treedepth bound) or some jump out of $(b^-,b^+)$. 
In the first case, the theorem follows, in the second we get a contradiction to the assumption that $G$ contains no such $T$-type subgraph. That is we assume that such a jump exists.
These jumps are used to construct three paths each of length at least~$2k$ sharing the single common vertex $x$, that is our theorem follows.

The paths are built inductively. At the end of each inductive step we obtain three paths which share the single common vertex $x$, together with a pair of vertices $b^-$, $b^+$ which satisfy the conditions of Lemma~\ref{lem-cut-pair-S}. Moreover, in every block of three steps each of the paths grows in length by at least one, as ensured by the addition of a new jump. It follows that after $6k$ steps each of these paths has length at least~$2k$, as desired.

To ensure that these paths have the desired properties, we label one of these paths as the active frontier, one of these as the inert frontier and the other as the candidate. Let $f^a$ be the final vertex of the active frontier, $f^b$ be the final vertex of the inert frontier, and $f^c$ be the final vertex of the candidate path. We find that $f^c \in P[f^a:x]$ and either, $f^a \in P[:x]$ and $f^b \in P[:x]$, or,  $f^a \in P[x:]$ and $f^b \in P[x:]$. See Figure~\ref{fig-s-sketch}.

During the inductive step, the inert frontier remains unchanged, the candidate path is extended by at least one edge (via the addition of a new jump), and the active frontier may either stay the same or increase in length. 
In the next step, the previous inert frontier becomes the new candidate, and one of the remaining two paths (either the former active frontier or the former candidate) takes the role of active frontier. If the previous active frontier does not increase in one step, it becomes the inert frontier in the step after, that is in the following step it is extended. This rotation ensures that after every three steps each of the three paths has grown by at least one.

To maintain disjointness, we introduce auxiliary properties ensuring that the extended paths intersect only at $x$ and that we can identify a pair of vertices $b^-$, $b^+$ such that, together with $x$, they satisfy the conditions of Lemma~\ref{lem-cut-pair-S}. Concretely, we find the candidate and active frontiers can be extended to one of $b^-$, $b^+$, while a subpath of the active frontier path, together with the inert frontier, extends to the other. These extensions of the active, inert and candidate paths remain internally disjoint, thus meeting the conditions for Lemma~\ref{lem-cut-pair-S}.  

We now present the proof of Theorem~\ref{thrm-S} in full detail.

\medskip
\noindent
{\bf Theorem~\ref{thrm-S} (restated).}
{\it
    For any $k \geq 1$, let $c = 16(2k-1)(k-1)$. If $G$ is proper bridgeless, contains some degree~$4$ vertex, has treedepth at least $8(7k^3+15k^2-\frac{4k}{9}+3)^2+6$ and no $T$-type subgraph, with treedepth bound~$c$, then $G$ contains $S_{1,k,k,k}$ as a subgraph.
}

\begin{proof}
    For $k \geq 1$, let $c = 16(2k-1)(k-1)$. We direct the reader to Section~\ref{sec-TL-jump-chain} for definitions and notation regarding $T$-type subgraphs, jumps and chain extensions. In this proof, whenever we refer to a $T$-type subgraph, it is understood to be with respect to the treedepth bound $c$ unless explicitly stated otherwise. Let $G$ be a graph with maximum degree at least~$4$, treedepth at least $8(7k^3+15k^2-\frac{4k}{9}+3)^2+6$ and no $T$-type subgraph.
    
    Let $\ell = 7k^3+15k^2-\frac{4k}{9}$. As $G$ has maximum degree at least~$4$, treedepth at least $2(2\ell+3)^2+6$ and is proper bridgeless, it follows from Theorem~\ref{thm-S_11kk}, that $G$ contains some $S_{1,1,\ell,\ell}$ subgraph. Let $x$ denote the vertex at the centre of this $S_{1,1,\ell,\ell}$ and let $A^+= (a^+_1, \ldots, a^+_\ell)$ and $A^-= (a^-_1, \ldots, a^-_\ell)$ be those paths such that $(x) + A^+$ and $(x) + A^-$ correspond to that pair of paths of length $\ell$ in $S_{1,1,\ell,\ell}$. Let $P = (a^-_\ell, \ldots, a^-_1,x,a^+_1, \ldots, a^+_\ell)$. Note we will consider all jumps with respect to $P$. As noted in Section~\ref{sec-TL-jump-chain}, for every pair $a,b \in V(P)$ we will fix some arbitrary jump between $a$ and $b$, if such a jump exists, and denote this by $Z^P(a,b)$.
    
    In the following we will show that $G$ contains a $S_{2k,2k,2k}$ subgraph with centre $x$. Observation~\ref{obs-S_3->S_1kkk} then implies that $G$ contains $S_{1,k,k,k}$ as a subgraph and thus the theorem follows. 
    To show the existence of a $S_{2k,2k,2k}$ subgraph, we will either describe some $S_{2k,2k,2k}$ subgraph with centre $x$, or inductively define pairs $(b^-_1,b^+_1), \ldots, (b^{-}_{6k},b^{+}_{6k})$ together with the jump sequences $I^i_1 \subseteq \ldots \subseteq I^i_{6k}$ for $i \in \{1,2,3\}$. We will then show that the jump sequences $I^1_{6k}$, $I^2_{6k}$, $I^3_{6k}$ each have a corresponding path and that these paths can be extended such that they have a single common vertex $x$. We will use these paths to construct a $S_{2k,2k,2k}$ subgraph.

    Intuitively, for $\delta \in \{1, \ldots, 6k-1\}$, in the inductive step from $\delta$ to $\delta+1$ we will add at least one jump to at least one of $I^1_\delta$, $I^2_\delta$ or $I^3_\delta$. We further show that after $3$ steps each of $I^1_\delta$, $I^2_\delta$ and $I^3_\delta$ has increased in length by at least~$1$ so after at most~$6k$ steps we find $S_{2k,2k,2k}$.
    
    The vertices $b^-_{\delta-1}$, $b^+_{\delta-1}$, $b^-_{\delta}$, $b^+_{\delta}$, $b^-_{\delta+1}$, $b^+_{\delta+1}$ will be such that $P[b^-_{\delta-1}:b^+_{\delta-1}] \subseteq P[b^-_{\delta}:b^+_{\delta}] \subseteq P[b^-_{\delta+1}:b^+_{\delta+1}]$. If a jump between $(a,b)$ is added to one of $I^1_\delta$, $I^2_\delta$ or $I^3_\delta$, then $a, b \notin V(P[b^-_{\delta-1}:b^+_{\delta-1}]) \setminus \{b^-_{\delta-1},b^+_{\delta-1}\}$.  
    
    Further, either every jump added to $I^1_\delta$, $I^2_\delta$ and $I^3_\delta$ will have an endpoint in $V(P[b^+_{\delta}:b^+_{\delta+1}]) \setminus \{b^+_{\delta+1}\}$, in which case we say that $b^+_\delta$ is \textit{active}, or every jump added to $I^1_\delta$, $I^2_\delta$ and $I^3_\delta$ will have an endpoint in $V(P[b^-_{\delta}:b^-_{\delta+1}]) \setminus \{b^-_{\delta+1}\}$, in which case we say that $b^-_\delta$ is \textit{active}.
    
    To formalise these ideas and allow us to do this, the pairs of vertices and jump sequences will satisfy the following properties.   
    \begin{description}
        \item[P1.] $\dist_P(x, b^+_1), \dist_P(x, b^-_1) \leq 4k^2-1$
                
        \item[P2.] Exactly one of $b^+_\delta$ and $b^-_\delta$ is active. If $b_\delta^+$ is active we have the following properties, else the properties hold after interchanging $b^+_\delta$, $b^+_{\delta+1}$ with $b^-_\delta$, $b^-_{\delta+1}$, and vice versa.
        \begin{enumerate}[a)]
            \item $b^-_{\delta+1}$ is active,
            \item $b^-_{\delta+1} = b^-_{\delta}$,
            \item there is no negative jump out of $(b^{-}_{\delta},b^{+}_{\delta})$, and
            \item $\dist_P(b^+_{\delta}, b^+_{\delta+1}) \leq 4k^2-\frac{k(4\delta-14)}{3}+\frac{\delta(\delta-7)+11}{9}-1$.
        \end{enumerate}
        These sub-properties will imply that if $b^+_\delta$ is active then for every jump with endpoints $(a,b)$ added to $I^1_\delta$, $I^2_\delta$ or $I^3_\delta$, we have $a \in V(P[b^+_{\delta}:b^+_{\delta+1}]) \setminus \{b^+_{\delta+1}\}$.

        \item[P3.] One of $I^1_\delta$, $I^2_\delta$ and $I^3_\delta$ is labelled the positive frontier, another the negative frontier and another the candidate. Again suppose $b^+_\delta$ is active, otherwise, in the following, we replace $b^+_\delta$ with $b^-_\delta$ and the positive with the negative frontier accordingly.

        As $b^+_\delta$ is active, we let the positive frontier be the active frontier and the negative frontier be the inert frontier. Let  $a \neq b \neq c \in \{1,2,3\}$ be such that, $I^a_\delta$ is the active frontier, $I^b_\delta$ is the inert frontier and $I^c_\delta$ is the candidate. Using the notation introduced regarding jump sequences, for $i \in \{1,2,3\}$, if $|I^i_\delta| \geq 1$ we let $s^i_\delta = (I^i_\delta)^x_{1}$ and $f^i_\delta = (I^i_\delta)^y_{-1}$.
        
        That is, $s^i_\delta$ and $f^i_\delta$ correspond to the \textit{starting} and \textit{final} vertex of the path described by $I^i_\delta$. If $I^i_\delta = []$, then let $s^i_\delta = f^i_\delta = x$. The following properties will hold.
        \begin{enumerate}[a)]
            \item $|I^a_\delta| \geq \lfloor \frac{\delta+2}{3} \rfloor$, $|I^b_\delta| \geq \lfloor \frac{\delta+1}{3} \rfloor$ and $|I^c_\delta| \geq \lfloor \frac{\delta}{3} \rfloor$.

            \item $f^a_\delta = y^+(b^-_{\delta}, b^+_{\delta})$ with $(I^a_\delta)^x_{-1} \in V(P[b^-_{\delta}: b^+_{\delta}])$ and $b^+_{\delta} \in V(P[f^c_\delta:f^a_\delta]) \setminus \{f^a_\delta\}$.

            \item Those paths described by $I^a_\delta$, $I^b_\delta$, $I^c_\delta$ can be extended to $x$ via $P[x:s^a_\delta]$, $P[x:s^b_\delta]$ and $P[x:s^c_\delta]$ respectively. Let $D_a$, $D_b$ and $D_c$ denote these extended paths.

            \item Let $\bar{D}_a$ denote the subpath of $D_a$ from $x$ to $(I^a_\delta)^y_{-2}$. If $|I^a_\delta| =1$, then let $\bar{D}_a = (x)$. The path $\bar{D}_a$ can be extended to $b^-_\delta$, we denote the resulting path by $D^{\ext}_a$. The path $D_b$ can also be extended to $b^-_\delta$, we denote the resulting path by $D^{\ext}_b$. 
                            
            \item The paths $P[:b^-_\delta]$, $P[b^+_\delta:]$, $P[f^c_\delta:f^a_\delta]$, $D^{\ext}_a$, $D^{\ext}_b$ and $D_c$ are pairwise internally disjoint.

            \item The path $Z^P((I^a_\delta)_{-1})$ is disjoint from both $D^{\ext}_b$ and $D_c$.
        \end{enumerate}
        In our inductive step the \textit{inert} frontier will remain the same, that is $I^b_{\delta+1} = I^b_\delta$. The \textit{candidate} is the next list of pairs which will necessarily increase in length, that is $|I^c_{\delta+1}| \geq |I^c_\delta| +1$, the \textit{active} frontier may remain the same.
    \end{description}

    \begin{figure}
        \centering
        \includegraphics[width=1\linewidth,page=8]{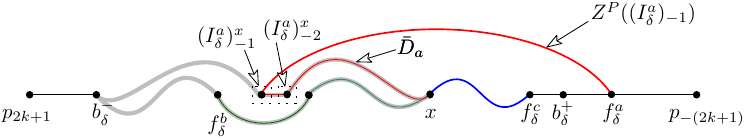}
        \caption{This diagram illustrates Properties 1-3. The paths $D_a$, $D_b$ and $D_c$ are shown in red, green and blue respectively. The paths $D^{\ext}_a$ and $D^{\ext}_b$ are both shown in gray. Note that $D^{\ext}_b$ is that gray path containing $D_b$ as a subpath and $D^{\ext}_a$ is that gray path containing $\bar{D}_a$ as a subpath.}
        \label{fig:properties}
    \end{figure}

    \noindent
    An illustration of these properties can be seen in Figure~\ref{fig:properties}. If $b^+_{6k}$ is active, from Properties P2a) and b), for every $\delta \in \{2,\ldots, 6k-1\}$, if $\delta$ is odd, then $b^+_{\delta}= b^+_{\delta}$. Further, from Property P2d), if $\delta$ is even then $\dist_P(b^+_{\delta}, b^+_{\delta+1}) \leq 4k^2-\frac{k(4\delta-14)}{3}+\frac{\delta(\delta-7)+11}{9}-1$. Note similar reasoning follows for $b^-_{6k}$ where we replace odd by even and visa versa. The case where $b^-_{6k}$ is active now follows symmetrically.
    
    Combining this with Property P1, it follows that 
    \begin{align*}
        \dist_P(x, b^+_{6k}), \dist_P(x, b^-_{6k}) 
        & \leq 4k^2 -1 + \sum^{3k}_{\delta=1} \left(4k^2-\frac{k(4\delta-14)}{3}+\frac{\delta(\delta-7)+11}{9}-1 \right) 
        \\ & = 7k^3+13k^2-\frac{4k}{9}-1.
    \end{align*}
    Likewise, for every $\delta \in \{1, \ldots, 6k\}$, we obtain that $\dist_P(x, b^+_{\delta}), \dist_P(x, b^-_{\delta}) \leq 7k^3+13k^2-\frac{4k}{9}-1$.    
    As $\dist_P(x, a^+_{\ell}) = \dist_P(x, a^-_{\ell}) = \ell$ and $\ell = 7k^3+15k^2-\frac{4k}{9}$, we get that $b^+_{\delta}, b^-_{\delta} \in V(P[2k+1 : - (2k+1)])$ for $\delta \in \{1, \ldots, 6k\}$.
    
    \medskip\noindent
    \textbf{Base case: Define $b^+_1$ and $b^-_1$ with properties P1 and P2a--d.}\\
    Let $i \in \{1, \ldots, \ell\}$ be the largest index, such that there is some jump from $x$ to $a^+_i$. That is, the largest $i$ such that there is a path from $x$ to $a^+_i$ which is both edge-disjoint and internally vertex-disjoint from $P$. If such an $i$ exists let $y^+_0 = a^+_i$, if no such $i$ exists let $y^+_0 = x$. Likewise, let $i' \in \{1, \ldots, \ell\}$ be the largest index such that there is some jump from $x$ to $a^-_{i'}$. Let $y^-_0 = a^-_{i'}$, if such an $i'$ exists, and $y^-_0 = x$ otherwise.
    
    If $\dist_P(x, y^+_0) > 2k$, then the paths $P[:x]$,  $P[x:a^+_{2k}]$ and $Z^P(x,a^+_i) + P[a^+_i:]$ each have length at least~$2k$ and share the single common vertex $x$. 
    That is $G$ contains a $S_{2k,2k,2k}$ subgraph with centre $x$ and the theorem follows. Hence, we may assume that $\dist_P(x, y^+_0) \leq 2k$ and symmetrically $\dist_P(x, y^-_0) \leq 2k$.

 \medskip
\noindent {\ensuremath{\vartriangleright}} {\sf \sffamily Claim~\ref{thrm-S}.1.}
        If $y^-_0= y^+_0 =x$, then either $\{x\}$ is the witness set for some $T$-type subgraph or $G$ contains some $S_{2k,2k,2k}$ with centre $x$. 
    \begin{claimproof}
        As $x$ has degree at least~$4$, there exist at least two distinct vertices $q$, $q' \in N(x) \setminus P$. Let $C$ and $C'$ be those components of $G-x$ containing $q$ and $q'$, respectively. Note that possibly $C = C'$. Since $y^-_0= y^+_0 =x$, every path from either $q$ or $q'$ to $P$ must contain~$x$ and hence, $(V(C) \cup V(C')) \cap V(P) = \varnothing$. Suppose $\td(C) \geq 4k+1$. From Theorem~\ref{thm-td-path}, $C$ contains some path $Q$ of length at least $4k+1$. As $C$ is connected and contains $q$, either $q \in V(Q)$ or there is some shortest path from $q$ to $Q$ in $C$. Note this path can be extended via a subpath of $Q$ with length at least $2k-2$. By definition $x$ and $q$ are adjacent, that is we obtain a path $Q'$ of length $2k$ from $x$ in $G[V(C) \cup \{x\}]$. As $(V(C) \cup V(C')) \cap V(P) = \varnothing$, $G$ contains $S_{2k,2k,2k}$ with centre $x$ and paths of length at least~$2k$ via $Q'$, $P[:x]$ and $P[x:]$. That is, our claim holds. Symmetrically, the same holds for $C'$. Hence, we may assume that $\td(C), \td(C') \leq 4k$. Note now $C \cup C'$ is a $T$-type subgraph with witness set $\{x\}$. Thus concluding the proof of this claim.
    \end{claimproof}

    \noindent
    Hence, by Claim~\ref{thrm-S}.1, if $y^-_0= y^+_0 =x$ then our theorem follows. Thus, we may assume that $y^+_0 \neq x$ or $y^-_0 \neq x$. Without loss of generality, suppose $y^+_0 \neq x$, else we exchange the roles of $A^+$ and $A^-$.
    
    We further claim that $y^-_0, y^+_0, x$ meet the conditions for Lemma~\ref{lem-cut-pair-S}. To see this, note first that, as $\dist_P(x, y^-_0), \dist_P(x, y^+_0) \leq 2k$, it follows that $y^-_0, y^+_0 \in V(P[2k+1 : - (2k+1)])$. Condition i) follows as $y^-_0, y^+_0$ were maximal. By definition of $y^+_0$ there exists a pair of internally disjoint paths from $x$ to $y^+_0$ via $P[x:y^+_0]$ and $Z^P(x,y^+_0)$, respectively, showing condition ii). 
    If $x \neq y^-_0$, then there also exist internally disjoint paths from $x$ to $y^-_0$ via $P[y^-_0: x]$ and $Z^P(y^-_0, x)$, respectively, showing condition iii) holds. Given the paths $P[x:y^+_0]$, $Z^P(y^-_0, x)$ and $P[y^-_0:x]$ are internally disjoint, condition iv) also follows. That is, by applying Lemma~\ref{lem-cut-pair-S} we find that there must exist some jump out of $(y^-_0,y^+_0)$. 
    Recalling the notation introduced regarding maximum jumps, we find that at least one of $(x^+(y^-_0,y^+_0), y^+(y^-_0,y^+_0))$ or $(x^-(y^-_0,y^+_0), y^-(y^-_0,y^+_0))$ must exist. Notice in both cases this corresponds to either a positive or negative chain extension of $(y^-_0,y^+_0)$ with length $1$.
    
    Suppose that $(x^+(y^-_0,y^+_0), y^+(y^-_0,y^+_0))$ exists. The case where $(x^-(y^-_0,y^+_0), y^-(y^-_0,y^+_0))$ is the only existing pair will follow symmetrically. Recall that $T^+(y^-_0,y^+_0)$ is the maximum positive chain extension of $(y^-_0,y^+_0)$. As $(x^+(y^-_0,y^+_0), y^+(y^-_0,y^+_0)) \in T^+(y^-_0,y^+_0)$, necessarily $|T^+(y^-_0,y^+_0)| \geq 1$. 
    Let $b^-_1 = y^-_0$ and $b^+_1 = T^+(y^-_0,y^+_0)^y_{-1}$. 
    Recall that $T^+(y^-_0,y^+_0)^y_{-1}$ corresponds to the second endpoint of the final jump of $T^+(y^-_0,y^+_0)$. Again we direct the reader to Section~\ref{sec-TL-jump-chain} for details of definitions and notation regarding jump sequences. 
    We say $b^-_1$ is active. By maximality of $T^+(y^-_0,y^+_0)$ there is no positive jump out of $(b^-_1, b^+_1)$, that is Property P2c) holds in the base case.
    
    \begin{figure}
        \centering
        \includegraphics[width=1\linewidth, page=10]{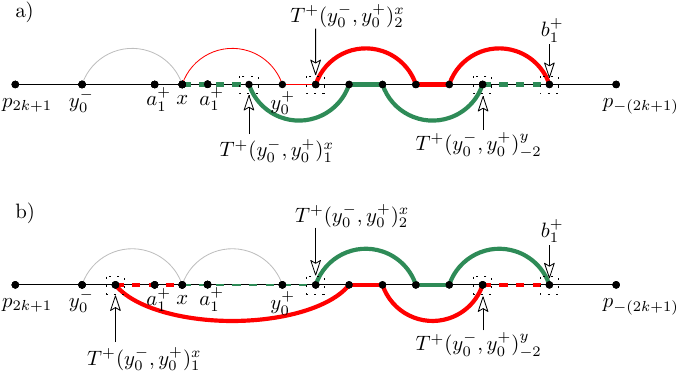}
        \caption{An illustration of the cases of Claim~\ref{thrm-S}.2. In Figure a) $T^+(y^-_0,y^+_0)^x_1 \in V(P[x:y^+_0])$ and in Figure b) $T^+(y^-_0,y^+_0)^x_1 \in V(P[y^-_0:x])$. The odd path of $T$ is drawn in solid red and the even path of $T$ is drawn in solid green. These paths can also be extended to $(x, b^+_1)$ via the dashed lines. The odd and the even path of $T^+(y^-_0,y^+_0)$ are emphasised using a thicker line. Note in Case a) the even path of $T^+(y^-_0,y^+_0)$ is a subpath of the odd path of $T$, that is, it is drawn in red. In Case b) $T = T^+(y^-_0,y^+_0)$ and so the even path of $T^+(y^-_0,y^+_0)$ is also the even path of $T$, that is, it is drawn in red.}
        \label{fig:base-c1-2}
    \end{figure}

    \smallskip
    \noindent
    Another useful idea to introduce is that of extending a path via $P$. Let $Q = (q_1, \ldots, q_r)$ be a path, and let $s, t \in P$. If $q_1, q_r \in P$, we call the path $P[s:q_1]+Q+P[q_r:t]$ \textit{$Q$ extended to $(s,t)$ via $P$}. Note we also require that the paths $P[s:q_1] - q_1$, $Q$ and $P[q_r:t] - q_r$ are disjoint. Additionally, may extend the empty path to $(s,t)$ via $P$, which will be the path~$P[s:t]$.

\medskip
\noindent {\ensuremath{\vartriangleright}} {\sf \sffamily Claim~\ref{thrm-S}.2.}
        There is some jump sequence $T$, such that either $T = T^+(y^-_0,y^+_0)$ or $T = (x,y^+_0) + T^+(y^-_0,y^+_0)$. We note that in the first case $T$ is the maximum positive chain extension of $(y^-_0,y^+_0)$ and in the second case $T$ is the maximum positive chain extension of $(a^-_1,a^+_1)$. Note that $a^-_1$ and $a^+_1$ are both adjacent to $x$ in $P$. Further, the odd path of $T$ (see Observation~\ref{obs-chain-paths} for odd and even paths) extended to $(x, b^+_1)$ via $P$ and the even path of $T$ extended to $(x, b^+_1)$ via $P$ are internally disjoint. In addition, these paths are disjoint from $P[:y^-_0]$ and either the path $P[y^-_0:x]$ or $Z^P(y^-_0, x)$.
    \begin{claimproof}
        We first note that both the odd path and the even path of $T^+(y^-_0,y^+_0)$ are disjoint from $Z^P(y^-_0,x)$ and $Z^P(x,y^+_0)$ else $y^+_0$ was not maximal. On the contrary, the paths $Z^P(y^-_0,x)$ and $Z^P(x,y^+_0)$ are not necessarily disjoint from one another.

        We first consider the case where $T^+(y^-_0,y^+_0)^x_1 \in V(P[x:y^+_0])$, see Figure~\ref{fig:base-c1-2}a). Recall that by the notation introduced regarding jump sequences, $T^+(y^-_0,y^+_0)^x_1$ is the first vertex of the odd path. Let $T = (x,y^+_0) + T^+(y^-_0,y^+_0)$. Note $T$ is a positive chain extension of $(a^-_1, a^+_1)$. Further, $T^x_1 = x$ and $T^x_2 = T^+(y^-_0,y^+_0)^x_1$, that is the odd path of $T^+(y^-_0,y^+_0)$ is a subpath of the even path of $T$, see Figure~\ref{fig:base-c1-2}a). If $|T|$ is odd, then the odd path of $T$ is a $(x, b^+_1)$ path. 
        Further, this path is disjoint from the paths $P[x:T^y_1] - \{x,T^y_1\}$ and $P[T^x_{-1}:b^+_1] - \{T^x_{-1},b^+_1\}$. By Observation~\ref{obs-chain-paths}, the odd and even path of $T$ are disjoint, that is the odd path is internally disjoint from the even path of $T$ extended to $(x, b^+_1)$ via $P$, that is our claim holds. Similarly, if $|T|$ is even, then the even path of $T$ extended via $P[x:T^x_2]$ and the odd path of $T$ extended via $P[T^y_{-2}:b^+_1]$ are internally disjoint from each other and the path $P[:x]$. That is our claim holds.

        If $T^+(y^-_0,y^+_0)^x_1 \in V(P[y^-_0:x])$ let $T = T^+(y^-_0,y^+_0)$, see Figure~\ref{fig:base-c1-2}b). By definition, $T$ is a positive chain extension of $(y^-_0,y^+_0)$. Extending the odd and the even path of $T$ to $(x,b^+_1)$ via $P$ we again find that these paths are internally disjoint. Note, if $|T| \geq 2$, then the odd path is extended by $P[x:T^x_1]$ and the even path is extended by $P[x:T^x_2]$, further, depending on if $|T|$ is odd or even one of these paths will also be extended via the path $P[T^y_{-2}:b^+_1]$. If $|T| =1$, then the even path of $T$ is empty. The odd and the even path of $T$ extended to $(x,b^+_1)$ via $P$ correspond to the paths $P[x:T^x_1] + Z^P(T^x_1,b^+_1)$ and $P[x:b^+_1]$, respectively. Note that both of them are disjoint from the path $Z^P( y^-_0, x)$, thus proving our claim.
    \end{claimproof}

\medskip
\noindent {\ensuremath{\vartriangleright}} {\sf \sffamily Claim~\ref{thrm-S}.3.}
        If $\dist_P(x, b^+_1) \geq 4k^2$, then $G$ contains some $S_{2k,2k,2k}$ subgraph with centre $x$.
    \begin{claimproof}
    We first claim, if $\dist_P(x, b^+_1) \geq 4k^2$, then at least one of the following must hold,
     \begin{itemize}
         \item $|T^+(y^-_0,y^+_0)| \geq 4k-2$,
         \item $\dist_P(T^+(y^-_0,y^+_0)^y_i, T^+(y^-_0,y^+_0)^y_{i+1}) \geq 2k- \lceil \frac{i}{2} \rceil$ for some $i \in \{1, \ldots, |T^+(y^-_0,y^+_0)|-1\}$, or
         \item  $\dist_P(y^+_0, T^+(y^-_0,y^+_0)^y_1) \geq 2k$.
     \end{itemize}
     Towards a contradiction, suppose $|T^+(y^-_0,y^+_0)| \leq 4k-3$, $\dist_P(y^+_0, T^+(y^-_0,y^+_0)^y_1) \leq 2k-1$ and $\dist_P(T^+(y^-_0,y^+_0)^y_i, T^+(y^-_0,y^+_0)^y_{i+1}) \leq 2k- \lceil \frac{i}{2} \rceil-1$ for every $i \in \{1, \ldots, |T^+(y^-_0,y^+_0)|-1\}$. Recall that $\dist_P(x, y^+_0) \leq 2k$ and $b^+_1 = T^+(y^-_0,y^+_0)^y_{-1}$.
     Hence, it follows that for every $i \in \{1, \ldots, |T^+(y^-_0,y^+_0)|\}$, we get that
     \[\dist_P(x, T^+(y^-_0,y^+_0)^y_i) \leq 2k + \sum^i_{j =1} \left(2k- \left\lceil \frac{j-1}{2} \right\rceil-1\right),\] that is, as $|T^+(y^-_0,y^+_0)| \leq 4k-3$, 
     \begin{align*}
     \dist_P(x, b^+_1) & \leq  2k + (4k-3)(2k-1) - \sum^{4k-3}_{j =1} \left \lceil \frac{j-1}{2} \right\rceil \\ &=  2k + (4k-3)(2k-1) - (2k-2)(2k-1) = 4k^2-2k+1,
     \end{align*}
     a contradiction.

     That is we assume that $|T^+(y^-_0,y^+_0)| \geq 4k-2$, $\dist_P(T^+(y^-_0,y^+_0)^y_i, T^+(y^-_0,y^+_0)^y_{i+1}) \geq 2k- \lceil \frac{i}{2} \rceil$ for some $i \in \{1, \ldots, |T^+(y^-_0,y^+_0)|-1\}$, or  $\dist_P(y^+_0, T^+(y^-_0,y^+_0)^y_1) \geq 2k$.
     
     From Claim~\ref{thrm-S}.2, for every $i \in  \{1,\ldots, |T^+(y^-_0,y^+_0)|-1\}$ there exist  paths from $x$ to $y^-_0$, $T^+(y^-_0,y^+_0)^y_{i}$ and $T^+(y^-_0,y^+_0)^y_{i+1}$ which share a single common vertex which is $x$. Further, by definition these paths have lengths at least $0$, $\lceil \frac{i}{2} \rceil +1$ and $\lceil \frac{i+1}{2} \rceil +1$ respectively. 

     Suppose that $|T^+(y^-_0,y^+_0)| \geq 4k-2$. 
     We extend that path from $x$ to $y^-_0$ via $P[:y^-_0]$. This path alongside those paths from $x$ to $T^+(y^-_0,y^+_0)^y_{4k-3}$ and $T^+(y^-_0,y^+_0)^y_{4k-2}$ describe three paths of length at least~$2k$ with a single common vertex $x$. That is $G$ contains some $S_{2k,2k,2k}$ subgraph with centre $x$.

     Suppose now that $\dist_P(T^+(y^-_0,y^+_0)^y_i, T^+(y^-_0,y^+_0)^y_{i+1}) \geq 2k- \lceil \frac{i}{2} \rceil$. We extend that path from $x$ to $y^-_0$ via the path $P[:y^-_0]$, that path from $x$ to $T^+(y^-_0,y^+_0)^y_{i}$ via the path $P[T^+(y^-_0,y^+_0)^y_{i}:T^+(y^-_0,y^+_0)^y_{i+1}] - T^+(y^-_0,y^+_0)^y_{i+1}$ and that path from $x$ to $T^+(y^-_0,y^+_0)^y_{i+1}$ via the path $P[T^+(y^-_0,y^+_0)^y_{i+1}:]$. We note that these three extended paths have length at least~$2k$ and have a single common vertex $x$.  That is $G$ contains some $S_{2k,2k,2k}$ subgraph with centre $x$. 
     
     This leaves only the case where $\dist_P(y^+_0, T^+(y^-_0,y^+_0)^y_1) \geq 2k$. Claim~\ref{thrm-S}.2 also implies that there exist paths from $x$ to $y^-_0$, $y^+_0$ and $T^+(y^-_0,y^+_0)^y_{1}$ which share a single common vertex which is $x$. If $\dist_P(y^+_0, T^+(y^-_0,y^+_0)^y_1) \geq 2k$, then extending the paths via $P[:y^-_0]$, $P[y^+_0:T^+(y^-_0,y^+_0)^y_{1}] - T^+(y^-_0,y^+_0)^y_{1}$ and $P[T^+(y^-_0,y^+_0)^y_{i+1}:]$, respectively, we obtain three paths of length at least~$2k$ with a single common vertex $x$. That is $G$ contains some $S_{2k,2k,2k}$ subgraph with centre $x$.
    \end{claimproof}

    \noindent
    It follows then from Claim~\ref{thrm-S}.3, that Property~P1 holds. 
    
    We now also claim that the vertices $b^-_1$ and $b^+_1$ meet the conditions for Lemma~\ref{lem-cut-pair-S} to be applied. To see this, we first note that from Claim~\ref{thrm-S}.3 we get $b^-_1, b^+_1 \in V(P[2k+1 : - (2k+1)])$. Further, by the maximality of $y^-_0$ and $y^+_0$, condition~i) of Lemma~\ref{lem-cut-pair-S} holds for $b^-_1$, $b^+_1$.
    By Claim~\ref{thrm-S}.2, there exists a pair of internally disjoint paths from $x$ to $b^+_1$ (via $T$) and another pair of internally disjoint paths from $x$ to $y^-_0$ (via $P[x:y^-_0]$ or $Z^P(x,y^-_0)$), implying condition~ii) and~iii) of Lemma~\ref{lem-cut-pair-S} also hold. 
    At least one of these $x$ to $y^-_0$ paths is disjoint from both $x$ to $b^+_1$ paths, that is, applying Lemma~\ref{lem-cut-pair-S}, either $y^+(b^-_1, b^+_1)$ or $y^-(b^-_1, b^+_1)$ must exist.
    As $b^-_1 = y^-_0$ and $b^+_1 = T^+(y^-_0, y^+_0)^y_{-1}$, it follows that $y^+(b^-_1, b^+_1)$ does not exist, else $T^+(y^-_0, y^+_0)$ was not maximum. That is $y^-(b^-_1, b^+_1)$ must exist.
    We let $(c^x_1, c^y_1) = (x^-(b^-_1, b^+_1), y^-(b^-_1, b^+_1))$
    By definition, $c^x_1 \in V(P[y^-_0:]) \setminus \{y^-_0\}$ 
    and $c^y_1 \in V(P[:b^-_1]) \setminus \{b^-_1\}$, that is intuitively $(x^-(b^-_1, b^+_1), y^-(b^-_1, b^+_1))$ \textit{cross} $b^-_1$.
    We say that $b^-_1$ is active. We highlight that Properties P2a--d are satisfied.

    We will now consider the following specific case to simplify the remainder of the base case.
    
    \begin{figure}
        \centering
        \includegraphics[page=2]{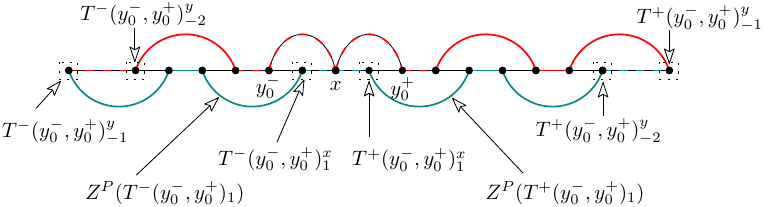}
        \caption{Illustration of Claim~\ref{thrm-S}.4. Here $T^-(y^-_0,y^+_0) = T^-(y^-_0,T^+(y^-_0,y^+_0)^y_{-1})$, $T^+(y^-_0,y^+_0) = T^+(y^-_0,T^+(y^-_0,y^+_0)^y_{-1})$ and the paths $Z^P(T^-(y^-_0,y^+_0)_1)$, $Z^P(T^+(y^-_0,y^+_0)_1)$ are disjoint. The solid green paths correspond to the odd path of $T^-(y^-_0,y^+_0)$ and the odd path of $T^+(y^-_0,y^+_0)$. The dashed green paths show these extended to $(x, T^-(y^-_0,y^+_0)^y_{-1})$ and $(x, T^+(y^-_0,y^+_0)^y_{-1})$ paths. The solid red paths correspond to the even path of $T^-(y^-_0,y^+_0)$ and the even path of $T^+(y^-_0,y^+_0)$. The dashed red paths show these extended to $(x, T^-(y^-_0,y^+_0)^y_{-1})$ and $(x, T^+(y^-_0,y^+_0)^y_{-1})$ paths.}
        \label{fig-base-switch}
    \end{figure}

\medskip
\noindent {\ensuremath{\vartriangleright}} {\sf \sffamily Claim~\ref{thrm-S}.4.}
        Suppose that $T^-(y^-_0,y^+_0) = T^-(y^-_0,T^+(y^-_0,y^+_0)^y_{-1})$ and $T^+(y^-_0,y^+_0) = T^+(T^-(y^-_0,y^+_0)^y_{-1}, y^+_0)$ and the paths $Z^P(T^-(y^-_0,y^+_0)_1)$, $Z^P(T^+(y^-_0,y^+_0)_1)$ are disjoint. Then either $G$ contains some $S_{2k,2k,2k}$ subgraph with centre $x$ or a $T$-type subgraph.
    \begin{claimproof}
        An illustration of the following is shown in Figure~\ref{fig-base-switch}. We note that as $T^+(y^-_0,y^+_0) = T^+(T^-(y^-_0,y^+_0)^y_{-1}, y^+_0)$ and these are maximal, there does not exist a positive jump out of the interval $(T^-(y^-_0,y^+_0)^y_{-1}, T^+(y^-_0,y^+_0)^y_{-1})$. Further, as $T^-(y^-_0,y^+_0) = T^-(y^-_0,T^+(y^-_0,y^+_0)^y_{-1})$ there does not exist a negative jump out of the interval $(T^-(y^-_0,y^+_0)^y_{-1}, T^+(y^-_0,y^+_0)^y_{-1})$.

        We now claim that the vertices $(T^-(y^-_0,y^+_0)^y_{-1}, T^+(y^-_0,y^+_0)^y_{-1})$ meet the criteria to apply Lemma~\ref{lem-cut-pair-S}. As there is no jump out of the interval $(T^-(y^-_0,y^+_0)^y_{-1}, T^+(y^-_0,y^+_0)^y_{-1})$, it follows that either $G$ contains some $S_{2k,2k,2k}$ subgraph with centre $x$ or a $T$-type subgraph.

        By Claim~\ref{thrm-S}.3, we get that either $G$ contains some $S_{2k,2k,2k}$ subgraph with centre $x$ or $\dist_P(x, T^+(y^-_0,y^+_0)^y_{-1}) < 4k^2$. Symmetrically, the same holds for $T^-(y^-_0,y^+_0)^y_{-1})$, that is, $T^-(y^-_0,y^+_0)^y_{-1}, T^+(y^-_0,y^+_0)^y_{-1} \in V(P[2k+1 : - (2k+1)])$. Condition~i) holds by maximality of $y^-_0$ and $y^+_0$.

        We note that both the odd and even paths of $T^+(y^-_0,y^+_0)$ and $T^-(y^-_0,y^+_0)$ are internally disjoint from the paths $Z^P(x, y^+_0)$ and $Z^P(x, y^-_0)$ else $y^-_0$ or $y^+_0$ weren't maximum. It follows that we can extend the odd path of $T^+(y^-_0,y^+_0)$ via the path $P[x:T^+(y^-_0,y^+_0)^x_1]$ and the even path of $T^+(y^-_0,y^+_0)$ via the path $Z^P(x, y^+_0)+ P[y^+_0:T^+(y^-_0,y^+_0)^x_2]$. We note that these are internally disjoint $(x, T^+(y^-_0,y^+_0)^y_{-1})$ paths. Note, symmetrically, we can define internally disjoint $(x, T^-(y^-_0,y^+_0)^y_{-1})$ paths. It follows that condition~ii) and~iii) hold.

        We now claim that the extended odd path of $T^-(y^-_0,y^+_0)$ is internally disjoint from both extended paths of $T^+(y^-_0,y^+_0)$. Suppose these paths are not internally disjoint. It follows that there exists some odd $i \in \{1,\ldots, |T^-(y^-_0,y^+_0)|\}$ and $j \in \{2,\ldots, |T^+(y^-_0,y^+_0)|\}$ such that the paths $Z^P(T^-(y^-_0,y^+_0)_i)$ and $Z^P(T^+(y^-_0,y^+_0)_j)$ intersect. It follows that there exists a jump between $T^-(y^-_0,y^+_0)^x_i$ and $T^+(y^-_0,y^+_0)^y_j$ hence $T^+(y^-_0,y^+_0)^y_1$ was not maximum. It follows that the extended odd path of $T^-(y^-_0,y^+_0)$ is internally disjoint from both extended paths of $T^+(y^-_0,y^+_0)$ and so condition~iv) holds. That is applying Lemma~\ref{lem-cut-pair-S}, we find that either $G$ contains some $S_{2k,2k,2k}$ subgraph with centre $x$ or a $T$-type subgraph, thus concluding our proof.
    \end{claimproof}

    \noindent
    Suppose that $T^+(y^-_0,y^+_0)^x_1 \in V(P[x : y^+_0]) \setminus \{y^+_0\}$, $c^x_1 \in V(P[y^-_0 : x]) \setminus \{y^-_0\}$ and the paths $Z^P(T^+(y^-_0,y^+_0)_1)$ and $Z^P(c^x_1, c^y_1)$ are disjoint. We note that this implies that the maximum negative chain extension of $(b^-_1,b^+_1)$ is also the maximum negative chain extension of $(y^-_0,y^+_0)$.
    
    In this case we will \textit{reverse} $P$, i.e. let $P = (a^+_\ell, \ldots, a^+_1, x, a^-_1, \ldots, a^-_\ell)$ and redefine vertices $y^-_0$, $y^+_0$, $b^-_1$, $b^+_1$, $c^x_1$ and $c^y_1$ based on this new path. 

    As the maximum negative chain extension of $(y^-_0,y^+_0)$ is also the maximum negative chain extension of $(y^-_0,T^+(y^-_0,y^+_0)^y_{-1})$, before reversing $P$, and $b^+_1$ becomes what was previously $T^-(y^-_0,y^+_0)^y_{-1}$, when $P$ is reversed. It follows that if, $T^+(y^-_0,y^+_0)^x_1 \in V(P[x : y^+_0]) \setminus \{y^+_0\}$, $c^x_1 \in V(P[y^-_0 : x]) \setminus \{y^-_0\}$ and the paths $Z^P(T^+(y^-_0,y^+_0)_1)$ and $Z^P(c^x_1, c^y_1)$ are disjoint, then the same must hold where $P$ is reversed. More formally, that is $T^-(y^-_0,T^+(y^-_0,y^+_0)^y_{-1})^x_1 \in V(P[y^-_0 : x]) \setminus \{y^-_0\}$ and $T^+(T^-(y^-_0,y^+_0)^y_{-1}, y^+_0)^x_1 \in V(P[x : y^+_0]) \setminus \{y^+_0\}$.

    As $Z^P(T^-(y^-_0,T^+(y^-_0,y^+_0)^y_{-1})_1)$, $Z^P(T^+(T^-(y^-_0,y^+_0)^y_{-1}, y^+_0)_1)$ are disjoint, it then follows by Claim~\ref{thrm-S}.4 that either $G$ contains some $S_{2k,2k,2k}$ subgraph with centre $x$ or a $T$-type subgraph. In both of these cases our Theorem holds.

    \medskip\noindent    
    \textbf{Base case: Define the jump sequences $I^1_1$, $I^2_1$, $I^3_1$ with properties P3a--f.}\\
    Recall that $b^-_1$ is active. Note we will define the jump sequences $I^1_1$, $I^2_1$, $I^3_1$ such that $I^1_1$ is the positive frontier, $I^2_1$ is the candidate and $I^3_1$ is the inert frontier.

    \medskip
    \noindent {\ensuremath{\vartriangleright}} {\sf \sffamily Claim~\ref{thrm-S}.5.}
        We can define jump sequences $I^1_1$, $I^2_1$, $I^3_1$ with the following properties.
        \begin{description}
            \item[B1.] $(I^1_1)^y_{-1} = c^y_1$, $(I^1_1)^x_{-1} \in V(P[b^-_{1}: b^+_{1}]) \setminus \{b^-_{1}, b^+_{1}\}$ and $b^-_{1} \in V(P[f^1_1:f^2_1]) \setminus \{f^1_1\}$.
            
            \item[B2.] The paths described by $I^1_1$, $I^2_1$, $I^3_1$ can be extended to $x$ via $P[x:s^1_1]$, $P[x:s^2_1]$ and $P[x:s^3_1]$ respectively. Let $D_1$, $D_2$ and $D_3$ denote these paths. Recall that for $i \in \{1,2,3\}$, if $I^i_1 = []$, then $s^i_1 = x$ and $D_i$ consists of the single vertex $x$. Else, $s^i_1 = (I^i_1)^x_1$.
                        
            \item[B3.] Let $\bar{D}_1$ be the subpath of $D_1$ from $x$ to $(I^a_\delta)^y_{-2}$. If $|I^a_\delta| =1$, then let $\bar{D}_1 = (x)$. The paths $\bar{D}_1$ and $D_3$ can both be extended to $b^-_\delta$. We denote these extended paths by $D^{\ext}_1$ and $D^{\ext}_3$ respectively.
                        
            \item[B4.] The paths $P[:b^-_1]$, $P[b^+_1:]$, $P[c^y_1: y_0^+]$, $D^{\ext}_1$, $D^{\ext}_3$ and $D_2$ are pairwise internally disjoint.
            
            \item[B5.] The path $Z^P((I^1_1)_{-1})$ is disjoint from both $D^{\ext}_3$ and $D_2$.
        \end{description}
    
    \begin{claimproof}
        We first note that, if the path $Z^P(c^x_1, c^y_1)$ intersects either $Z^P(y^-_0, x)$ or $Z^P(x, y^+_0)$ then the path $Z^P(x, c^y_1)$ exists. This is a contradiction as it implies that $y^-_0$ was not maximal. 
        
        We first consider the case where $c^x_1 \in V(P[y^-_0 : y^+_0]) \setminus \{y^+_0\}$. In this case, the path $Z^P(c^x_1, c^y_1)$ is internally disjoint from the even path of $T^+(y^-_0,y^+_0)$ else $T^+(y^-_0,y^+_0)^y_1$ was not maximal. Further, if $Z^P(c^x_1, c^y_1)$ is not disjoint from the odd path of $T^+(y^-_0,y^+_0)$ it may only intersect with vertices in $Z^P(T^+(y^-_0,y^+_0)_1)$, else again $T^+(y^-_0,y^+_0)^y_1$ was not maximal. 
        
        We now further consider the following subcases based on the position of $c^x_1$ and $T^+(y^-_0,y^+_0)^x_{1}$. That is, the cases of: $c^x_1$, $T^+(y^-_0,y^+_0)^x_{1} \in V(P[x : y^+_0])$; $c^x_1, T^+(y^-_0,y^+_0)^x_{1} \in V(P[y^-_0:x])$; or one of $c^x_1$ and $T^+(y^-_0,y^+_0)^x_{1}$ is in $V(P[y^-_0:x])$ and the other in $V(P[x : y^+_0])$. In the final of these cases, we will also differentiate between where $V(Z^P(c^x_1, c^y_1)) \cap V(Z^P(T^+(y^-_0,y^+_0)_1)) = \varnothing$ and where $V(Z^P(c^x_1, c^y_1)) \cap V(Z^P(T^+(y^-_0,y^+_0)_1)) \neq \varnothing$. 
        Each of these cases are illustrated in Figure~\ref{fig-base-intervals} with Table~\ref{tab-base-intervals} giving the specific jump sequences $I^1_1$, $I^2_1$, $I^3_1$, paths $D_1$, $D_2$, $D_3$ and chain extension $\hat{T}$ corresponding to each case.
        
        \begin{figure}
            \centering
            \includegraphics[page=3]{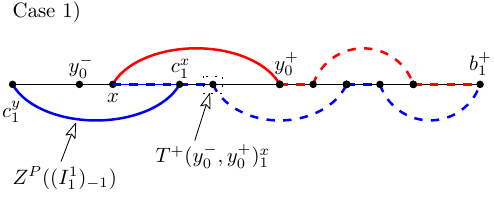}
            \includegraphics[page=4]{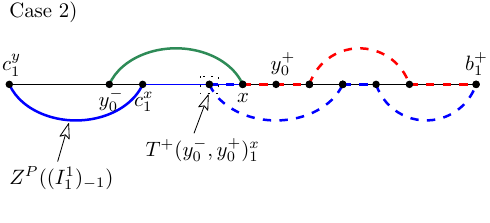}
            \includegraphics[page=5]{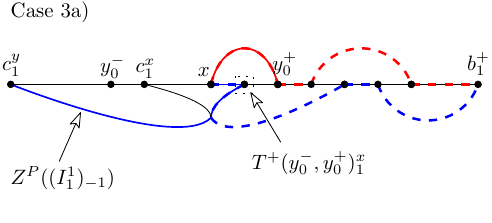}
            \includegraphics[page=6]{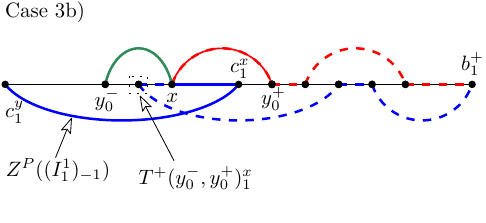}
            \includegraphics[page=7]{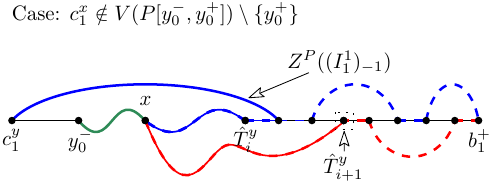}
            \caption{An illustration of the cases of Claim~\ref{thrm-S}.5.
            The path $D_2$ is drawn in green, $D_1$ is drawn in solid blue and $D_3$ is drawn in solid red. The path $D^{\ext}_1$ is drawn in dashed blue and $D^{\ext}_3$ is drawn in dashed red. Note where a subpath is contained in both $D_1$ and $D^{\ext}_1$ (likewise for $D_3$ and $D^{\ext}_3$) a dashed line is drawn on top of the solid line.
            See also Table~\ref{tab-base-intervals} for a summary of the jump sequences $I^1_1$, $I^2_1$, $I^3_1$, paths $D_1$, $D_2$, $D_3$ and chain extension $\hat{T}$ corresponding to each case.}
            \label{fig-base-intervals}
        \end{figure}
                
        \begin{itemize}
            \item \textbf{Subcase 1.} If $c^x_1$, $T^+(y^-_0,y^+_0)^x_{1} \in V(P[x : y^+_0])$, let $I^1_1= [(c^x_1, c^y_1)]$, $I^2_1= []$ and $I^3_1 = [(x,y^+_0)]$. Note, $V(P[f^1_1:f^2_1]) \setminus \{f^1_1\} = V(P[c^y_1: x]) \setminus \{c^y_1\}$ so by definition, $b^-_{1} \in V(P[f^1_1:f^2_1]) \setminus \{f^1_1\}$ and so Properties B1 and B2 hold by definition. Let $D_1$, $D_2$ and $D_3$ be those paths as defined in the claim statement.

            Let $\hat{T}$ denote the jump sequence $(x,y^+_0) + T^+(y^-_0,y^+_0)$. We note that $\hat{T}$ is the maximum positive chain extension of $(a^-_1,a^+_1)$. Further, $\hat{T}^x_2 = T^+(y_0^-, y_0^+)^x_1$ and $\hat{T}^y_{-1} = b_1^+$. We extend the even path of $\hat{T}$ by $P[x:\hat{T}^x_2]$. Observe that the odd path of $\hat{T}$ and this extended even path of $\hat{T}$ both begin with the vertex $x$. 
            
            We now show that these paths can be further extended such that they both end with the vertex $\hat{T}^y_{-1}$, recall that $b^+_1 = \hat{T}^y_{-1}$. If $|\hat{T}|$ is odd then we obtain these $(x,b^+_1)$ paths by extending the even path of $\hat{T}$ by $P[\hat{T}^y_{-2}:\hat{T}^y_{-1}]$, else if $|\hat{T}|$ is even, we extend the odd path of $\hat{T}$ by $P[\hat{T}^y_{-2}:\hat{T}^y_{-1}]$. We denote the extended odd path by $D^{\ext}_3$ and the extended even path by $D^{\ext}_1$. We note that $D^{\ext}_1$ and $D^{\ext}_3$ are internally disjoint.

            We note that $D_3 = Z^P(x,y^+_0)$ is a subpath of $D^{\ext}_3$ and $\bar{D}_1 = (x)$ is a subpath of $D^{\ext}_1$. That is, Property~B3 of our claim also holds. 
            
            As $\hat{T} = T^+(a^-_1,a^+_1)$ and $\hat{T}^y_{-1} = b^+_1$, it follows that $D^{\ext}_1$, $D^{\ext}_3$ and $D_2$ are internally disjoint from each other, $P[:x]$ and $P[b^+_1:]$. Hence, Property~B4 also holds. Finally, as $Z^P(c^x_1, c^y_1)$ is internally disjoint from the even path of $T^+(y^-_0,y^+_0)$, Property~B5 also holds.\\

            \item \textbf{Subcase 2.} If $c^x_1, T^+(y^-_0,y^+_0)^x_{1} \in V(P[y^-_0:x])$, let $I^1_1= [(c^x_1, c^y_1)]$, $I^2_1 = [(x,y^-_0)]$ and $I^3_1= []$. Again, by definition, Properties B1 and B2 of our claim hold.

            Let $\hat{T} = T^+(y^-_0,y^+_0)$. If $|\hat{T}| = 1$, then we let $D^{\ext}_1 = P[x:\hat{T}^y_{1}]$ and $D^{\ext}_3 = Z^P(\hat{T}_1)$. Otherwise, if $|\hat{T}| \geq 2$, then we extend the even path of $\hat{T}$ by $P[x:\hat{T}^x_2]$ and the odd path of $\hat{T}$ by $P[x:\hat{T}^x_1]$. In addition, if $|\hat{T}|$ is odd, then we extend the even path of $\hat{T}$ by $P[\hat{T}^y_{-2}:\hat{T}^y_{-1}]$, else we extend the odd path of $\hat{T}$ by $P[\hat{T}^y_{-2}:\hat{T}^y_{-1}]$. We denote the extended even path by $D^{\ext}_3$ and the extended odd path by $D^{\ext}_1$. We note that $D^{\ext}_1$ and $D^{\ext}_3$ are internally disjoint $(x,b_1^+)$ paths.

            As $D_3 = (x)$ is a subpath of $D^{\ext}_3$ and $\bar{D}_1 = (x)$ is a subpath of $D^{\ext}_1$, Property~B3 of our claim also holds. Further, by definition of these paths, Property B4 also holds. Now as $Z^P(c^x_1, c^y_1)$ is internally disjoint from the even path of $T^+(y^-_0,y^+_0)$, Property B5 holds.\\
                        
            \item \textbf{Subcase 3a.} Suppose one of $c^x_1$ and $T^+(y^-_0,y^+_0)^x_{1}$ is in $V(P[y^-_0:x])$ and the other in $V(P[x: y^+_0])$ and $V(Z^P(c^x_1, c^y_1)) \cap V(Z^P(T^+(y^-_0,y^+_0)_1)) \neq \varnothing$. Let $I^1_1= [(T^+(y^-_0,y^+_0)^x_1, c^y_1)]$, $I^2_1 = []$ and $I^3_1= [(x,y^+_0)]$. By definition, Properties~B1 and B2 hold.
            
            Recall that by definition $c^x_1$ minimises $\dist_P(c^x_1,c^y_1)$ and $T^+(y^-_0,y^+_0)^x_{1}$ minimises\\ $\dist_P(T^+(y^-_0,y^+_0)^x_{1},T^+(y^-_0,y^+_0)^y_{1})$. It follows that $c^x_1 \in V(P[y^-_0:x])$ and $T^+(y^-_0,y^+_0)^x_{1} \in V(P[x:y^+_0])$.

            Let $\hat{T} = (x,y^+_0) + T^+(y^-_0,y^+_0)$. The even path of $\hat{T}$ can be extended via the path $P[x:\hat{T}^x_2]$. In addition, if $|\hat{T}|$ is odd, then we extend the even path of $\hat{T}$ by $P[\hat{T}^y_{-2}:\hat{T}^y_{-1}]$, else we extend the odd path of $\hat{T}$ by $P[\hat{T}^y_{-2}:\hat{T}^y_{-1}]$. That is, we obtain a pair of internally disjoint $(x, b^+_1)$ paths. We will denote the extended odd path by $D^{\ext}_3$ and that extended even path by $D^{\ext}_1$. Given $D_3$ is a subpath of $D^{\ext}_3$ and $\bar{D}_1 = (x)$ is a subpath of $D^{\ext}_1$, Property~B3 of our claim also holds.
                        
            Note that $Z^P(T^+(y^-_0,y^+_0)^x_1, \ c^y_1)$ is internally disjoint from both $Z^P(x,y^-_0)$ and $Z^P(x,y^+_0)$, else $y^-_0$ was not maximal. Further, as $Z^P(T^+(y^-_0,y^+_0)^x_1, c^y_1)$ is internally disjoint from the even path of $T^+(y^-_0,y^+_0)$ it is also internally disjoint from $D^{\ext}_3$. That is properties B4 and B5 also hold.

            \item \textbf{Subcase 3b.} Again suppose one of $c^x_1$ and $T^+(y^-_0,y^+_0)^x_{1}$ is in $V(P[y^-_0:x])$ and the other in $V(P[x : y^+_0])$. We now consider the case where $V(Z^P(c^x_1, c^y_1)) \cap V(Z^P(T^+(y^-_0,y^+_0)_1)) = \varnothing$.

            Suppose first that $T^+(y^-_0,y^+_0)^x_1 \in V(P[x : y^+_0]) \setminus \{y^+_0\}$ and $c^x_1 \in V(P[y^-_0 : x]) \setminus \{y^-_0\}$. We recall that as $V(Z^P(c^x_1, c^y_1)) \cap V(Z^P(T^+(y^-_0,y^+_0)_1)) = \varnothing$, $T^+(y^-_0,y^+_0)^x_1 \in V(P[x : y^+_0]) \setminus \{y^+_0\}$ and $c^x_1 \in V(P[y^-_0 : x]) \setminus \{y^-_0\}$, we must have chosen to \textit{reverse} $P$. That following our observation resulting from Claim~\ref{thrm-S}.4, $T^-(y^-_0,T^+(y^-_0,y^+_0)^y_{-1})^x_1 \in V(P[y^-_0 : x])$, $T^+(T^-(y^-_0,y^+_0)^y_{-1}, y^+_0)^x_1 \in V(P[x : y^+_0])$ and so either $G$ contains some $S_{2k,2k,2k}$ subgraph with centre $x$ or a $T$-type subgraph. In both cases this implies that our Theorem holds, that is we assume that $T^+(y^-_0,y^+_0)^x_1 \in V(P[y^-_0 : x]) \setminus \{y^-_0\}$ and $c^x_1 \in V(P[x : y^+_0]) \setminus \{y^+_0\}$.
            
            Let $\hat{T} = (x,y^+_0) + T^+(y^-_0,y^+_0)$, $I^1_1= [(c^x_1, c^y_1)]$ and $I^2_1= [(x,y^-_0)]$, $I^3_1= [\hat{T}_2]$.

            We extend the even path of $\hat{T}$ by $P[x:\hat{T}^x_2]$. If $|\hat{T}|$ is odd then we extend the even path of $\hat{T}$ by $P[\hat{T}^y_{-2}:\hat{T}^y_{-1}]$, else we extend the odd path of $\hat{T}$ by $P[\hat{T}^y_{-2}:\hat{T}^y_{-1}]$. We denote the extended odd path by $D^{\ext}_3$ and the extended even path by $D^{\ext}_1$. We note that $D^{\ext}_1$ and $D^{\ext}_3$ are internally disjoint $(x,b_1^+)$ paths. Given $D_3$ is a subpath of $D^{\ext}_3$ and $\bar{D}_1 = (x)$ is a subpath of $D^{\ext}_1$, Property~B3 of our claim also holds. Now, Property B4 holds by definition. 
            
            Finally, Property 5 holds as $Z^P(c^x_1, c^y_1)$ is internally disjoint from the even path of $T^+(y^-_0,y^+_0)$ and so also $D^{\ext}_1$.
        \end{itemize}

    \begin{table}[t]
    \begin{center}
     \aboverulesep=0ex
     \belowrulesep=0ex
    \begin{tabular}{m{1.2em}|m{9.25em}|m{14.5em}|m{10em}}
    \toprule
     & Jump sequences & Paths $D_1$, $D_2$, $D_3$ & $\hat{T}$ \\\midrule
    1. & $I^1_1= [(c^x_1, c^y_1)]$ & $D_1 = P[x:c^x_1] + Z^P(c^x_1, c^y_1)$ & $\hat{T} = (x,y^+_0)+ T^+(y^-_0,y^+_0)$\\ 
    & $I^2_1= []$ & $D_2 = (x)$ &\\
    & $I^3_1 = [(x,y^+_0)]$ & $D_3 = Z^P(x, y^+_0)$ & \\\midrule
    2. & $I^1_1= [(c^x_1, c^y_1)]$ & $D_1 = P[x:c^x_1] + Z^P(c^x_1, c^y_1)$ & $\hat{T} = T^+(y^-_0,y^+_0)$\\ 
    & $I^2_1 = [(x,y^-_0)]$ & $D_2 = Z^P(x, y^-_0)$ &\\
    & $I^3_1= []$ & $D_3 = (x)$ & \\\midrule
    3a. & $I^1_1= [(T^+(y^-_0, y^+_0), c^y_1)]$ & $D_1 = P[x:c^x_1] + Z^P(T^+(y^-_0, y^+_0), c^y_1)$ & $\hat{T} = (x,y^+_0)+ T^+(y^-_0,y^+_0)$\\ 
    & $I^2_1= []$ & $D_2 = (x)$ & \\
    & $I^3_1 = [(x,y^+_0)]$ & $D_3 = Z^P(x, y^+_0)$ & \\\midrule
    3b. & $I^1_1= [(c^x_1, c^y_1)]$ & $D_1 = P[x:c^x_1] + Z^P(c^x_1, c^y_1)$& $\hat{T} = (x,y^+_0)+ T^+(y^-_0,y^+_0)$\\ 
    & $I^2_1= []$ & $D_2 = (x)$ & \\
    & $I^3_1 = [(x,y^+_0)]$ & $D_3 = Z^P(x, y^+_0)$ &\\\midrule
    \bottomrule
    \end{tabular}
    
    \caption{\label{tab-base-intervals}A summary of Subcases 1, 2, 3a, and 3b of Claim~\ref{thrm-S}.5. For each case we give the jump sequences $I^1_1$,$I^2_1$, $I^3_1$, the paths $D_1$, $D_2$, $D_3$ and the chain extension $\hat{T}$. A visual depiction of these cases can be seen in Figure~\ref{fig-base-intervals}.}
    \end{center}
    
    \end{table}

        \noindent
        In each of the subcases above, we defined $I_1^1, I_1^2, I_1^3$ satisfying the properties B1--B5. Hence, if $c_1^x \in V(P[y_0^- : y_0^+])\setminus \{ y_0^+\}$, then our claim holds.
        
        We may now assume that $c^x_1 \notin V(P[y^-_0,  y_0^+]) \setminus \{ y_0^+\}$, that is $c^x_1 \in P[\hat{T}^y_i:\hat{T}^y_{i+1}] \setminus \{\hat{T}^y_{i+1}\}$ for some $i \in \{1, \ldots, |\hat{T}|-1\}$. Let $T$ be that chain extension resulting from Claim~\ref{thrm-S}.2. Recall that the odd and the even path of $T$ can be extended to a pair of internally disjoint $(x, b^+_1)$ paths, via the path $P[x:T^x_1]$ and the paths $P[x:T^x_2]$, $P[T^y_{-2}:T^y_{-1}]$, if $|T| \geq 2$, and $P[x:T^y_{1}]$ if $|T| =1$. We will denote these paths by $D^{\ext}_1$ and $D^{\ext}_3$ but not yet fix which of $D^{\ext}_1$ and $D^{\ext}_3$ refers to the odd path and which refers to the even path. 
        
        From Claim~\ref{thrm-S}.2, either the path $P[y^-_0:x]$ is internally disjoint from both $D^{\ext}_1$ and $D^{\ext}_3$ or the path $Z^P(y^-_0, x)$ is internally disjoint from both $D^{\ext}_1$ and $D^{\ext}_3$. If $P[y^-_0:x]$ is internally disjoint from both $D^{\ext}_1$ and $D^{\ext}_3$, then let $I^2_1 = []$, else, let $I^2_1 = [(y^-_0, x)]$.
                
        If $i$ is odd, let $D^{\ext}_1$ be that odd extended path and $D^{\ext}_3$ be that even extended path. Let $I^3_1 = [\hat{T}_j: 1\leq j \leq i$, $j\mod 2 =1]$ and $I^1_1 = [\hat{T}_j: 1\leq j \leq i$, $j\mod 2 =0] + (c^x_1,c^y_1)$. If $i$ is even, let $D^{\ext}_1$ be that even extended path and $D^{\ext}_3$ be that odd extended path, $I^3_1 = [\hat{T}_j: 1\leq j \leq i$, $j\mod 2 =0]$ and $I^1_1 = [\hat{T}_j: 1\leq j \leq i$, $j\mod 2 =1] + (c^x_1,c^y_1)$. By definition Properties B1 and B2 hold. Further, the subpath of $D_1$ to $(I^a_\delta)^y_{-2}$ is also a subpath of $D^{\ext}_1$ and $D_3$ is a subpath of $D^{\ext}_3$, that is, Property B3 holds. As the paths $D^{\ext}_1$, $D^{\ext}_3$, $D_2$ are internally disjoint from each other and by definition also the paths $P[:b^-_1]$, $P[b^+_1:]$, $P[c^y_1: y_0^+]$, Property B4 also holds.

        Note, the path $Z^P(c^x_1, c^y_1)$ is internally disjoint from the even path of $T$, if $i$ is odd, and the odd path of $T$ otherwise, else, $T^+(y^-_0,y^+_0)^y_{i+1}$ was not maximal. $Z^P(c^x_1, c^y_1)$ is also disjoint from the path $D_2$, else $T^+(y^-_0,y^+_0)^y_{1}$ was not maximal, that is Property B5 holds thus concluding the proof of this claim.
    \end{claimproof}

    \noindent
    Let $I^1_1$, $I^2_1$, $I^3_1$ be those jump sequences resulting from Claim~\ref{thrm-S}.5. We let $I^1_1$ be the active frontier, $I^2_1$ be the candidate and $I^3_1$ be the inert frontier. We note that Property B1 implies that $|I^1_1| \geq 1$ and so Property P3a holds in the base case. Further, the Properties B1-B5 of Claim~\ref{thrm-S}.5 correspond to Properties P3b-f. That is the vertices $b^-_1$, $b^+_1$ and jump sequences $I^1_1$, $I^2_1$, $I^3_1$ satisfy Properties P1, P2a--d and P3a--f, thus concluding the base case.

    \medskip \noindent
    \textbf{Inductive case}\\
    We now suppose that for some $1 \leq \delta \leq 6k-1$ there exist the pairs $(b^-_1,b^+_1), \ldots, (b^{+}_{\delta},b^{-}_{\delta})$ and jump sequences $I^i_1 \subseteq \ldots \subseteq I^i_{\delta}$ for $i \in \{1,2,3\}$ satisfying Properties P1-P3. Let $a \neq b \neq c \in \{1,2,3\}$ be such that $I^a_\delta$ is the active frontier, $I^b_\delta$ is the inert frontier and $I^c_\delta$ is the candidate. We consider the case where $b^{+}_{\delta}$ is active, the case where  $b^{-}_{\delta}$ is active follows symmetrically. From Properties~P3a and P3c we have paths $D_a$, $D_b$ and $D_c$ with length at least $\lfloor \frac{\delta+2}{3} \rfloor$, $\lfloor \frac{\delta+1}{3} \rfloor$ and $\lfloor \frac{\delta}{3} \rfloor$ respectively. Further, combining properties P3d-f, it follows that the paths $D_a$, $D_b$, $D_c$ and $P[f^c_\delta:f^a_\delta]$ are internally disjoint. We also assume by Property P3d that there exist the paths $D^{\ext}_a$ and $D^{\ext}_b$ as described.

    We will now define vertices $b^+_{\delta+1}, b^-_{\delta+1}$ and construct jump sequences $I^1_{\delta+1}$, $I^2_{\delta+1}$, $I^3_{\delta+1}$ that will again satisfy the properties. Further, $I^b_{\delta+1}$ will be the candidate and one of $I^a_{\delta+1}$, $I^c_{\delta+1}$ will be the active frontier, the other will be the inert frontier.

Note that Property P1 holds by the base case. We say $b^-_{\delta+1}$ is active and $b^-_{\delta+1} = b^-_{\delta}$, as $b_\delta^+$ is active, Properties P2a and P2b hold by definition.  If $T^+(b^{-}_{\delta}, f^a_{\delta})$ exists, let $b^+_{\delta+1} = T^+(b^{-}_{\delta}, f^a_{\delta})^y_{-1}$, else let $b^+_{\delta+1} = f^a_{\delta}$. Either by definition or maximality of $T^+(b^{-}_{\delta}, f^a_{\delta})$, it follows that $y^+(b^-_{\delta+1}, b^+_{\delta+1})$ does not exist. That is, Property P2c holds for $\delta+1$. It remains to show that Property P2d is also satisfied. 

    \begin{figure}
         \centering
         \includegraphics[width=1\linewidth,page=14]{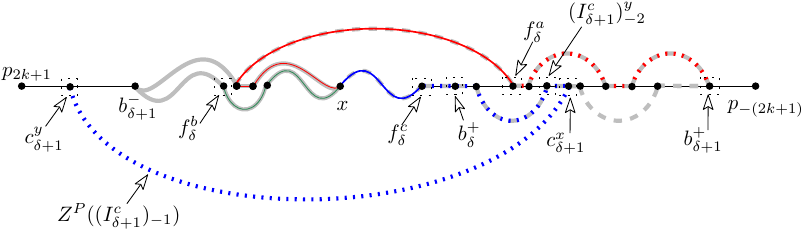}
         \caption{Here the paths $D_a$, $D_b$, $D_c$ are shown in \textbf{solid} red, green and blue respectively. The paths $D^{\ext}_a$ and $D^{\ext}_b$ from $x$ to $b^-_{\delta+1}$ are shown in \textbf{solid} gray. In the case depicted $I^c_{\delta+1}$ is the active frontier, $I^b_{\delta+1}$ is the candidate and $I^a_{\delta+1}$ is the inert frontier. The path $D^{ind}_a$ consists of both the solid and dotted red paths, the path $D^{ind}_c$ consist of the solid and dotted blue paths anf $D^{ind}_b= D_b$. The path $\hat{D}^{\ext}_a$ is that dashed gray path containing $D^{ind}_a$ as a subpath and the path $\hat{D}^{\ext}_c$ is that dashed gray path containing that subpath of $D^{ind}_c$ from $x$ to $(I^c_{\delta+1})^y_{-2}$.}
         \label{fig:induction}
     \end{figure}
     We will now construct a pair of internally disjoint paths from $x$ to $b^+_{\delta+1}$ which we will call, $\hat{D}^{\ext}_a$ and $\hat{D}^{\ext}_c$. $D_a$ and $D_c$ will form subpaths of $\hat{D}^{\ext}_a$ and $\hat{D}^{\ext}_c$, respectively. Further, the paths $\hat{D}^{\ext}_a$, $\hat{D}^{\ext}_c$, $D^{\ext}_b$, $P[:f^b_{\delta}]$ and $P[b^+_{\delta+1}:]$ are pairwise internally disjoint.

     We first consider the case where $T^+(b^{-}_{\delta}, f^a_{\delta})$ exists. By Property P3b, $f^a_{\delta} = y^+(b^-_{\delta},b^+_{\delta})$ and so $T^+(b^{-}_{\delta}, f^a_{\delta})^x_1 \in V(P[b^+_\delta:f^a_{\delta}]) \setminus \{f^a_{\delta}\}$. Further, both the odd and the even path of $T^+(b^{-}_{\delta}, f^a_{\delta})$ are disjoint from the paths $D^{\ext}_a$, $D^{\ext}_b$, $D_a$, $D_b$ and $D_c$ else $f^a_\delta$, i.e. $y^+(b^-_{\delta},b^+_{\delta})$, was not maximal.
     This implies that the path $D_a$ can be extended via $P[f^a_\delta: T^+(b^{-}_{\delta}, f^a_{\delta})^x_2]$ and the even path of $T^+(b^{-}_{\delta}, f^a_{\delta})$ and $D_c$ can be extended via $P[f^c_\delta: T^+(b^{-}_{\delta}, f^a_{\delta})^x_1]$ and the odd path of $T^+(b^{-}_{\delta}, f^a_{\delta})$, to obtain a pair of internally disjoint paths from $x$ to $b^+_{\delta+1}$ which are also disjoint from $D^{\ext}_b$, $P[:f^b_{\delta}]$ and $P[b^+_{\delta+1}:] - \{b^+_{\delta+1}\}$.
     We denote these paths by $\hat{D}^{\ext}_a$ and $\hat{D}^{\ext}_c$ respectively. See Figure~\ref{fig:induction} for an illustration.
     
     If $T^+(b^{-}_{\delta}, f^a_{\delta})$ does not exist, then we let $\hat{D}^{\ext}_a = D_a$ and $\hat{D}^{\ext}_c = D_c + P[f^c_{\delta}:f^a_{\delta}]$.
    \medskip
    \noindent {\ensuremath{\vartriangleright}} {\sf \sffamily Claim~\ref{thrm-S}.6.}
        If $\dist_P(b^+_{\delta}, b^+_{\delta+1}) \geq 4k^2-\frac{k(4\delta-14)}{3}+\frac{\delta(\delta-7)+11}{9}$, then $G$ contains some $S_{2k,2k,2k}$ subgraph with centre $x$.
    \begin{claimproof}
        Recall that there exist paths $D_a$, $D_b$ and $D_c$. We let $\ell_a$ denote the length of the path $D_a$ and $\ell_c$ denote the length of the path $D_c$.
    
        For every $i \in \{1,\ldots,|T^+(b^{-}_{\delta}, f^a_{\delta})|\}$, there exist paths from $x$ to $f^b_\delta$, $T^+(b^{-}_{\delta}, f^a_{\delta})^y_{i}$ and $T^+(b^{-}_{\delta}, f^a_{\delta})^y_{i+1}$ which share the single common vertex $x$. That path to $f^b_\delta$ corresponds to $D_b$. If $i$ is even, then that path to $T^+(b^{-}_{\delta}, f^a_{\delta})^y_{i}$ is a subpath of $\hat{D}^{\ext}_a$ and that path to $T^+(b^{-}_{\delta}, f^a_{\delta})^y_{i+1}$ is a subpath of $\hat{D}^{\ext}_c$. It follows that this path from $x$ to $T^+(b^{-}_{\delta}, f^a_{\delta})^y_{i}$ has length at least $\ell_a + \frac{i}{2}$ and that path from $x$ to $T^+(b^{-}_{\delta}, f^a_{\delta})^y_{i+1}$ has length at least $\ell_c + \frac{i}{2}+1$. Similarly, if $i$ is odd, then that path to $T^+(b^{-}_{\delta}, f^a_{\delta})^y_{i}$ is a subpath of $\hat{D}^{\ext}_c$ and that path to $T^+(b^{-}_{\delta}, f^a_{\delta})^y_{i+1}$ is a subpath of $\hat{D}^{\ext}_a$. These paths have length at least $\ell_c + \frac{i+1}{2}$ and $\ell_a + \frac{i+1}{2}$ respectively.
            
        Suppose $\dist_P(b^+_{\delta+1}, b^+_{\delta}) \geq 4k^2-\frac{k(4\delta-14)}{3}+\frac{\delta(\delta-7)+11}{9}$.
        We now claim that at least one of the following length or distance inequalities must hold,
        \begin{align*}
            |T^+(b^{-}_{\delta}, f^a_{\delta})^y_{-1}| > & \ 2(2k-\ell_c)-1\\
            |T^+(b^{-}_{\delta}, f^a_{\delta})^y_{-1}| > & \ 2(2k-\ell_a)\\
            \dist_P(f^c_\delta, f^a_\delta) > & \ 2k- \ell_c\\
            \dist_P(f^a_\delta, T^+(b^{-}_{\delta}, f^a_{\delta})^y_{1}) > & \ 2k- \ell_a\\
        %\end{align*}
        %\begin{align*}
            \dist_P(T^+(b^{-}_{\delta}, f^a_{\delta})^y_{i}, T^+(b^{-}_{\delta}, f^a_{\delta})^y_{i+1}) > & \ 2k - \left(\ell_c+\frac{i+1}{2}\right) \text{ for some }\\
            &  \textbf{ odd } i \in \{1,\ldots,|T^+(b^{-}_{\delta}, f^a_{\delta})^y_{-1}|-1\}, \text{ or}\\
            \dist_P(T^+(b^{-}_{\delta}, f^a_{\delta})^y_{i}, T^+(b^{-}_{\delta}, f^a_{\delta})^y_{i+1}) > & \ 2k - \left(\ell_a+\frac{i}{2}\right) \text{ for some }\\
            & \text{\textbf{ even} } i \in \{1,\ldots,|T^+(b^{-}_{\delta}, f^a_{\delta})^y_{-1}|-1\}.\\  
        \end{align*}
        \noindent
        Assume towards a contradiction that none of the above length or distance inequalities holds. This implies that for every $i \in \{1,\ldots,|T^+(b^{-}_{\delta}, f^a_{\delta})^y_{-1}|\}$,
        \begin{equation*}
          \dist(b^+_{\delta}, T^+(b^{-}_{\delta}, f^a_{\delta})^y_{i}) \leq \begin{cases}
            \sum^{\frac{i-1}{2}}_{j=0} (2k-(\ell_c+j)) + \sum^{\frac{i-1}{2}}_{j=0} (2k-(\ell_a+j)), \text{ \hspace{37pt}if $i$ is odd}. \\
            \\
            \dist(b^+_{\delta}, T^+(b^{-}_{\delta}, f^a_{\delta})^y_{i-1}) + 2k-1 -(\ell_a+\frac{i}{2}) =   \\
            \sum^{\frac{i-2}{2}}_{j=0} (2k-(\ell_c+j)) + \sum^{\frac{i-2}{2}}_{j=0} (2k-(\ell_a+j)) + 2k-1 -(\ell_a+\frac{i}{2}),  \\
             \text{ \hspace{226pt} if $i$ is even.}
          \end{cases}
        \end{equation*}

        \noindent
        Further, as $|T^+(b^{-}_{\delta}, f^a_{\delta})^y_{-1}| \leq 2(2k-\ell_c)-1$, $|T^+(b^{-}_{\delta}, f^a_{\delta})^y_{-1}| \leq 2(2k-\ell_a)$ and $b^+_{\delta+1} = T^+(b^{-}_{\delta}, f^a_{\delta})^y_{-1}$, it follows that,
        \begin{align*}
            \dist(b^+_{\delta}, b^+_{\delta+1}) & \leq \sum^{\frac{(2(2k-\ell_c)-2}{2}}_{j=0} (2k-(\ell_c+j)) + \sum^{\frac{(2(2k-\ell_c)-2}{2}}_{j=0} (2k-(\ell_a+j))\text{, and,}\\
            \dist(b^+_{\delta}, b^+_{\delta+1}) & \leq \sum^{\frac{2(2k-\ell_a)-2}{2}}_{j=0} (2k-(\ell_c+j)) \\
            & + \sum^{\frac{2(2k-\ell_a)-2}{2}}_{j=0} (2k-(\ell_a+j)) + 2k-1 -\left(\ell_a+\frac{2(2k-\ell_a)}{2}\right)\\
            & = \sum^{\frac{2(2k-\ell_a)-2}{2}}_{j=0} (2k-(\ell_c+j)) + \sum^{\frac{2(2k-\ell_a)-2}{2}}_{j=0} (2k-(\ell_a+j)) -1.
        \end{align*}

        \noindent
        By Property P3a and P3b, $\ell_a \geq \lfloor \frac{\delta+2}{3} \rfloor \geq \frac{\delta-1}{3}$ and $\ell_c \geq \lfloor \frac{\delta}{3} \rfloor \geq \frac{\delta-3}{3}$. It follows that $\dist(b^+_{\delta}, b^+_{\delta+1}) \leq 4k^2-\frac{k(4\delta-14)}{3}+\frac{\delta(\delta-7)+11}{9}-1$,
        a contradiction. That is we may assume that at least one of the length or distance conditions are met. 
        
        We will now show that in each of these cases $G$ contains some $S_{2k,2k,2k}$ subgraph with centre $x$.

        If $|T^+(b^{-}_{\delta}, f^a_{\delta})^y_{-1}| \geq 2(2k-\ell_c)$ then there exist paths from $x$ to $f^b_\delta$, $T^+(b^{-}_{\delta}, f^a_{\delta})^y_{2(2k-\ell_c)-1}$ and $T^+(b^{-}_{\delta}, f^a_{\delta})^y_{2(2k-\ell_c)}$ which share the single common vertex $x$. Further, that path from $x$ to $T^+(b^{-}_{\delta}, f^a_{\delta})^y_{2(2k-\ell_c)-1}$ has length at least $\ell_c + \frac{2(2k-\ell_c)-1+1}{2} =2k$, that is, extending that path from $x$ to $f^b_\delta$ by $P[:f^b_\delta]$ and that path from $x$ to $T^+(b^{-}_{\delta}, f^a_{\delta})^y_{2(2k-\ell_c)}$ by $P[T^+(b^{-}_{\delta}, f^a_{\delta})^y_{2(2k-\ell_c)}:]$, we obtain three paths of length at least~$2k$ with a single common vertex $x$. That is $G$ contains some $S_{2k,2k,2k}$ subgraph with centre $x$.

        Similarly, if $|T^+(b^{-}_{\delta}, f^a_{\delta})^y_{-1}| \geq 2(2k-\ell_a)+1$ then there exist paths from $x$ to $f^b_\delta$, $T^+(b^{-}_{\delta}, f^a_{\delta})^y_{2(2k-\ell_a)}$ and $T^+(b^{-}_{\delta}, f^a_{\delta})^y_{2(2k-\ell_a)+1}$ which share the single common vertex $x$. As that path from $x$ to $T^+(b^{-}_{\delta}, f^a_{\delta})^y_{2(2k-\ell_a)}$ has length at least $\ell_a + \frac{2(2k-\ell_a)}{2} = 2k$, extending that path from $x$ to $f^b_\delta$ by $P[:f^b_\delta]$ and that path from $x$ to $T^+(b^{-}_{\delta}, f^a_{\delta})^y_{2(2k-\ell_a)+1}$ by $P[T^+(b^{-}_{\delta}, f^a_{\delta})^y_{2(2k-\ell_a)+1}:]$, we obtain three paths of length at least~$2k$ with a single common vertex $x$. That is $G$ contains some $S_{2k,2k,2k}$ subgraph with centre $x$.

        If $\dist_P(f^c_\delta, f^a_\delta) \geq 2k+1- \ell_c$, then extending those paths from $x$ to $f^a_\delta$, $f^b_\delta$ and $f^c_\delta$ via the paths $P[f^a_\delta:]$, $P[:f^b_\delta]$ and $P[f^c_\delta:f^a_\delta] - \{f^a_\delta\}$, we again obtain three paths of length at least~$2k$ with a single common vertex $x$. That is $G$ contains some $S_{2k,2k,2k}$ subgraph with centre $x$. If $\dist_P(f^a_\delta, T^+(b^{-}_{\delta}, f^a_{\delta})^y_{1}) \geq 2k+1- \ell_a$, then similarly extending those paths from $x$ to $f^b_\delta$, $f^a_\delta$ and $T^+(b^{-}_{\delta}, f^a_{\delta})^y_{1}$ via the paths $P[:f^b_\delta]$, $P[f^a_\delta:T^+(b^{-}_{\delta}, f^a_{\delta})^y_{1}] - \{T^+(b^{-}_{\delta}, f^a_{\delta})^y_{1}\}$ and $P[T^+(b^{-}_{\delta}, f^a_{\delta})^y_{1}:]$, we obtain three paths of length at least~$2k$ with a single common vertex $x$. That is $G$ contains some $S_{2k,2k,2k}$ subgraph with centre $x$.

        Finally, suppose for some $i \in \{1,\ldots,|T^+(b^{-}_{\delta}, f^a_{\delta})^y_{-1}|-1\}$, we either have that $i$ is odd and $\dist_P(T^+(b^{-}_{\delta}, f^a_{\delta})^y_{i}, T^+(b^{-}_{\delta}, f^a_{\delta})^y_{i+1}) \geq 2k+1 - (\ell_c+\frac{i+1}{2})$ or we have that $i$ is even and $\dist_P(T^+(b^{-}_{\delta}, f^a_{\delta})^y_{i}, T^+(b^{-}_{\delta}, f^a_{\delta})^y_{i+1}) \geq 2k+1 - (\ell_a+\frac{i}{2})$. 
        We can once again extend those paths from $x$ to $f^b_\delta$, $T^+(b^{-}_{\delta}, f^a_{\delta})^y_{i}$ and $T^+(b^{-}_{\delta}, f^a_{\delta})^y_{i+1}$ via the paths $P[:f^b_\delta]$, $P[T^+(b^{-}_{\delta}, f^a_{\delta})^y_{i}:T^+(b^{-}_{\delta}, f^a_{\delta})^y_{i+1}] - \{T^+(b^{-}_{\delta}, f^a_{\delta})^y_{i+1}\}$ and $P[T^+(b^{-}_{\delta}, f^a_{\delta})^y_{i+1}:]$. Each of these paths have length at least~$2k$ with a single common vertex $x$. That is $G$ contains some $S_{2k,2k,2k}$ subgraph with centre $x$. 
    \end{claimproof}

    \noindent
    It now follows from Claim~\ref{thrm-S}.6, that Properties P2a--d hold for $b^-_{\delta+1}$ and $b^+_{\delta+1}$. We now claim that the vertices $b^-_{\delta+1}$, $b^+_{\delta+1}$ and $x$ meet the conditions for the application of Lemma~\ref{lem-cut-pair-S}. 
    
    Note that, Claim~\ref{thrm-S}.6 also implies $b^-_{\delta+1}, b^+_{\delta+1} \in V(P[2k+1:-(2k+1)])$. By maximality of the vertices $y^+_0$ and $y^-_0$, condition i) is satisfied. As $b^-_{\delta} = b^-_{\delta+1}$, by P3d there exist paths  $D^{\ext}_a$ and $D^{\ext}_b$ from $x$ to $b^-_{\delta+1}$ and by P3e these are internally disjoint, that is condition ii) is satisfied. There also exist internally disjoint paths $\hat{D}^{\ext}_a$ and $\hat{D}^{\ext}_c$ from $x$ to $b^+_{\delta+1}$. Further, the paths $D^{\ext}_a$, $D^{\ext}_b$ and $\hat{D}^{\ext}_c$ are also internally disjoint, that is conditions iii) and iv) are satisfied. Applying Lemma~\ref{lem-cut-pair-S} either $y^+(b^-_{\delta+1}, b^+_{\delta+1})$ or $y^-(b^-_{\delta+1}, b^+_{\delta+1})$ must exist. By the maximality of $b^+_{\delta+1}$, it follows that $y^+(b^-_{\delta+1}, b^+_{\delta+1})$ cannot exist and so $y^-(b^-_{\delta+1}, b^+_{\delta+1})$ must exist. Let $(c^x_{\delta+1}, c^y_{\delta+1}) = (x^-(b^-_{\delta+1}, b^+_{\delta+1}), y^-(b^-_{\delta+1}, b^+_{\delta+1}))$.
    
    By P2c, $y^-(b^-_{\delta}, b^+_{\delta})$ does not exist and $b^-_{\delta+1} = b^-_{\delta}$. It follows that $Z^P(c^x_{\delta+1}, c^y_{\delta+1})$ is internally disjoint from the paths $D_a$, $D_b$ and $D_c$ and $c^x_{\delta+1} \in V(P[b^+_{\delta}: b^+_{\delta+1}]) \setminus \{b^+_{\delta+1}\}$. Further, either $c^x_{\delta+1} \in V(P[b^+_{\delta}: f^a_{\delta}]) \setminus \{f^a_{\delta}\}$, $c^x_{\delta+1} \in V(P[f^a_{\delta}: T^+(b^{-}_{\delta}, f^a_{\delta})^y_{1}]) \setminus \{T^+(b^{-}_{\delta}, f^a_{\delta})^y_{1}\}$ or $c^x_{\delta+1} \in V(P[T^+(b^{-}_{\delta}, f^a_{\delta})^y_{i}: T^+(b^{-}_{\delta}, f^a_{\delta})^y_{i+1}]) \setminus T^+(b^{-}_{\delta}, f^a_{\delta})^y_{i+1}$ for some $i \in \{1, \ldots, |T^+(b^{-}_{\delta}, f^a_{\delta})|-1\}$.

    Suppose that $c^x_{\delta+1} \in V(P[b^+_{\delta}: f^a_{\delta}]) \setminus \{f^a_{\delta}\}$. It follows that the path $Z^P(c^x_{\delta+1}, c^y_{\delta+1})$ is internally disjoint from $\hat{D}^{\ext}_a$ and $\hat{D}^{\ext}_b \setminus Z^P(T^+(b^{-}_{\delta}, f^a_{\delta})_{1})$, else either $T^+(b^{-}_{\delta}, f^a_{\delta})^y_{1}$ was not maximal or $c^x_{\delta}$ was not minimal with respect to distance from $c^y_{\delta}$. We let $I^a_{\delta+1} = I^a_{\delta}$, $I^b_{\delta+1} = I^b_{\delta}$ and $I^c_{\delta+1} = I^c_{\delta} + (c^x_{\delta+1}, c^y_{\delta+1})$. We say $I^a_{\delta+1}$ is the inert frontier, $I^c_{\delta+1}$ is the active frontier and $I^b_{\delta+1}$ is the candidate.

    Suppose instead $c^x_{\delta+1} \in V(P[f^a_{\delta}: T^+(b^{-}_{\delta}, f^a_{\delta})^y_{1}]) \setminus \{T^+(b^{-}_{\delta}, f^a_{\delta})^y_{1}\}$. Similarly the path $Z^P(c^x_{\delta+1}, c^y_{\delta+1})$ is internally disjoint from $\hat{D}^{\ext}_a \setminus  Z^P(T^+(b^{-}_{\delta}, f^a_{\delta})_{2})$ and $\hat{D}^{\ext}_c$, else either $T^+(b^{-}_{\delta}, f^a_{\delta})^y_{2}$ was not maximal or $c^x_{\delta}$ was not minimal with respect to distance from $c^y_{\delta}$. We let $I^a_{\delta+1} = I^a_{\delta} + (c^x_{\delta+1}, c^y_{\delta+1})$, $I^b_{\delta+1} = I^b_{\delta}$ and $I^c_{\delta+1} = I^c_{\delta} + (T^+(b^{-}_{\delta}, f^a_{\delta})^x_{1}, T^+(b^{-}_{\delta}, f^a_{\delta})^y_{1})$. We say $I^a_{\delta+1}$ is the active frontier, $I^c_{\delta+1}$ is the inert frontier and $I^b_{\delta+1}$ is the candidate.

    Finally, suppose $c^x_{\delta+1} \in V(P[T^+(b^{-}_{\delta}, f^a_{\delta})^y_{i}: T^+(b^{-}_{\delta}, f^a_{\delta})^y_{i+1}]) \setminus T^+(b^{-}_{\delta}, f^a_{\delta})^y_{i+1}$ for some $i \in \{1, \ldots, |T^+(b^{-}_{\delta}, f^a_{\delta})|\}$. If $i$ is odd, then the path $Z^P(c^x_{\delta+1}, c^y_{\delta+1})$ is internally disjoint from $\hat{D}^{\ext}_a$ and $\hat{D}^{\ext}_c \setminus  Z^P(T^+(b^{-}_{\delta}, f^a_{\delta})_{i+1})$. 
    We let $I^a_{\delta+1} = I^a_{\delta} + [T^+(b^{-}_{\delta}, f^a_{\delta})_j: 1 \leq i \leq j,  \ j \mod 2 = 0]$, $I^b_{\delta+1} = I^b_{\delta}$ and $I^c_{\delta+1} = I^c_{\delta} + [T^+(b^{-}_{\delta}, f^a_{\delta})_j: 1 \leq i \leq j, \ j \mod 2 = 1] + (c^x_{\delta+1}, c^y_{\delta+1})$. 
    We say $I^a_{\delta+1}$ is the inert frontier, $I^c_{\delta+1}$ is the active frontier and $I^b_{\delta+1}$ is the candidate. 
    
    If $i$ is even, then the path $Z^P(c^x_{\delta+1}, c^y_{\delta+1})$ is internally disjoint from $\hat{D}^{\ext}_a \setminus  Z^P(T^+(b^{-}_{\delta}, f^a_{\delta})_{i+1})$ and $\hat{D}^{\ext}_c$. 
    We let $I^a_{\delta+1} = I^a_{\delta} + [T^+(b^{-}_{\delta}, f^a_{\delta})_j: 1 \leq i \leq j, j \mod 2 = 0] + (c^x_{\delta+1}, c^y_{\delta+1})$, $I^b_{\delta+1} = I^b_{\delta}$ and $I^c_{\delta+1} = I^c_{\delta} + [T^+(b^{-}_{\delta}, f^a_{\delta})_j: 1 \leq i \leq j, j \mod 2 = 1]$. We say $I^a_{\delta+1}$ is the active frontier, $I^c_{\delta+1}$ is the inert frontier and $I^b_{\delta+1}$ is the candidate.

    We note given $|I^c_{\delta+1}| = |I^c_{\delta}| +1$ Property P3a holds for $\delta+1$. Further, Property P3b holds by definition. Let us denote those paths described by $I^a_{\delta+1}$, $I^b_{\delta+1}$, $I^c_{\delta+1}$ extended to $x$ by $P[x:s^a_{\delta+1}]$, $P[x:s^b_{\delta+1}]$ and $P[x:s^c_{\delta+1}]$, respectively, by $D^{ind}_a$, $D^{ind}_b$, $D^{ind}_c$. Note the paths $D^{ind}_a \setminus (c^x_{\delta+1}, c^y_{\delta+1})$ and $D^{ind}_c \setminus (c^x_{\delta+1}, c^y_{\delta+1})$ form subpaths of $\hat{D}^{\ext}_a$ and $\hat{D}^{\ext}_c$. That is properties P3c-f also hold and concluding inductive case.

    We now assume that we have the jump sequences $I^1_{6k}$, $I^2_{6k}$ and $I^3_{6k}$. By properties P3a-f, there exist three paths of length at least~$2k$ from $x$, each sharing only the common vertex $x$. That is, $G$ contains a $S_{2k,2k,2k}$ with centre $x$.
\end{proof}

\section{The Algorithms}\label{sec-problem-specific}

We now apply Theorems~\ref{thrm-H} and~\ref{thrm-S} to show that \col{} is polynomial-time solvable
for $\mathbb{H}_m^{1,k,k,k}$-subgraph-free graphs and $S_{1,k,k,k}$-subgraph-free graphs. 
In both cases we need another structural result that allows us to assume that the graphs in the input do not contain protected fans, $T$-type and $L$-type subgraphs.
We first consider $\mathbb{H}_m^{1,k,k,k}$-subgraph-free graphs.

\begin{lemma}\label{lem-col-H}
    Let $(G,r)$ be an instance of \col\ and let $c,m,k \geq 1$ be constants. We can in polynomial time obtain a graph $G'$ such that:
    \begin{itemize}
        \item $G'$ has minimum degree at least~$3$;
        \item $G'$ contains no protected fan of order $m+k+2$;
        \item $G'$ contains no $L$-type subgraph, with respect to the bound $c$ and length $m+k$;
        \item $G'$ can be coloured with $r$ colours, if and only if $G$ can;
        \item if $G$ is $\mathbb{H}_m^{1,k,k,k}$-subgraph-free, then so is $G'$.
    \end{itemize}
\end{lemma}

\begin{proof}
    Towards this we will define three operations. The first we call degree~$2$ vertex removal, the second we call the fan-contraction operation and the third is the $L$ removal operation. 
    We will first describe each of these operations, showing that the resulting graph can be coloured using $r$ colours if, and only if, the original graph could. 
    We then also claim that if the original graph was $\mathbb{H}_m^{1,k,k,k}$-subgraph-free then so is the resulting graph. Whenever we refer to a $L$-type subgraph, it is always understood to be with respect to this pair of constants $c$ and $m+k$.

    \smallskip\noindent
    \textbf{Degree $2$ vertex removal:}\\    
    If $|V(G)| \geq 2$, $r\geq 3$ and $G$ contains some degree~$2$ vertex, $v$, then $(G,r)$ is a yes-instance of \col{} if, and only if, $(G-v,r)$ is a yes-instance. That is we first exhaustively remove degree~$2$ vertices. We highlight that the family of $\mathbb{H}_m^{1,k,k,k}$-subgraph-free graphs is closed under vertex deletion.

    \smallskip\noindent
    \textbf{Fan-contraction operation:}\\    
    Suppose some $Q \subseteq V(G)$ is a protected fan with centre vertex $z$ and end vertices $x,y$. We note if $r \geq 4$, then $G$ can be coloured with $r$ colours if and only if $G - (Q\setminus \{x,y,z\})$ can. 
    We let $G^F(Q) = G - (Q\setminus \{x,y,z\})$, and we call this a \emph{type $1$ fan-contraction}. Assume now $r=3$. 
    If $m+k+2$ is odd then $x,y,z$ must each take different colours in every $3$-colouring of $G$. Let $G^F(Q)$ be the graph obtained from $G$ by deleting the vertices $Q \setminus \{x,y,z\}$ and adding the edge $xy$. That is, the vertices $\{x,y,z\}$ induce a triangle in $G^F(Q)$. We call this a \emph{type $2$ fan-contraction}. 
    If $m+k+2$ is even, then $x$ and $y$ must take the same colour. We let $G^F(Q)$ be the graph obtained from $G$ by replacing the vertices $Q\setminus \{z\}$ by a new vertex $w$ such that $N(w) = N(x) \cup N(y)$. We call this a \emph{type $3$ fan-contraction}. We note this operation takes polynomial time to apply and $|V(G^F(Q))| \leq |V(G)|-1$.

    \medskip
\noindent {\ensuremath{\vartriangleright}} {\sf \sffamily Claim~\ref{lem-col-H}.1.}
        Let $G^F(Q)$ be that graph obtained from $G$ by applying the fan-contraction operation to some protected fan on vertices $Q \subseteq V(G)$, with $|Q| = m+k+2$. If $G^F(Q)$ contains some $\mathbb{H}_m^{1,k,k,k}$ subgraph then so must $G$.
    \begin{claimproof}
        Let $z,x,y \in Q$ denote the centre and end vertices of this fan, respectively. We let $R=(r_1,\ldots, r_{m+k+1})$ denote that path of length $m+k$ between $x$ and $y$ in $G[Q]-z$, we highlight that $r_1 = x$ and $r_{m+k+1} = y$.
                
        Suppose $G^F(Q)$ contains some $\mathbb{H}_m^{1,k,k,k}$ subgraph $F$. Let $f,f'$ denote the degree~$3$ vertices of $F$, let $f+ M + f'$ denote that path of length $m$ between $f$ and $f'$. Let $f^1_1$ be that isolated vertex in $F- (V(M) \cup \{f,f'\})$ and $F_2$, $F_3$, $F_4$ denote those paths of length $k-1$. For each $i \in \{2,3,4\}$ and $j \in \{1,\ldots,k\}$, let $f^i_j$ denote the $j$th vertex of $F_i$. Without loss of generality $f^1_1$ and $f^2_1$ are adjacent to $f$ and $f^3_1$ and $f^4_1$ are adjacent to $f'$. See Figure~\ref{fig-Hstar}.

        If $G^F(Q)$ has been obtained by a type $1$ fan-contraction, then $F \subseteq G$ and so $G$ contains $\mathbb{H}_m^{1,k,k,k}$ as a subgraph. Likewise, if $G^F(Q)$ has been obtained via a type $2$ operation but the edge $xy \notin E(F)$ or if $G^F(Q)$ has been obtained via a type $3$ operation but that vertex $w \in V(G^F(Q)) \setminus V(G)$ is not in $F$, then $F \subseteq G$. That is in each case we find $\mathbb{H}_m^{1,k,k,k}$ as a subgraph. We are now left with two cases: $G^F(Q)$ has been obtained via a type $2$ operation and $F$ contains the edge $xy$ or $G^F(Q)$ has been obtained via a type $3$ operation $F$ contains that new vertex $w \notin V(G)$.

        Suppose $G^F(Q)$ has been obtained via a type $2$ fan-contraction. If both $x,y \in V(F) \setminus M$, then without loss of generality either $x = f$ and
        $y \in \{f^1_1, f^2_1\}$;
        or $x = f'$ and $y \in \{f^3_1, f^4_1\}$; or $x = f^i_j$ and $y = f^i_{j+1}$, for some $i \in \{2,3,4\}$, $j \in \{1,\ldots,k-1\}$. Let $F'$ be that graph obtained by replacing the edge $xy \in E(F)$ by the path $(r_1,\ldots, r_{m+k+1})$. Note, $F'$ contains $\mathbb{H}_m^{1,k,k,k}$ as a subgraph and $F' \subseteq G$, that is $G$ contains $\mathbb{H}_m^{1,k,k,k}$ as a subgraph.

        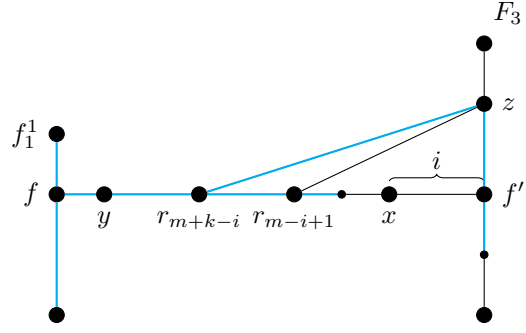
\begin{figure}
    \centering
    \begin{tikzpicture}[x=1.25cm, y=0.8cm]
        \node[draw, circle, fill, inner sep=2pt] (A1) at (0.5,0) [label={left:$f^1_1$}]{};
        \node[draw, circle, fill, inner sep=2pt] (y) at (1,-1) [label={below:$y$}]{};
        \node[draw, circle, fill, inner sep=2pt] (r) at (2,-1) [label={below:$r_{m+k-i}$}]{};
        \node[draw, circle, fill, inner sep=2pt] (r') at (3,-1) [label={below:$r_{m-i+1}$}]{};
        \node[draw, circle, fill, inner sep=1pt] (n) at (3.5,-1){};
        \node[draw, circle, fill, inner sep=2pt] (x) at (4,-1) [label={below:$x$}]{};
        
        \node[draw, circle, fill, inner sep=2pt] (B1) at (0.5,-3) {};
        \node[draw, circle, fill, inner sep=2pt] (C) at (5,1.5) {};
        \node[draw, circle, fill, inner sep=2pt] (z) at (5,0.5) [label={right:$z$}]{};
        \node[draw, circle, fill, inner sep=2pt] (E) at (0.5,-1) [label={left:$f$}]{};
        \node[draw, circle, fill, inner sep=2pt] (F) at (5,-1) [label={right:$f'$}]{};
        \node[draw, circle, fill, inner sep=1pt] (n') at (5,-2){};
        \node[draw, circle, fill, inner sep=2pt] (F1) at (5,-3) {};
        \node[] (lB) at (5.25,2) {$F_3$};

        \draw[decoration={brace,raise=5pt},decorate] (4,-1) -- node[above=6pt] {$i$} (5,-1);
    
        \draw (C) -- (F1);
        \draw (E) -- (F);
        \draw (r) -- (z);
        \draw (r') -- (z);

        \draw [cyan, thick] (A1) -- (E) -- (B1);
        \draw [cyan, thick] (E) -- (y);
        \draw [cyan, thick] (y) -- (r);
        \draw [cyan, thick] (r) -- (r');
        \draw [cyan, thick] (r') -- (n);
        \draw [cyan, thick] (r) -- (z);
        \draw [cyan, thick] (z) -- (F);
        \draw [cyan, thick] (F) -- (n');
        \end{tikzpicture}
        \caption{An illustration of the graph $F'$ from Claim~\ref{lem-col-H}.1, where $G^F(Q)$ has been obtained via a type $2$ fan-contraction and $z \in V(F_3)$. That $\mathbb{H}_m^{1,k,k,k}$ subgraph of $F'$ is highlighted in blue.}
        \label{fig:fan-safety-claim}
    \end{figure}

        We now assume that at least one of $x$ or $y \in M$. Recall that the case where $(x,y) \notin E(F)$ was covered previously, that is we assume that $x$ and $y$ appear consecutively in the path $(f) + M + (f')$. Without loss of generality, $\dist_F(x,f') \leq \dist_F(y,f')$ and $\dist_F(x,f') = i$, for some $i \in \{0, \ldots, m-1\}$.
        Let $F'$ be that graph obtained from $F$ by replacing the edge $xy \in E(F)$ by the path $(r_1,\ldots, r_{m+k+1})$ and adding the edges $zr_{m-i+1}$ and $zr_{m+k-i}$. 
        We note that $\dist_F(f',r_{m-i+1}) = m$, $\dist_F(f,r_{m+k-i}) = m$
        and $F' \subseteq G$. We now claim that $F'$ contains $\mathbb{H}_m^{1,k,k,k}$ as a subgraph. 
        If $z \notin V(M[1:i-1]) \cup F_3 \cup F_4 \cup \{f'\}$, then $F'$ contains $\mathbb{H}_m^{1,k,k,k}$ as a subgraph with degree~$3$ vertices $f'$, $r_{m-i+1}$, a path of length $1$ via the edge $zr_{m-i+1}$ and paths of length $k$ via $R[m-i+1:m+k-i+1]$, $f' + F_3$ and $f' + F_4$. 
        If $z \in V(M[1:i-1]) \cup F_3 \cup F_4 \cup \{f'\}$, then we note there exists some path $Z$ of length at least $k-1$ from $z$ in $F[V(M[1:i-1]) \cup F_3 \cup F_4 \cup \{f'\}]$. 
        It follows that $F'$ contains $\mathbb{H}_m^{1,k,k,k}$ as a subgraph with degree~$3$ vertices $f$ and $r_{m+k-i}$, a path of length $1$ via the edge $f f^1_1$ and paths of length $k$ via $f + F_2$, $f'+ R[m+k-(i+1):m-i]$, and $f' + Z$. 
        See Figure~\ref{fig:fan-safety-claim} for an illustration where $z \in V(F_3)$. That is, if $G^F(Q)$ has been obtained via a type $2$ fan-contraction then our claim holds.

        The case where $G^F(Q)$ has been obtained via a type $3$ fan-contraction follows similarly. Recall that we may assume that $F$ contains some vertex $w \notin V(G)$, else $F$ is also a subgraph of $G$. Further, $N(w) = N(x) \cup N(y)$, meaning if either $N(w) \cap V(F) \subseteq N(x)$ or $N(w) \cap V(F) \subseteq N(x)$, then $F$ is again a subgraph of $G$. 
        That is $G$ contains $\mathbb{H}_m^{1,k,k,k}$ as a subgraph. We therefore assume that there is some pair $r',r''$ such that $r' \in N(x) \cap V(F)$ and $r'' \in N(x) \cap V(F)$. If $w \notin V(M) \cup \{f,f'\}$, then let $F'$ be that graph obtained by replacing the path $(r',w,r'') \subseteq F$ by the path $(r', r_1,\ldots, r_{m+k+1},r'')$. 
        Note, $F'$ contains $\mathbb{H}_m^{1,k,k,k}$ as a subgraph and $F' \subseteq G$, that is $G$ contains $\mathbb{H}_m^{1,k,k,k}$ as a subgraph. 
        
        If $w \in V(M) \cup \{f,f'\}$, then $\dist_F(w,f') = i$ for some $i \in \{0, \ldots, m\}$. Let $F'$ be that graph obtained from $F$ by replacing the vertex $w$ with the path $(r_1,\ldots, r_{m+k+1})$ and adding the edges $zr_{m-i+1}$ and $zr_{m+k+1-i}$. Note that $\dist(f', r_{m-i+1}) = m$, $\dist(f, r_{m+k+1-i}) = m$ and $F' \subseteq G$.
        If $z \notin V(M[1:i-1]) \cup F_3 \cup F_4 \cup \{f'\}$, then $F'$ contains $\mathbb{H}_m^{1,k,k,k}$ as a subgraph with degree~$3$ vertices $f'$ and $r_{m-i+1}$, a path of length $1$ via the edge $zr_{m-i+1}$ and paths of length~$k$ via $R[m-i+1:m+k-i+1]$, $F_3$, and $F_4$. If $z \in V(M[1:i-1]) \cup F_3 \cup F_4 \cup \{f'\}$, then we note there exists some path $Z$ of length at least $k-1$ from $z$ in $F[V(M[1:i-1]) \cup F_3 \cup F_4 \cup \{f'\}]$. It follows that $F'$ contains $\mathbb{H}_m^{1,k,k,k}$ as a subgraph with degree~$3$ vertices $f$ and $r_{m+k+1-i}$, a path of length $1$ via the edge $ff^1_1$ and paths of length $k$ via $f+F_2$, $R[m+k+1-i:m+1-i]$, and $r_{m+k+1-i} + Z$.
    \end{claimproof}

    \smallskip\noindent
    \textbf{$L$ removal operation:}\\
    Suppose there is some $C \subseteq V(G)$ such that $G[C]$ is a minimal $L$-type subgraph. By Lemma~\ref{lem-T-type-smalltd}, $G[C \cup S]$ has treedepth at most~$3k+2$. That is, by Courcelle's theorem~\cite{Co90} we can decide, in polynomial time, if $G[C \cup S]$ can be coloured using $r$ colours. If not, then $(G,r)$ is a no-instance of \col. Assume now $G[C \cup S]$ can be coloured using $r$ colours. 

    Let $S = \{u, v\}$. We branch on the two non-isomorphic colourings of $u$ and $v$. For each branch, we consider the {\sc Precolouring Extension} problem on the graph $G[C \cup \{u, v\}]$, with precoloured vertices $u$ and $v$. Since $\td(G[C \cup {u, v}]) \leq 3k+2$, it follows from \cite{JS97}, who showed that even the more general problem {\sc List Colouring} is polynomial-time solvable on graphs of bounded treewidth, that each branch can be processed in polynomial time.

    If there is both a proper colouring where $u,v$ are coloured the same and where they are different, then $G$ can be coloured with $r$ colours if and only if, $G-C$ can be coloured. Let $G^T(C \cup \{u,v\}) = G-C$. This is again a type $1$ $L$-removal operation.
    
    If there is only a proper colouring of $G[C \cup \{u,v\}]$ where $u,v$ are coloured differently, then let $G^T(C \cup \{u,v\})$ be the graph obtained from $G$ by removing the vertices of $C$ and adding the edge $uv$. We call this a type $2$ operation. The remaining case is were $u$ and $v$ must take the same colour in every colouring of $G[C \cup S]$. Let $G^T(C \cup S)$ be the graph obtained from $G$ by removing the vertices of $C \cup \{u,v\}$ and adding a new vertex $w$ such that $N(w) = N(u) \cup N(v)$. We call this a type $3$ operation. Note, $G$ can be coloured using $r$ colours if, and only if, $G^T(C \cup S)$ can be coloured using $r$ colours.

    \medskip
\noindent {\ensuremath{\vartriangleright}} {\sf \sffamily Claim~\ref{lem-col-H}.2.}
        Let $G^T(C \cup S)$ be that graph obtained from $G$ by applying the $L$ removal operation to minimal $L$-type subgraph $G[C]$ with a witness set $S = \{u,v\}$. If $G^T(C \cup S)$ contains some $\mathbb{H}_m^{1,k,k,k}$ subgraph then so must~$G$.
    \begin{claimproof}
in        Note that $G^T(C \cup S)$ has been obtained via either a type $1$, $2$, $3$ operation. If it was a type $1$ operation, then $G^T(C \cup S)$ is a subgraph of $G$, that is our claim holds by definition. We therefore assume that $G^T(C \cup S)$ has been obtained via either an operation of type $2$ or~$3$.
        
        By definition, for some $\ell \geq m+k+1$, there is some induced path $R = (r_1,\ldots, r_\ell)$ where $u = r_1$ and $v = r_\ell$ in $G[C \cup S]$. Further, $u$ and $v$ have at least~$2$ neighbours in $C \cup S$. As $G$ has minimum degree at least~$3$ and $N(C) \setminus C = \{u,v\}$, for every $i \in \{1,\ldots, \ell\}$, $r_i$ has some neighbour in $C \setminus V(R)$ which we will denote by $z_i$. We highlight that $\{z_1,\ldots, z_\ell\} \cap V(G^T(C \cup S)) = \varnothing$.
        
        Suppose $G^T(C \cup S)$ contains some $\mathbb{H}_m^{1,k,k,k}$ subgraph $F$. Let $f,f'$ denote the degree~$3$ vertices of $F$, let $f+ M + f'$ denote that path of length $m$ between $f$ and $f'$. Let $f^1_1$ be that isolated vertex in $F- (V(M) \cup \{f,f'\})$ and $F_2$, $F_3$, $F_4$ denote those paths of length $k-1$. For each $i \in \{2,3,4\}$ and $j \in \{1,\ldots,k\}$, let $f^i_j$ denote the $j$th vertex of $F_i$. Without loss of generality $f^1_1$ and $f^2_1$ are adjacent to $f$ and $f^3_1$ and $f^4_1$ are adjacent to $f'$. See Figure~\ref{fig-Hstar}.

        If $G^T(C \cup S)$ has been obtained via a type $2$ operation but the edge $uv \notin E(F)$, then $F \subseteq G$. Likewise, if $G^T(C \cup S)$ has been obtained via a type $3$ operation but that vertex $w \in V(G^T(C \cup S)) \setminus V(G)$ is not in $F$, then $F \subseteq G$.
        
        If $G^T(C \cup S)$ has been obtained via a type $2$ operation, let $F'$ be that graph obtained from $F$ by replacing the edge $uv$ with the path $R$. If $G^T(C \cup S)$ has been obtained via a type $3$ operation, let $F'$ be that graph obtained from $F$ by replacing the vertex $w$ with the path $R$. Note in each case $F' \subseteq G$. If $u,v,w \in V(F) \setminus M$, then $F'$ and so also $G$ contains $\mathbb{H}_m^{1,k,k,k}$ as a subgraph.

       Suppose now that at least one of $u$, $v$ or $w \in M$. If $G^T(C \cup S)$ has been obtained via a type~$2$ operation, without loss of generality $\dist_F(u,f') \leq \dist_F(v,f')$ and $\dist_F(u,f') = i$ for some $i \in \{0, \ldots, m-1\}$. If the operation had type $3$, then $\dist_F(w,f') = i$ for some $i \in \{0, \ldots, m\}$. Let $F'$ be that graph obtained from $F$ by replacing the edge $uv \in E(F)$ or vertex $x$ by the path $(r_1,\ldots, r_{m+k+1})$ and adding the edge $z_{m-i+1}r_{m-i+1}$. We note that $\dist_R(f',r_{m-i+1})=m$. Now $F'$ contains $\mathbb{H}_m^{1,k,k,k}$ as a subgraph with degree~$3$ vertices $f'$, $r_{m-i+1}$, a path of length $1$ via the edge $z_{m-i+1}r_{m-i+1}$ and paths of length $k$ via $R[m-i+1:m+k-i+1]$, $f' + F_3$ and $f' + F_4$.
    \end{claimproof}

    \noindent
    We have now shown that, after applying one of these operations, we obtain a resulting graph which can be coloured using $r$ colours if, and only if, the original graph could. Further, if the original graph was $\mathbb{H}_m^{1,k,k,k}$-subgraph-free then so is the resulting graph. We will now show how these operations can be applied to obtain that graph $G'$ described in the theorem statement.
    
    We first, in polynomial time, exhaustively apply the degree~$2$ removal operation to obtain a graph with minimum degree at least~$3$. Further, in polynomial time, for every $Q \subseteq V(G)$ of size $m+k+2$ we can decide if $G[Q]$ is a protected fan in $G$. 
    If $G[Q]$ is a protected fan we can construct a graph $G^F(Q)$, such that $|V(G^F(Q))| < |V(G)|$, $(G^F(Q),r)$ is a yes-instance of \col{} if, and only if, $(G,r)$ is. Further, if $G$ is $\mathbb{H}_m^{1,k,k,k}$-subgraph-free, then so is $G^F(Q)$.

    Likewise, in polynomial time, for every $S \subseteq V(G)$, such that $|S| \leq 2$, we can decide if $S$ is the witness set for some minimal $L$-type subgraph. If so we can construct a graph $G^T(C \cup S)$ such that $|G^T(C \cup S)| < |V(G)|$, $(G^T(C \cup S),r)$ is a yes-instance of \col{} if and only if $(G,r)$ is and if $G$ is $\mathbb{H}_m^{1,k,k,k}$-subgraph-free, then so is $G^T(C \cup S)$.
    
    We will first apply the degree~$2$ removal operation. We will then alternate between exhaustively applying the fan contraction operation and applying the $L$ removal operation. After each application of either the fan contraction operation or the $L$ removal operation we again remove any degree~$2$ vertices. Doing this exhaustively, we obtain a graph $G'$ such that $G'$ contains no protected fan of order $m+k+2$, $G'$ contains $L$-type subgraph, $G'$ can be coloured with $r$ colours, if and only if $G$ can, and if $G$ is $\mathbb{H}_m^{1,k,k,k}$-subgraph-free, then so is~$G'$.
\end{proof}

\noindent
Given an instance $(G,r)$ of \col, we now apply Lemma~\ref{lem-col-H} and obtain a modified instance $(G',r)$. By Theorem~\ref{thrm-H}, the treedepth of $G'$ is bounded by a constant and since \col\ is polynomial-time solvable for graphs of bounded treewidth~\cite{KR03} and therefore also for graphs of bounded treedepth we obtain the following.

\begin{theorem}\label{thm-col-H}
    \col\ is solvable in polynomial time for $\mathbb{H}_m^{1,k,k,k}$-subgraph-free graphs, for all $m, k \geq 1$.
\end{theorem}

\begin{proof}
    Let $(G,r)$ be an instance of \col, where $G$ is a $\mathbb{H}_m^{1,k,k,k}$-subgraph-free graph.
    We apply Lemma~\ref{lem-col-H} and obtain in polynomial time an instance $(G',r)$ of \col\ such that $G'$ has minimum degree at least~$3$; $G'$ contains no protected fan of order $m+k+2$; and $G'$ contains no $L$-type subgraph, with treedepth bound $(4k - 3)(8k^2 - 6k + 2m + 8) - 1$ and length bound $m+k$. Further, $G'$ can be coloured with $r$ colours, if and only if $G$ can and as $G$ is $\mathbb{H}_m^{1,k,k,k}$-subgraph-free, then so is $G'$.

    By Theorem~\ref{thrm-H}, the treedepth of $G'$ is bounded by a constant. As \col\ is polynomial-time solvable for graphs of bounded treewidth~\cite{KR03} and therefore also for graphs of bounded treedepth the theorem follows.
\end{proof}

\noindent
We continue with the application of Theorem~\ref{thrm-S} to solve \col{} on $S_{1,k,k,k}$-subgraph-free graphs. 
Similarly to Theorem~\ref{thrm-H}, in order to apply Theorem~\ref{thrm-S}, we first show one last structural result, Lemma~\ref{lem-col-S}. In order to prove it, we start with two more lemmas.

\begin{lemma}\label{lem-safety-add-edge-H}
    Let $G$ be a graph and $c,k \geq 1$. Suppose there is some $C \subseteq V(G)$ such that $G[C]$ is a minimal $T$-type subgraph with respect to $c$ and with some witness set $S$ of size~$2$. Let $G'$ be the graph obtained from $G$ by removing the vertices of $C$ and adding an edge between the two vertices of $S$. If $G$ is $S_{1,k,k,k}$-subgraph-free then so is the graph $G'$.
\end{lemma}

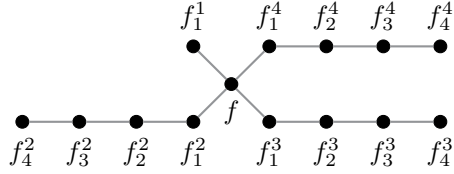
\begin{figure}[ht]
    \centering
    \begin{tikzpicture}
    \tikzstyle{vertex} = [draw, circle, inner sep=1.7pt,fill];
    \tikzstyle{edge} = [draw=Gray, thick];

    \node[vertex, label = below:$f$](x) at (0,0){};
    \node[vertex, label = above:$f_1^1$](f11) at (-0.5,0.5){};

    \node[vertex, label = below:$f_1^2$](f21) at (-0.5,-0.5){};
    \node[vertex, label = below:$f_2^2$](f22) at (-1.25,-0.5){};
    \node[vertex, label = below:$f_3^2$](f23) at (-2,-0.5){};
    \node[vertex, label = below:$f_4^2$](f24) at (-2.75,-0.5){};

    \node[vertex, label = below:$f_1^3$](f31) at (0.5,-0.5){};
    \node[vertex, label = below:$f_2^3$](f32) at (1.25,-0.5){};
    \node[vertex, label = below:$f_3^3$](f33) at (2,-0.5){};
    \node[vertex, label = below:$f_4^3$](f34) at (2.75,-0.5){};

    \node[vertex, label = above:$f_1^4$](f41) at (0.5,0.5){};
    \node[vertex, label = above:$f_2^4$](f42) at (1.25,0.5){};
    \node[vertex, label = above:$f_3^4$](f43) at (2,0.5){};
    \node[vertex, label = above:$f_4^4$](f44) at (2.75,0.5){};

    \draw[edge](x) -- (f11);
    \draw[edge](x) -- (f21);
    \draw[edge](x) -- (f31);
    \draw[edge](x) -- (f41);

    \draw[edge](f21) -- (f22);
    \draw[edge](f22) -- (f23);
    \draw[edge](f23) -- (f24);

    \draw[edge](f31) -- (f32);
    \draw[edge](f32) -- (f33);
    \draw[edge](f33) -- (f34);

    \draw[edge](f41) -- (f42);
    \draw[edge](f42) -- (f43);
    \draw[edge](f43) -- (f44);
\end{tikzpicture}
    \caption{An illustration of $S_{1,4,4,4}$ with the corresponding vertex names.}
    \label{fig:star}
\end{figure}

\begin{proof}
    We will prove the contrapositive. Let $S= \{u,v\}$. Suppose there is some $F \subseteq G'$ such that $F$ is isomorphic to $S_{1,k,k,k}$, then we claim $G$ must also contain some subgraph isomorphic to $S_{1,k,k,k}$. Let $f$ denote the degree~$4$ vertex of $F$. Let $f^1_1$ be that vertex lying on the branch of length $1$ and for each $i \in \{2,3,4\}$ and $j \in \{1, \ldots, k\}$, let $f^i_j$ denote the $j$th vertex along the $i$th branch, see Figure~\ref{fig:star}.
    
    By definition of a $T$-type subgraph, we have that $|N(u) \cap( C \cup \{u,v\})|, |N(v) \cap (C \cup \{u,v\})| \geq 2$. Let $u',u''$ be two neighbours of $u$ in $C \cup \{u,v\}$ and let $v',v''$ be two neighbours of $v$ in $C \cup \{u,v\}$. Note, if $u$ and $v$ are not connected in $G[C \cup \{u,v\}]$ then there is some $C' \subsetneq C$ such that $v, v' \notin C'$ and $G[C']$ is a $T$-type subgraph with witness set $\{u\}$. By assumption $G[C]$ is minimal, that is without loss of generality there is a path from $u$ to $v$ via $u'', v''$ and vertices in $C$. We denote this path by $P^C = (u,u'', \ldots, v'',v)$.

    If $(u,v) \notin E(F)$ then $F$ is a subgraph of $G$, that is, $G$ contains a subgraph isomorphic to $S_{1,k,k,k}$. Let $F'$ be the graph obtained from $F$ by replacing the edge $uv$ by the path $P^C$. Given the vertices of $P^C - \{u,v\}$ do not appear in $G'$, $F'$ is (possibly the supergraph of) some $S_{1,k,k,k}$ subgraph in $G$.
\end{proof}

\begin{lemma}\label{lem-safety-contract-H}
    Let $G$ be a graph and $c, k \geq 1$. Suppose there is some $C \subseteq V(G)$ such that $G[C]$ is a minimal $T$-type subgraph with respect to $c$ and with some witness set $S$ of size~$2$. Let $G'$ be the graph obtained from $G$ by removing the vertices of $C \cup S$ and adding a new vertex $w$ such that $N(w) = N(S)$. If $G$ is $S_{1,k,k,k}$-subgraph-free then so is the graph $G'$.
\end{lemma}

\begin{proof}
    We will prove the contrapositive. Let $S= \{u,v\}$. Suppose there is some $F \subseteq G'$ such that $F$ is isomorphic to $S_{1,k,k,k}$, then we claim $G$ must also contain some subgraph isomorphic to $S_{1,k,k,k}$. Let $f$ denote the degree~$4$ vertex of $F$. Let $f^1_1$ be that vertex lying on the branch of length $1$ and for each $i \in \{2,3,4\}$ and $j \in \{1, \ldots, k\}$, let $f^i_j$ denote the $j$th vertex along the $i$th branch, see Figure~\ref{fig:star}.

    By definition of a $T$-type subgraph, $|N(u) \cap( C \cup \{u,v\})|, |N(v) \cap (C \cup \{u,v\})| \geq 2$ holds. Let $u',u''$ be two neighbours of $u$ in $C \cup \{u,v\}$ and let $v',v''$ be two neighbours of $v$ in $C \cup \{u,v\}$. 
    Note, if $u$ and $v$ are not connected in $G[C \cup \{u,v\}]$ then there is some $C' \subsetneq C$ such that $v, v' \notin C'$ and $G[C']$ is a $T$-type-subgraph with witness set $\{u\}$. 
    By assumption $G[C]$ is minimal, that is, without loss of generality, there is a path from $u$ to $v$ via $u'', v''$ and vertices in $C$. We denote this path by $P^C = (u,u'', \ldots, v'',v)$.

    If $w \notin V(F)$, then $F$ is a subgraph of $G$ and so $G$ contains a subgraph isomorphic to $S_{1,k,k,k}$. If $w \neq f$, then either $w = f^1_1$ or $w = f^i_j$ for some $i \in \{2,3,4\}$ and $j \in \{1, \ldots, k\}$. 
    Recall there exists the path $P^C \subseteq G$ with $V(P^C) \cap V(G') = \varnothing$. 
    Let $F'$ be the graph obtained by replacing the vertex $w$ with the path $P^C$. Now $F'$ is (possibly the supergraph of) some $S_{1,k,k,k}$ subgraph in $G$. Suppose now $w = f$.
    
    It follows that either $u$ or $v$ is adjacent to at least~$2$ of $f^2_1,f^3_1,f^4_1$ in $G$. Without loss of generality, say $v$ is adjacent to $f^3_1$ and $f^4_1$. 
    If $v$ is also adjacent to $f^2_1$ in $G$, then replacing the vertex $f^1_1$ by $v''$ in $F$ we find a $S_{1,k,k,k}$ subgraph in $G$. 
    From this, we assume that $v$ is not adjacent to $f^2_1$ and so $u$ must be adjacent to $f^2_1$. 
    Let $\hat{P}^C$ be a shortest path from~$u$ to~$v$ in $G[V(C) \cup \{u,v\}]$. 
    Note $\hat{P}^C$ cannot contain both $v'$ and $v''$ else there is a shorter path containing only one of them. 
    Without loss of generality, we may assume that $\hat{P}^C$ does not contain $v'$. 
    Let $F'$ be the graph obtained by replacing the edge $wf^2_1$ by the path $\hat{P}^C$ and the edge $wf^1_1$ by the edge $vv'$. Once again $F'$ is (possibly the supergraph of) some $S_{1,k,k,k}$ subgraph in~$G$.
\end{proof}

\noindent
Using the previous two lemmas we can now show the following lemma.

\begin{lemma}\label{lem-col-S}
    Let $(G,r)$ be an instance of \col\ and let $c,k \geq 1$ be constants. We can in polynomial time obtain a graph $G'$ such that:
    \begin{itemize}
        \item $G'$ contains no $T$-type subgraph (with respect to $c$);
        \item $G'$ can be coloured with $r$ colours, if and only if $G$ can;
        \item if $G$ is $S_{1,k,k,k}$-subgraph-free, then so is $G'$.
    \end{itemize}
\end{lemma}

\begin{proof}
    Each of our $T$-type subgraphs will be with respect to the constant $c$. We now define a $T$ removal operation as follows. Suppose $G$ contains some set $C \subseteq V(G)$ such that $G[C]$ is a minimal $T$-type subgraph with witness set $S$. By Lemma~\ref{lem-T-type-smalltd}, $\td(G[C \cup S]) \leq 3k+2$. We branch on the (at most) two non-isomorphic colourings of $S$. For each branch, we consider the precolouring extension problem on the graph $G[C \cup S]$, with the set $S$ of precoloured vertices. Since $\td(G[C \cup S]) \leq 3k+2$, we can use the result of~\cite{JS97} again to process each branch in polynomial time.

    If there is no proper $r$-colouring of $G[C \cup S]$, then $(G,r)$ is a no-instance of \col{} and so we return no. Assume now $G[C \cup S]$ can be coloured using $r$ colours. If $|S| = 1$, then $G$ can be coloured using $r$ colours, if, and only if, both $G[C \cup S]$ and $G-C$ can be coloured using $r$. Let $G^T(C \cup S) = G-C$. As $|N(S) \cap C| \geq 2$, we get that $|V(G^T(C \cup S))| < |V(G)|$. In addition, if $G$ is $S_{1,k,k,k}$-subgraph-free then so is $G^T(C \cup S)$. 

    That is, we now assume that $|S| =2$. Let $S = \{u,v\}$. If there is both a proper colouring where $u,v$ are coloured the same and where they are different, then $G$ can be coloured with $r$ colours if and only if, $G-C$ can be coloured. Let $G^T(C \cup \{u,v\}) = G-C$. As the class of $S_{1,k,k,k}$-subgraph-free graphs is closed under vertex deletion, if $G$ is $S_{1,k,k,k}$-subgraph-free, then so is $G^T(C \cup \{u,v\})$. If there is only a proper colouring of $G[C \cup \{u,v\}]$ where $u,v$ are coloured differently, then let $G^T(C \cup \{u,v\})$ be the graph obtained from $G$ by removing the vertices of $C$ and adding the edge $uv$. It follows from Lemma~\ref{lem-safety-add-edge-H}, that if $G$ is $S_{1,k,k,k}$-subgraph-free then so is $G^T(C \cup S)$. 
    As $|N(S) \cap C| \geq 2$, we get that $|V(G^T(C \cup S))| < |V(G)|$.
    The remaining case is where $u$ and $v$ must take the same colour in every colouring of $G[C \cup S]$. Let $G^T(C \cup S)$ be the graph obtained from $G$ by removing the vertices of $C$ and adding a new vertex $w$ such that $N(w) = N(u) \cup N(v)$. Now, $G$ can be coloured using $r$ colours if, and only if, $G^T(C \cup S)$ can be coloured using $r$ colours. As $|N(S) \cap C| \geq 2$, $|V(G^T(C \cup S))| < |V(G)|$. Further, from Lemma~\ref{lem-safety-contract-H}, if $G$ is $S_{1,k,k,k}$-subgraph-free then so is $G^T(C \cup S)$.

    In polynomial time, for every $S \subseteq V(G)$, such that $|S| \leq 2$, we can decide if $S$ is the witness set for some minimal $T$-type subgraph. If so we can construct a graph $G^T(C \cup S)$ such that $|G^T(C \cup S)| < |V(G)|$ and $(G^T(C \cup S),r)$ is a yes-instance of \col{} if, and only if, $(G,r)$ is. Applying this operation exhaustively, we obtain a graph $G'$ such that $G'$ contains no $T$-type subgraph, $G'$ can be coloured with $r$ colours, if and only if $G$ can and if $G$ is $S_{1,k,k,k}$-subgraph-free, then so is $G'$.
\end{proof}

\noindent
Given an instance $(G,r)$ of \col, we now apply Lemma~\ref{lem-col-S} and obtain a modified instance $(G',r)$. We may assume by \cite{JMPPSV23} that $G'$ contains no bridges and further by~\cite{Brooks41} that $G'$ contains some degree~$4$ vertex.
By Theorem~\ref{thrm-S}, the treedepth of $G'$ is bounded by a constant and since \col\ is polynomial-time solvable for graphs of bounded treewidth~\cite{KR03} and therefore also for graphs of bounded treedepth, we obtain the following.

\begin{theorem}\label{thm-col-S}
    \col\ is solvable in polynomial time for $S_{1,k,k,k}$-subgraph-free graphs, for all $k\geq 1$.
\end{theorem}

\begin{proof}
    Let $(G,r)$ be an instance of \col, where $G$ is a $S_{1,k,k,k}$-subgraph-free graph.
    We apply Lemma~\ref{lem-col-S} and obtain in polynomial time an instance $(G',r)$ of \col\ such that $G'$ contains no $T$-type subgraph with treedepth bound $16(2k-1)(k-1)$, and $G'$ is $S_{1,k,k,k}$-subgraph-free. Further, $G'$ can be coloured using $r$ colours if, and only if, $G$ can be coloured using $r$ colours.
    In~\cite{JMPPSV23}, it was shown that the graph obtained from $G'$ by deleting all bridges can be coloured using $r$ if, and only if $G'$ can be coloured using $r$ colours. That is, we assume that $G'$ contains no bridges.
    If $G'$ has maximum degree at most~$3$, by Brooks' Theorem~\cite{Brooks41}, we solve \col\ in polynomial time. That is we assume that $G'$ contains some degree~$4$ vertex.
    Hence, from Theorem~\ref{thrm-S}, as $G'$ is $S_{1,k,k,k}$-subgraph-free, $\td(G) < 8(7k^3+15k^2-\frac{4k}{9}+3)^2+6$. As \col\ is polynomial-time solvable for graphs of bounded treewidth~\cite{KR03} and therefore also for graphs of bounded treedepth the theorem follows.
\end{proof}

\section{Wider Applicability of Our Techniques}\label{s-wider}

Our main structural results, Theorems~\ref{thrm-H} and~\ref{thrm-S}, can be applied to other problems, as we explain below.

\medskip
\noindent
{\it Set of Conditions I.}
Let $\Pi$ be a problem that is polynomial-time solvable for graphs of bounded 
treedepth. If, for every $k,m\geq 1$, there exist constants $c$, $\ell$ such that we are able to preprocess the input graph in polynomial time to a polynomial number of smaller graphs of minimum degree at least~$3$ that contain neither a large protected fan nor an $L$-type subgraph with length bound $\ell$ and treedepth bound $c$, then $\Pi$ is polynomial-time solvable on $\mathbb{H}_m^{1,k,k,k}$-subgraph-free graphs for all $k,m\geq 1$.

\medskip
\noindent
{\it Set of Conditions II.}
Let $\Pi$ be a problem that is polynomial-time solvable for both subcubic graphs and graphs of bounded treedepth. If, for every $k\geq 1$, there exists a constant $c$ such that we are able to preprocess the input graph in polynomial time to a polynomial number of smaller graphs that each contain some vertex of degree at least~$4$ but neither contain a proper bridge nor a $T$-type subgraph with treedepth bound $c$, then
$\Pi$ is polynomial-time solvable for $S_{1,k,k,k}$-subgraph-free graphs for all $k\geq 1$.

\smallskip
\noindent
We note that if a problem satisfies the above two sets of conditions then the algorithm of~\cite{JMPPSV23} can be applied directly
to obtain an algorithm for $S_{1,1,r,r}$-subgraph-free graphs.
Conversely, the above two sets of conditions are not satisfied by every problem
considered in~\cite{JMPPSV23}: {\sc Matching Cut} is \NP-complete
for bipartite graphs in which one partition class has maximum degree~$2$~\cite{Mo89} and thus on $(\mathbb{H}_1,\mathbb{H}_3,\mathbb{H}_5,\ldots)$-subgraph-free graphs.
It is still open whether the polynomial-time result for $S_{1,1,r,r}$-subgraph-free graphs for {\sc Matching Cut} from~\cite{JMPPSV23} can be generalized to 
$S_{1,k,k,k}$-subgraph-free graphs. 
This is in contrast to a related problem: {\sc Stable Cut}, which we recall is the problem of deciding if a connected graph $G$ has an independent set $I$ such that $G-I$ is disconnected.  Below we apply our techniques to show the same results for {\sc Stable Cut} as we did for {\sc Colouring}:

\begin{theorem}\label{thm-sc-S}
{\sc Stable Cut} is solvable in polynomial time for $\mathbb{H}_m^{1,k,k,k}$-subgraph-free graphs for all $k,m\geq 1$, and for $S_{1,k,k,k}$-subgraph-free graphs for all $k\geq 1$.
\end{theorem}

\begin{proof}
    Suppose first that $G$ contains some cutset $S$ of size at most~$2$. Note that we can find such a cutset in polynomial time by iterating over all possible sets. If either $S$ has size $1$ or $G[S]$ consists of a pair of non-adjacent vertices, then $G$ contains a stable cut. Suppose now $G[S]$ consists of two adjacent vertices. As $S$ is a clique, it was shown by Le, Mosca and M{\"{u}}ller~\cite{LeMM08} that $G$ has a stable cutset if and only if there is some connected component $G[C]$ of $G - S$ such that $G[C \cup S]$ has a stable cutset. We therefore consider each such connected component of $G-S$ in turn. 
    
    That is, we now assume that $G$ does not contain any cutset with size at most~$2$. Note this implies that $G$ has minimum degree at least~$3$ as the neighbourhood of any vertex with degree at most~$2$ is a cutset with size at most~$2$. Further, we may assume $G$ does not contain a bridge, a $T$-type subgraph or an $L$-type subgraph (for any treedepth and length bounds), since again this would give a cutset of size at most~$2$.

    Suppose $G$ is $S_{1,k,k,k}$-subgraph-free. As \stcut\ is polynomial-time solvable for graphs with maximum degree at most~$3$ \cite{Chvatal84}, we assume that $G$ contains some vertex with degree at least~$4$. As $G$ is bridgeless, contains no $T$-type subgraph and is $S_{1,k,k,k}$-subgraph-free, by Theorem~\ref{thrm-S}, $G$ has treedepth less than $2(7k^3+15k^2-\frac{4k}{9} + 3)^2 + 6$. As \stcut\ is readily seen to be definable in monadic
second-order logic, it is
    polynomial-time solvable on graphs of bounded treewidth due to Courcelle's Theorem~\cite{Co90}, and therefore also for graphs of bounded treedepth. Hence, we can decide if $G$ contains a stable cut in polynomial time.

    Suppose now $G$ is $\mathbb{H}_m^{1,k,k,k}$-subgraph-free. Recall that $G$ has minimum degree at least~$3$ and no $L$-type subgraph. Suppose $G$ contains some protected fan with vertex set $F$, centre $z$ and ends $x$ and $y$. For every vertex $v \in F \setminus \{x,y,z\}$, $N(v) \subseteq F$ and so $v$ is not contained in any minimal stable cut of $G$. It follows that $v$ can be deleted from $G$. We note that the class of $\mathbb{H}_m^{1,k,k,k}$-subgraph-free graphs is closed under vertex deletion. We exhaustively delete such vertices to obtain a graph $G'$ such that $G'$ contains no protected fan (of size at least~$3$). It now follows from Theorem~\ref{thrm-H} that $G'$ has treedepth bounded by some constant. Hence, we can apply Courcelle's Theorem~\cite{Co90} again to decide in polynomial time if $G'$ contains a stable cut. This concludes the proof of the theorem.
\end{proof}

\noindent
The {\sc Feedback Vertex Set} problem is to decide if a graph has a set $S$ ({\it feedback vertex set}) of size at most~$k$ for some given integer~$k$ such that $G-S$ contains no cycles. The {\sc Feedback Vertex Set} problem does not satisfy Set of Conditions~I, as we may not assume an instance of it has minimum degree at least~$3$. However, we can show that {\sc Feedback Vertex Set} satisfies Set of Conditions~II.

\begin{lemma}\label{lem-fvs}
   Let $(G,r)$ be an instance of \fvs{} and let $c,k \geq 1$ be constants. We can in polynomial time obtain a graph $G'$ and a constant $r' \geq 1$ such that:
    \begin{itemize}
        \item $G'$ contains no $T$-type subgraph (with respect to $c$);
        \item $(G',r')$ is a yes-instance of \fvs{} if and only if $(G,r)$ is a yes-instance of \fvs;
        \item if $G$ is $S_{1,k,k,k}$-subgraph-free, then so is $G'$.
    \end{itemize}
\end{lemma}
\begin{proof}
   Throughout this proof, all $T$-type subgraphs are with respect to the constant $c$.
    Let $(G,r)$ be an instance of \fvs. Suppose that $G$ contains some set $C \subseteq V(G)$ such that $G[C]$ is a minimal $T$-type subgraph with witness set $S$. Recall that $S$ has size at most~$2$.
    We introduce a  
    \emph{$T$-removal operation} resulting in a modified instance $(G^T(C \cup S), \hat{r})$ such that $(G,r)$ is a yes-instance of \fvs{} if and only if $(G^T(C \cup S), \hat{r})$ is a yes-instance of \fvs.

 As a first step, we find a minimum feedback vertex set~$A$ of $G[C \cup S]$. Note that this can be done in polynomial time as follows. We consider every $S' \subseteq S$ and solve {\sc Feedback Vertex Set} on 
 $G[C\cup (S\setminus S')]$ and afterwards we add $S'$ to the feedback vertex set that we found. As $|S|\leq 2$, and moreover, {\sc Feedback Vertex Set} is polynomial-time solvable for graphs of bounded treewidth~\cite{ALS91} and $\td(G[C \cup S]) \leq 3c+2$, this takes polynomial time.
  
  Let $A$ be a minimum feedback vertex set of $G[C \cup S]$. We distinguish two cases based on the intersection of $A$ and $S$.

Suppose first that $S \cap A  \neq \varnothing$. Since $|S| \leq 2$, this implies that $|S\setminus A| \leq 1$. Further, recall that by definition, $S$ is a cutset separating $C$ from $G-(C \cup S)$. Thus, no vertex of $C$ is contained in a cycle in $G-A$.
Let $\hat{A}$ be a minimum feedback vertex set of $G$. Since $A$ is a minimum feedback vertex set of $G[C \cup S]$, we may assume that $A \subseteq \hat{A}$.
This implies that $G$ contains a feedback vertex set of size $r$, if, and only if, $G - (C \cup S)$ contains a feedback vertex set of size $r - |A|$. Let $G^T(C \cup S) = G - (C \cup S)$ and $\hat{r} = r - |A|$.
    
   Consider now the other case, that is when $S \cap A = \varnothing$. If no vertex in $C$ is contained in a cycle in $G-A$, we proceed as in the case before. So we may assume that some vertex in $C$ is contained in some cycle in $G-A$. This implies that $|S| > 1$, so $|S| = 2$. Let $S = \{u,v\}$, for $u,v \in V(G)$. 
    It follows that both $u$ and $v$ are contained in this cycle. If there does not exist a path from $u$ to $v$ in $G[C \cup S]$, then $G$ contains a feedback vertex set of size $r$, if, and only if, $G - (C \cup S)$ contains a feedback vertex set of size $r - |A|$. We set $G^T = G-(C\cup S)$ and $\hat{r} = r - |A|$.
    Otherwise, if there is a path from $u$ to $v$ in $G[C \cup S]$, let $G^T(C \cup \{u,v\})$ denote the graph obtained from $G$ by removing the vertices of $C$ and adding the edge $(u,v)$ and $\hat{r} = r - |A|$. Now $G$ contains a feedback vertex set of size $r$, if, and only if, $G^T(C \cup S)$ contains a feedback vertex set of size $\hat{r}$.
    
    This operation therefore allows us to remove $C \cup S$. Note that if $G$ is $S_{1,k,k,k}$-subgraph-free then also the resulting graph $G^T(C \cup S)$ is $S_{1,k,k,k}$-subgraph-free.  
    This follows from the fact that $S_{1,k,k,k}$-subgraph-freeness is preserved under vertex deletion and from Lemma~\ref{lem-safety-add-edge-H} in the last case.

    \medskip
    \noindent
    In order to obtain $G'$ as given in the lemma statement, we apply the following exhaustively.    
    In polynomial time, for every $S \subseteq V(G)$, such that $|S| \leq 2$, decide if $S$ is the witness set for some minimal $T$-type subgraph. If so we apply the $T$ removal operation as described above and obtain a graph $G^T(C \cup S)$ such that $|G^T(C \cup S)| < |V(G)|$, and $(G,r)$ is a yes-instance of \fvs{} if and only if $(G^T(C \cup S), \hat{r})$ is a yes-instance of \fvs{} and if $G$ is $S_{1,k,k,k}$-subgraph-free, then so is $G^T(C\cup S)$. 
    
    Applying this operation exhaustively, we obtain a graph $G'$ such that $G'$ contains no $T$-type subgraph, $G'$ contains a feedback vertex set of size $r'$, if and only if $G$ contains a feedback vertex set of size $r$, and if $G$ is $S_{1,k,k,k}$-subgraph-free, then so is $G'$.
    Note that each application of the $T$ removal operation can be done in polynomial time. Further, since for every such application we get that $|V(G^T(C \cup S))| \leq |V(G)|$, we apply the $T$ removal operation at most $O(n)$ times. Thus, $G'$ can be obtained in polynomial time.
\end{proof}

\noindent
We are now ready to prove Theorem~\ref{thm-FVS}. 

\begin{theorem}\label{thm-FVS}
    \fvs{} is solvable in polynomial time for $S_{1,k,k,k}$-subgraph-free graphs.
\end{theorem}
\begin{proof}
  Let $(G,r)$ be an instance of \fvs, where $G$ is a $S_{1,k,k,k}$-subgraph-free graph. Let $c = 16(2k-1)(k-1)$ be a constant.
    We apply Lemma~\ref{lem-fvs} and obtain in polynomial time an instance $(G',r')$ of \fvs{} such that $G'$ is $S_{1,k,k,k}$-subgraph-free, contains no type-$T$ subgraph with respect to $c$ and $G'$ is $S_{1,k,k,k}$-subgraph-free. Further, $(G,r)$ is a yes-instance of \fvs{} if and only if $(G',r')$ is a yes-instance of \fvs.

    It is known that \fvs{} can be solved in polynomial time for cubic graphs~\cite{UKG88} and for graphs of bounded treewidth~\cite{ALS91}, and thus for graphs of bounded treedepth.
    As a bridge is not used in a cycle, we may safely remove it (without destroying the $S_{1,k,k,k}$-subgraph-freeness).
    Hence, we may assume that $G'$ is   
    bridgeless, and thus proper bridgeless, and also that $G'$ contains a vertex of degree~$4$, has treedepth at least~$8(7k^3 + 15k^2 - \frac{4k}{9} + 3)^2 + 6$ and no $T$-type subgraph with respect to $c$.
    We get from Theorem~\ref{thrm-S} that $G'$ contains $S_{1,k,k,k}$ as a subgraph, a contradiction to the fact that $G'$ is $S_{1,k,k,k}$-subgraph-free. So we conclude that we are always in one of the cases above and the theorem follows.
\end{proof}

\section{Conclusions}\label{sec-conc}

In contrast to $H$-free graphs~\cite{KKTW01} and $H$-minor-free graphs (see e.g.~\cite{JMOPPSV25}), for which {\sc Colouring} is fully classified, even formulating a dichotomy for $H$-subgraph-free graphs is still challenging. However, our new dichotomies for {\sc Colouring}, {\sc Stable Cut} and {\sc Feedback Vertex Set}, being the first of their kind, open up the way for further progress for many graph problems. In this section we will summarize our findings and explore some new directions.
%JL:but where did we say anything about hardness results for SC, FVC?
%DP*: we dont

\begin{figure}[t]
    \centering
   \begin{tikzpicture}
    \tikzstyle{vertex} = [draw, circle, inner sep=1.7pt,fill];
    \tikzstyle{edge} = [draw=Gray, thick];

\begin{scope}[scale = 0.7]
\begin{scope}[shift = {(0,0)}]
    \node[vertex](v1) at (0,0){};
    \node[vertex](v2) at (1,0){};
    \node[vertex](v3) at (2,0){};
    \node[vertex](v4) at (3,0){};
    \node[vertex](v5) at (0,1){};
    \node[vertex](v6) at (0,-1){};
    \node[vertex](v7) at (2,1){};
    \node[vertex](v8) at (2,-1){};

    \draw[edge](v1) -- (v2);
    \draw[edge](v1) -- (v5);
    \draw[edge](v1) -- (v6);
    \draw[edge](v2) -- (v3);
    \draw[edge](v3) -- (v4);
    \draw[edge](v3) -- (v7);
    \draw[edge](v3) -- (v8);

\end{scope}

\begin{scope}[shift = {(4.5,0)}]
    \node[vertex](v1) at (0,0){};
    \node[vertex](v2) at (1,0){};
    \node[vertex](v3) at (2,0){};
    \node[vertex](v4) at (3,0){};
    \node[vertex](v5) at (4,0){};
    \node[vertex](v6) at (1,1){};
    \node[vertex](v7) at (1,-1){};
    \node[vertex](v8) at (2,1){};

    \draw[edge](v1) -- (v2);
    \draw[edge](v2) -- (v3);
    \draw[edge](v3) -- (v4);
    \draw[edge](v4) -- (v5);
    \draw[edge](v3) -- (v8);
    \draw[edge](v2) -- (v7);
    \draw[edge](v2) -- (v6);

\end{scope}

\begin{scope}[shift = {(10,0)}]
    \node[vertex](v1) at (0,0){};
    \node[vertex](v2) at (1,0){};
    \node[vertex](v3) at (2,0){};
    \node[vertex](v4) at (3,0){};
    \node[vertex](v5) at (4,0){};
    \node[vertex](v6) at (1,1){};
    \node[vertex](v7) at (2,1){};
    \node[vertex](v8) at (2,-1){};

    \draw[edge](v1) -- (v2);
    \draw[edge](v2) -- (v3);
    \draw[edge](v3) -- (v4);
    \draw[edge](v4) -- (v5);
    \draw[edge](v3) -- (v8);
    \draw[edge](v3) -- (v7);
    \draw[edge](v2) -- (v6);

\end{scope}

\begin{scope}[shift = {(15.5,0)}]
    \node[vertex](v1) at (0,0){};
    \node[vertex](v2) at (1,0){};
    \node[vertex](v3) at (2,0){};
    \node[vertex](v4) at (0,1){};
    \node[vertex](v5) at (0,-1){};
    \node[vertex](v6) at (1,1){};
    \node[vertex](v7) at (2,1){};
    \node[vertex](v8) at (2,-1){};

    \draw[edge](v1) -- (v2);
    \draw[edge](v1) -- (v4);
    \draw[edge](v1) -- (v5);
    \draw[edge](v2) -- (v3);
    \draw[edge](v2) -- (v6);
    \draw[edge](v3) -- (v7);
    \draw[edge](v3) -- (v8);

\end{scope}
\end{scope}
    
\end{tikzpicture}
    \caption{The four open cases of graphs $H$ on eight vertices.}
    \label{fig:opencases}
\end{figure}
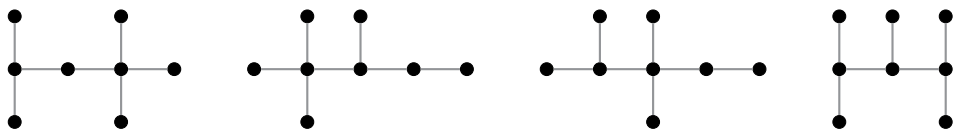

We proved that {\sc Colouring} on $H$-subgraph-free graphs is polynomial-time solvable on 
$\mathbb{H}_m^{1,k,k,k}$-subgraph-free graphs for all $k,m\geq 1$ and on
$S_{1,k,k,k}$-subgraph-free graphs for all $k\geq 1$. Combining these results with known \NP-completeness results yields 
a complete complexity classification of {\sc Colouring} on $H$-subgraph-free graphs whenever $H$ is a subdivided~$\mathbb{H}_0$ or a subdivided~$\mathbb{H}_1$. As mentioned in Section~\ref{sec:intro}, combining our new results with known results leaves open the cases where\\[-10pt]
\begin{itemize}
\item [(i)] $H$ is a tree of maximum degree~$4$ with exactly one vertex of degree~$4$ and at least one vertex of degree~$3$; or
\item [(ii)] $H$ is a subcubic tree with at least three vertices of degree~$3$.\\[-10pt] 
\end{itemize}
\noindent
Some graphs $H$ of type~(i) and~(ii) are covered by existing \NP-completeness results for smaller graphs~$H'$ (e.g. when $H$ contains an $S_{2,2,2,2}$). Nevertheless, there still exist exactly four open cases of graphs $H$ on eight vertices, three of type~(i) and one of type~(ii); see Figure~\ref{fig:opencases}.
Our structural analysis does not apply to these four open cases and new ideas are needed. In fact, even the case when $H$ is the $7$-vertex graph obtained from $\mathbb{H}_1$ by adding a pendant vertex to one of the degree-$3$ vertices is not covered by our new technique.
Indeed, there exists a family of $H$-subgraph-free graphs demonstrating that a modified Theorem~\ref{thrm-H}, where $\mathbb{H}_1$ is replaced by $H$ does not hold. Likewise, a modified of Theorem~\ref{thrm-S} with $S_{1,k,k,k}$ replaced by $H$ also does not hold.
This case, which is of type~(i), is shown to be polynomial-time solvable in~\cite{GPR15} by a rather involved tailor-made algorithm. Studying the proof technique in more detail would be a natural starting point.

Finally, we also propose to investigate in a systematic way which other graph problems satisfy Sets of Conditions I and II from Section~\ref{s-wider}. Moreover, we ask if {\sc Matching Cut} is polynomial-time solvable for $S_{1,k,k,k}$-subgraph-free graphs for all $k\geq 1$ 
and if {\sc Feedback Vertex Set} is polynomial-time solvable for $\mathbb{H}_m^{1,k,k,k}$-subgraph-free graphs for all $k,m\geq 1$.

\bibliography{ref}

\end{document}